# Admissible Permutations and the HCP, the AP and the TSP

Howard Kleiman

Prof. Emeritus, QCC, City University of New York

**Copyright 2002**

# Table of Contents

# Preface

My purpose in writing this book is to more widely disseminate the applications of H-admissible and admissible permutations to obtain Hamilton circuits in graphs and directed graphs. In particular, Conjecture 1.1 hypothesizes that Algorithm G or Algorithm $G_{no\ r-vertices}$ always obtains a hamilton circuit in a an arbitrary graph containing one in polynomial time. Wider research on whether a counter-example exists would be useful. I also would like to know how often the algorithm given in chapter 3, section 3.2 gives an optimal solution to the traveling salesman problem in polynomial time. Secondly, how often does it give an approximate solution to the TSP that is very close to an optimal one.

# Introduction

Chapter 1.

Let $C$ be a circle. Assume that n equally-spaced points have been placed on the circle at the points $\frac{2\pi i}{n}$ (i = 1, 2, ... , n) with $2\pi$ at 12 o'clock. Now place the number i on the circle at the point $\frac{2\pi i}{n}$. Assume that we have represented the vertices 1, 2, ... , n of a graph or digraph , $G$, on the circle as the consecutive points of an *n*-cycle, $H = (1\ 2\ 3\ ...\ n)$, where the edges (arcs) of $H$ don't necessarily lie in $G$. Call a 3-cycle (*a b c*) *H*-admissible if $H(a\ b\ c) = H'$ where $H'$ is also an n-cycle. In theorem 1.1, we prove that (*a b c*) is *H*-admissible if and only if the chords (*a, H(b)*) and (*b, H(c)*) properly intersect in *C*. If *G* has edges that have been randomly chosen, then the probability, $P_3$, that *(a b c)* is *H*-admissible is $\frac{n-3}{2(n-2)}$. In particular, as n $\to \infty$, then $P_3 \to \frac{1}{2}$. Now consider the permutation (*a c*))(*b d*)) where (a, *H(c)*), (c, *H(a)*),(b, *H(d)*), (d, *H(b)*) are all proper, distinct chords in *C*. *(a c)(b d)* is *H*-admissible if $H(a\ c)(b\ d) = H''$ where $H''$ is also an *n*-cycle. An *H*-admissible product of two disjoint cycles is called a *POTDTC*. Theorem 1.2 proves that necessary and sufficient conditions for *(a c)(b d)* to be *H*-admissible are that *a, b, c* and *d* are distinct points such that the chords *(a, c)* and *(b, d)* properly intersect in the interior of *C*. Furthermore, if the edges of the graph *G* are randomly chosen, then the probability, $P_{2,2}$, that *(a c)(b d)* will be *H*-admissible is $\frac{n-3}{3(n-2)}$. Thus, as $n \to \infty$, $P_{2,2} \to \frac{1}{3}$. A *pseudo-hamilton circuit* is an *n*-cycle not all of whose edges belong to *G*. A vertex, *v*, in *G* such that (*v, H(v)*) doesn't belong to *G* is a *pseudo-arc vertex*. A *pseudo-arc* is an arc of *H* whose initial vertex is a pseudo-arc vertex. Let *G* be an arbitrary graph. In the above discussion, we assumed that *H* is a pseudo-hamilton circuit*;* *a*, a pseudo-arc vertex of *(a b c)*; *a* and *b* are both pseudo-arc vertices of *(a c)(b d)*. The *in-degree*, $d^-(v)$, of a vertex *v* of a directed graph *D* is the number of arcs terminating in *v*. The *out-degree*, $d^+(v)$, of *v* is the number of vertices emanating from *v*. The *total degree* of *v* equals $d^-(v) + d^+(v)$. If *G* is a graph, then the *total degree* of *v* is $d(v)$. The minimum degree, $\delta(G)$, of a graph or digraph is min { total degree of *v*}. Let *G* contain the path $P = [a, b, c, ... , r]$ such that all of the "interior" vertices $b, c, ... , q$ have total degree 2. Suppose that we replace *P* by the vertex, $v_{abc...qr}$, such that we've deleted the edge $[a, r]$ (if it exists). Assume that we've done this for every path similar to *P*. Then the new graph we've obtained is the contracted graph of G'. In particular, if *G* contains a hamilton circuit, then $\delta(G') \geq 3$. If *D* is a digraph, then both the minimum in-degree and the minimum out-degree is always greater than 0 in the contracted digraph, D'. Let *G* be a graph in which all edges have been randomly chosen. We prove in theorem 1.6 that if a vertex, *v*, of G' is of degree 3 and $(v, H(v))$ is not an edge of *H*, then the minimal probability, $p_3(2)$, that there exists at least two *H*-admissible 3-cycles containing *v* is

$$\frac{286n^6 - 4326n^5 + 23489n^4 - 80546n^3 + 190342n^2 - 112242n - 27624}{360n^6 - 5040n^5 + 29160n^4 - 89280n^3 + 15264n^2 - 138240n + 51840}$$

As $n \to \infty$, $p_3(2) \to \frac{143}{180}$. Corresponding theorems are proved for random digraphs (directed graphs). In Section 1.3, we give preliminaries such as theorems and definitions useful in constructing the algorithms, Algorithm G (for graphs) and Algorithm D (for digraphs) in Section 1.4. Before discussing 1.3, we give the following definitions:

(i) Let *H* be a pseudo-hamilton circuit of a graph, *G*, while *s* is an *H*-admissible 3-cycle or POTDTC. Suppose that $Hs = H_1$. Given that we have a table giving the ordinal values of *H* beginning with $ORD(1) = 1$, we can represent $H_1$ by an *abbreviation* consisting of no more than three numbers. In general, we use abbreviations until we reach one containing $n^{.5}$ vertices. We then explicitly write out $H_i$ and construct new ordinal numbers, $ORD(i)$, as well as their inverses, $ORD^{-1}(i)$ (i). Thus, given a vertex, we can obtain its ordinal value to use in testing for $H_i$-admissibility. On the other hand, we can obtain the unique vertex that has $ORD(i)$ as its ordinal value.

(ii) Suppose that $S = [a, b, c,..., v, w]$ is a subpath of $H$ where $(a, v)$ is an edge of $G$. Then a *rotation* with respect to a and v is a path, $R = [a, v, u, ..., c, b, w]$ that replaces $S$ in $H$.

In [9], Erdös and Renyi give conditions under which a random graph is 2-connected. In [21], Palásti gives conditions for a directed graph to be strongly-connected. In [18], Kömlos and Szemerédi prove that if a random graph, $G$, is constructed by randomly choosing edges from the complete graph, $K_n$, until every vertex has a minimum degree of at least two, then $G$ contains a hamilton circuit with probability approaching 1 as $n \to \infty$. A similar theorem was proven by Bollobás in [4]. In [13], Frieze proved an analogous theorem for random directed graphs. The random graph constructed by Bollobás in [4], we call a *Boll graph*; that constructed by Frieze, a *Frieze-Boll digraph*. Let $D_{\text{2-in, 2-out}}$ be a random directed graph constructed by randomly choosing two arcs out of each if n vertices and two arcs into each vertex. Define the *random regular i-outgraph* $R_i$ in the following way: (i) Randomly select $i$ arcs out of each of $n$ vertices. (ii) Change each of the arcs chosen to an edge. The following conjectures are due to Frieze:

(1) $D_{\text{2-in, 2-out}}$ almost always contains a hamilton circuit as $n \to \infty$.

(2) $R_3$ almost always contains a hamilton circuit as $n \to \infty$.

In [12] Fenner and Frieze introduced $D_{\text{k-in, k-out}}$ and proved that it is strongly-connected is strongly-connected for $k \geq 2$. The hypothesis of (2) was proved valid for $k \geq 5$ by Frieze and Luczak in [14]. Algorithms G and D depend upon the ability to construct *H*-admissible 3-cycles and pairs of disjoint 2-cycles (*POTDTC*) that pass through all of the vertices of $G$ or $D$ provided that $G$ is at least 2-connected or $D$ is strongly-connected. We have previously obtained the following probabilities:

(i) The probability that a randomly chosen pair of arcs of the form $\{(a, H(b)), (b, H(c))\}$ intersect in the circle $C$ approaches $\frac{1}{2}$ as $n \to \infty$.

(ii) The probability that a randomly chosen pair of arcs *{(a, H(b)), (c, H(d))}* intersect in the circle *C* approaches $\frac{1}{3}$ as $n \to \infty$.

(iii) If G has the property that $\delta(G) \geq 3$, and *m* is an arbitrary pseudo-arc vertex, then as $n \to \infty$, the minimal probability that there will be at least two pairs of intersections of arcs from among those pairs of form *{(m, k), (k-1, r)}, {(m, l), (l-1, s)}, {(m, k), (i, 1)},*

*{(m, k), (j, 1)}, {(m, l), (i, 1)}, (m, l), (j, 1)}* is $\frac{143}{180}$. Here $H = (1, 2, ..., m)$ where *m* is a pseudo-arc vertex of *H*. In Section 1.3, we formally present Algorithm G for random graphs *G* and Algorithm D for random digraphs *D*. We first contract *G* to eliminate vertices of degree two or vertices which have in-degree or out-degree less than two in digraph *D*. Assume that the contracted graphs obtained are *G'* and *D'*, respectively. We use a balanced, binary search tree to represent *G'* with each row a branch. For *D'*, we require two balanced, binary search trees. The first, *ROWS*, represents arcs emanating from vertices. The other, *COLUMNS*, represents arcs terminating in vertices. We next obtain a pseudo-hamilton circuit, $H_0$, in graph *G'* or *D'*. Let *a* be an arbitrary pseudo-arc vertex. We pick up to *[log n]* + 1 arcs of form *(a, $H_0$(b))*. Given each arc chosen, we construct up to *[log n]* + 1 arcs *(b, $H_0$(c))*. We check each pair *{(a, $H_0$(b)), (b, $H_0$(c))}* to see if the two arcs intersect and thus define an $H_0$-admissible 3-cycle. We then select up to *{[log n] + 1}* arcs of the form *(c, $H_0$(a))* and check each pair of the form *{(a, $H_0$(b)), (c, $H_0$(a))}* to see if the two arcs intersect. If we are unable to obtain any pair of intersecting arcs, we then test for intersection of up to *[(log n)$^2$]* pairs of randomly chosen edges incident to pseudo-arc vertices *a* and d (assuming that we still have more than one pseudo-arc vertex). We then give each $H_0$-admissible set of vertices $S_i$ = *{(a, $H_0$(b)), (b, $H_0$(c)), (c, $H_0$(a))}* (or *{(a, $H_0$(b)), (b, $H_0$(a), (d, $H_0$(e)), (e, $H_0$(d))})* a value called its *score*. Here

$$SCORE(S_i) = \text{number of arcs in } S_i - \text{number of arc vertices.}$$

We choose a set, $S_i$, from among those sets having the highest score, which has a vertex of maximum degree. Call this set $S_M$. If more than one set satisfies these conditions, we randomly pick one of them. Let $s_0$ be the permutation defined by $S_M$. Then $H_1 = H_0 s_0$. If $G$ is a graph, and we are unable to obtain an $H_0$-admissible permutation, we randomly choose an arc *(a b)* out of a such that the degree of *b* in *G'* is the greatest possible; we then use *(a b)* to construct a rotation. If the terminal point of the rotation is a pseudo-arc vertex, we continue the algorithm from there. Otherwise, we randomly choose a pseudo-arc vertex of greatest degree from *PSEUDO*, the set of pseudo-arc vertices that we update following each iteration. If the graph is a digraph, we choose arcs emanating from *a* or terminating in $H_0(a)$. Since we cannot construct a rotation out of *a*, we backtrack in Algorithm D. If we have progressed to $H_i$, $H_{i+1} = H_i s_{i-1}^{-1} = H_{i-1}$. We use abbreviations to reduce the running time of the algorithm. In Section 1.4, we discuss the probability of success in obtaining a hamilton circuit in an extremal graph or digraph as well as in conjectures (1) and (2) of Frieze. In all cases, the probability of success approaches 1 as $n \to \infty$. In Section 1.6, we prove that if a graph or digraph contains a hamilton circuit, $H_C$, then there always exists a sequence of $H_i$-admissible 3-cycles or *POTDTC* which lead in a finite number of steps to $H_C$. We also prove that if a finite random graph or random digraph contains a hamilton circuit (hamilton cycle), then Algorithm G, Algorithm $G_{no\ r-vertices}$, or Algorithm D, respectively, almost always obtains a hamilton circuit (hamilton cycle) in at most $O(n^{3.5} (\log n)^4)$ running time. The striking thing about the probabilities obtained here is that the probability of failure decreases exponentially as the running time increases polynomially. Thus, if *M* is the number of graphs of degree *n* containing hamilton circuits, the expected number of failures in *M* is less than 1 when the running time is greater than than a calculable number *N(M)*. This leads to the following conjecture:

**Conjecture 1.1**. *Let G be an arbitrary graph containing a hamilton circuit. Then Algorithm G or Algorithm* $G_{no\ r-vertices}$ *always obtains a hamilton circuit of G in polynomial time.*

In Section 1.8, we give a heuristic algorithm for obtaining an approximate solution to the traveling salesman problem. Lastly, we give examples of the algorithms given in chapter 1,

## Chapter 2.

In [13], Frieze gives a procedure for obtaining a hamilton circuit in an extremal digraph, *D*. His algorithm has a probability approaching 1 as $n \to \infty$. In the first step, he constructs a vertex set of about *log n* disjoint cycles that span *D*. In the second step, he patches the cycles together using 2-edge exchanges. In step 3, he patches the cycles, one by one, into the largest cycle by a complex process of double rotations. The respective running times of the phases are $O(n^{1.5})$, $O(n(\log n))$, $O(n^{\frac{4}{3}+o(1)})$.

Before going on, a *derangement* is a permutation of the points in *V = {1,2,...,n}* in which none of the points are fixed. Let *d* be a derangement and s a permutation such that *ds = d'* where *d'* is also a derangement. Then *s* is an *admissible permutation*. Let *H* be a pseudo-hamilton circuit in *D*, a random digraph containing n vertices such that the arcs are chosen randomly until each vertex has both in-degree and out-degree greater than 0. In Section 2, we give an algorithm in which, as $n \to \infty$, a set of derangements sequentially applied to *H* yield a set of approximately *log n* disjoint cycles spanning *D*. The running time of this algorithm is $O(n(\log n)^3)$. Thus, the running time of Frieze's algorithm is reduced to $O(n^{\frac{4}{3}+o(1)})$.

## Chapter 3.

The first part of this chapter obtains a solution to the Assignment Problem, say $\sigma_{APOPT}$, using admissible permutations and a variant of the Floyd-Warshall algorithm given by Papadimitrou and Steiglitz in [6]. The second part uses another variant of F-W to obtain non-negatively-valued cycles, $s_i$ *(i = 1,2, ..., r)* such that $\sigma_{APOPT} s_1 s_2 ... s_r = \sigma_{FWTSPOPT}$, an approximation to an optimal solution to the Traveling Salesman Problem. We also give a sufficient condition for the approximation to be an optimal solution. Let *M* be an *n X n* cost matrix. We sort each row of *M* in increasing order of entries. We then place the column in which the entry appears in *M* in place of its entry. The matrix obtained in this manner is called *MIN(M)*. We next randomly construct a derangement, say *D*. We then apply $D^{-1}$ to the columns of *M* to obtain $D^{-1}M$. This matrix has the property that its diagonal elements correspond to the row costs of D, i. e., if *(a, D(a))* = 17 in *M*, then the entry *(a, a)* has the value 17 in $D^{-1}M$. Furthermore, if *(a b c)* is a cycle in $D^{-1}M$, then *(a, D(b)), (b, D(c)), (c, D(a))* are corresponding entries in *M*. We thus can obtain admissible permutations in $D^{-1}M$. We define the value of the arc *(a, b)* of $D^{-1}M$ with that of the entry *(a, D(a))* of *M*. In general, we will be using the first column entries of *MIN(M)* unless they already belong to the current derangement with which we are working. In that case, we pick the second column entry. W. l. o. g., suppose that we are trying to find an admissible permutation, s, such that $D_1$ = *Ds* has a smaller value than *D*, i.e., the sum of the values of the arcs of $D_1$ is less than the corresponding sum for *D*. In order to ascertain whether the cycle s will decrease the value of *D*, we define the *DIFF* function. In general,

*DIFF( a ) = d( a , MIN(M)(a , 1)) - d(a , D(a ))*. Thus, *DIFF( a )* ≤ 0. We assign *DIFF(a )* to each vertex *a* of *D*. In Phase 1, our goal in each iteration is to find a permutation, s, whose value is the smallest negatively-valued among all tested. We then define $D_1$ = *Ds* Let *a* be the vertex of *D* with the smallest negative *DIFF* value. We choose the first, second, ..., ([*log n*] + 1)-th smallest entries in row a of *MIN(M)* to obtain the smallest negatively-valued permutation *s*. If we can obtain no negatively-valued permutation, we go back to the *DIFF* values of the vertices of *D* and test the next [*log n*] smallest negatively-valued *DIFF* values. If we still cannot obtain a negatively-valued permutation, we go on to Phase 2. Assume that we have obtained the derangement $D_i$ in Phase 1. Before continuing, we state theorem 3.1 and its corollary.

**Theorem 3.1** *Let C be cycle containing n arcs. Assume that each arc ( a , b) has been assigned a real value, d( a , b). Then if*

$$V = \sum_{i=1}^{i=n} d(i, C(i)) < 0$$

*there exists at least one value i = i', with 1 ≤ $i_1$ ≤ n such that*

$$\sum_{j=0}^{j=m} d(i' + j, C(i' + j)) < 0$$

where $m = 0, 1, ..., n-1$ and $i' + j$ is represented by its value modulo n.

**Corollary 3.1.** *Suppose that C is a cycle containing n arcs such that*

$$W = \sum_{i=1}^{i=n} d(i, C(i)) < N$$

*Then there exists a value $i = i'$ such that each partial sum $S_m$ has the property that*

$$S_m = \sum_{j=0}^{j=m} d(i' + j, C(i' + j)) < N$$

*always holds*.

We now consider $D_i^{-1} M$. We first subtract $d(a, a)$ from each entry $d(a, j)$ (j=1,2,...,n) of row $a$ where $a$ runs through 1, 2, ..., n. This yields $D_i^{-1} M^-$. Starting with column 1, we search for a negatively-valued entry, say (i, j). Using *MIN(M),* we check to see if

(i) $d(i, j) + d(j, k) < 0$,

(ii) $d(i, j) + d(i, k) < d(i, k)$.

If both conditions are satisfied, we substitute $d(i, j) + d(j, k)$ for $d(i, k)$. Every time we obtain a new addition *(i, k)* to a path of negative value, *P,* such that $|P| + d(i, k)$ is still negative, we place the new arc in an $n \times n$ matrix called *PATH*. After each substitution, we check to see if $d(i, k) + d(k, i)$ is negative. If it is, we have a negative cycle. If $d(i, k) + d(k, i)$ is negative, we can always reconstruct the cycle directly from *PATH*. Once we have passed through a new $d(i, k)$, we write its negative value in italics. Also, if a negative path can no longer be continued, we write its last (and smallest) negative value in italics. We go through columns 1,2,...,*n* consecutively. If we obtain a new negative value at some entry *(i, k)* where *k* is less than the current column *j'* that we are working on, we underline its value. This indicates that we cannot further extend the path defined by *(i, k)* using the remaining *n - j'* columns. If we have gone through all columns from $j = 1$ to $j = n$, we have finished one iteration of this portion of the algorithm. If we have obtained a $n \times n$ negative cycle, say *C,* Let $D_{i+1} = D_i C$. We can obtain $D_{i+1}^{-1} M$ by applying the action of $C^{-1}$ to the columns of $D_i^{-1} M$. If it requires more than one iteration to obtain a negative cycle, we work only with the underlined entries. We can add only one arc obtained from $D_i^{-1} M^-$ to each underlined negative path during each of the second, third, ... iterations. This portion of the algorithm continues until there exists no negative path that can be extended. Suppose that we have obtained negative cycles $S = \{ C_i \mid i=1,2,...,r\}$ where $D_{i+r} = D_i C_1 C_2 ... C_r$. Then $D_{i+r}$ is an optimal solution to the Assignment Problem, i.e., $\sigma_{APOPT} = D_{i+r}$. This signals the end of Phase 2 of the algorithm. In Phase 3 of the algorithm, we wish to obtain an approximation to an optimal tour of $M$. Since we use a variant of the usual F-W algorithm, we name this approximation $\sigma_{FWTSPOPT}$, while an optimal tour is denoted by $\sigma_{TSPOPT}$. The first thing we must do is to obtain an upper bound for $|\sigma_{TSPOPT}|$. We do this by checking to see if any of $D_{i+j}$ (j=r,r-1,r-2,1,0) is an n-cycle. If an n-cycle is obtained, the first one obtained is an upper bound for $\sigma_{FWTSPOPT}$. If no n-cycle is obtained, we go back to the trials that yielded $D_i$ from $D_{i-1}$. We then test to see if any of the negative cycles obtained (other than the smallest negatively-valued one which gave us $D_i$), say $C'$, has the property that $D_{i-1} C''$ is an n-cycle. If not, we go back to $D_{i-2}$ and follow the same procedure. If necessary, we go back to $D$ itself : $D$ is an n-cycle. If none of $D_{i+j}$ (j=0,1,2...,r) is an n-cycle, we repeat Phase 1 n times. We then choose whichever tour obtained has the smallest value. W.l.o.g., assume that $D^*$ is the tour of smallest value that we've been able to obtain. Furthermore, let $m_0 = |D^*| - |\sigma_{APOPT}|$. Once we have obtained an optimal assignment, $D_{i+r}$, $(D_{i+r})^{-1} M$ contains no negative cycles. Thus, we can use corollary 3.1 to obtain cycles, $T_1$, $T_2$, ..., $T_s$ each subpath of which has value no greater than $m_0$. Each time we obtain a new cycle, we test to see if some subset of the cycles obtained, say $p$, $D_{i+r} p$ is a tour, $T_1$. If it is, we replace $m_0$ by $m_1 = |T_1| - |\sigma_{APOPT}|$. We then obtain cycle each of whose subset of arcs is no greater than $m_1$ in value. We continue the algorithm obtaining $m_i$ (i=2,3,...) until we can no longer obtain a cycle less than $m_{i'}$ for some natural number $i'$. The algorithm continues through at most $n-1$ iterations. The smallest tour obtained up until that point is designated as $\sigma_{FWTSPOPT}$.

The following corollary of Theorem 3.2 gives a sufficient condition for $\sigma_{TSTOPT} = \sigma_{FWTSPOPT}$.

**Theorem 3.7** *Suppose that $m_{i'}$ is the last upper bound obtained in Step 3a. Let $\sigma_{TSPOPT} = \sigma_{APOPT} s$ where $C = (a\ a_1\ a_2\ ...\ a_r\ b)$ be an arbitrary disjoint cycle of s where a is a determining vertex of C. Then there always exists a cycle of value no greater than $m_{i'}$ obtained in Step 3(a) that is of the form $C' = (a\ b_1\ b_2\ ...\ b_{r'}\ b)$.*

**Corollary 3.1** *If we can obtain no cycle in Step 3a less in value than $m_0$, then*

$$\sigma_{TSPOPT} = \sigma_{FWTSPOPT}.$$

# Chapter 1

## H-admissible Permutatutions and the HCP.

### Howard Kleiman

### HowardKleiman@qcc.cuny.edu

## 1.1 Introduction

Let *h* and *h'* be two n-cycles in $S_n$, the symmetric group on *n* points. Consider the permutation $s = h^{-1}h'$. Since *h* and *h'* both have the same number of points, *s* is an even permutation, i.e., a permutation containing an even number of cycles of even length. From elementary group theory, every even permutation can be represented as a product of (not necessarily distinct) 3-cycles. Let *G* be a randomly chosen graph containing *n* vertices. Assume that *G* contains a hamilton circuit, i.e., a cycle made from arcs of *G* containing each vertex of *G* precisely once. Call it $H_C$. On the other hand, let $h_0$ be a randomly chosen n-cycle in $S_n$, while $H_0$ is a corresponding cycle whose arcs lie in $K_n$, the complete graph on n vertices. $H_0$ is a pseudo-hamilton circuit of *G*. An $H_0$-admissible permutation, *s*, is a permutation representable as the product $h_0^{-1}h$ where $h_0$ and $h$ are n-cycles in $S_n$ and $H$ contains at least as many arcs as $H_0$ in *G*. Note that if (*a b c*) is an $H_i$-admissible 3-cycle, then $H_i(a)=b$, $H_i(b)=c$, $H_i(c)=a$. Thus, the action of *(a b c)* transforms the arcs (or pseudo-arcs) $(a,H_i(a))$ into $(a,H_i(b))$, $(b,H_i(b))$ into $(b,H_i(c))$, $(c,H_i(c))$ into $(c,H_i(a))$. In theorem 1.10, we prove that there exists a sequence of $H_i$-admissible permutations, $s_i$, *(i = 1,2, ... ,r)* such that $H_i s_i = H_{i+1}, ... , H_r = H_C$.

Copyright 2002

Here each $s_i$ is either a 3-cycle or a product of two disjoint 2-cycles. If $e = (a, H_i(a))$ is an arc of $H_i$ that lies in $K_n - G$, $e$ is a pseudo-arc of $H_i$ and $a$ is a pseudo-arc vertex. Theorem 1.10 is an existence theorem. In theorem 1.6, we prove that the probability that a random graph of minimal degree 3 contains fewer than two H-admissible 3-cycles having a pseudo-arc vertex is

$$p = \frac{286n^6 - 4326n^5 + 23489n^4 - 80546n^3 + 190342n^2 - 112242n - 27624}{360n^6 - 5040n^5 + 29160n^4 - 89280n^3 + 15264n^2 - 138240n + 51840}.$$

As $n \to \infty$, this probability of success approaches $\frac{143}{180}$. We can eliminate all vertices of degree two in $G$ by constructing the *contracted graph* of $G$, $G'$, where $\delta(G') \geq 3$. The contracted graph is constructed by deleting all edges of $G$ incident to vertices of degree 2. Furthermore, if $P = [a, v_1, v_2, ... v_r, b]$ is a path in $G$ such that $v_1, v_2, ... , v_r$ are all of degree 2, if the edge $[b, a]$ exists, we delete it when constructing $G'$. If $D$ is a directed graph, a vertex , $v$, whose in-degree=out-degree equals one, is equivalent to a vertex of degree 2 in $G$. We say that the *total degree* of $v$ is two. If $P_D = [a, v_1, v_2, ... , v_r, b]$ is a path in $D$ corresponding to $P$, we must delete all arcs emanating from from $a$ or terminating in $b$ while constructing $D'$. If arc $(b, a)$ exists, it must also be deleted. A reasonable question is: Are $G'$ and $D'$ random? In the usual sense, no. Consider the sets $E$ and $N$ where $E$ is the set of edges of degree 2 removed from $G$ during the construction of $G'$, while the elements of $N$ are the remaining edges of $G$ removed from $G$ while constructing $G'$. The elements of $E$ have the property that each edge lies on *every* hamilton circuit in $G$. On the other hand, $N$ contains edges of $G$ each of which lies on *no* hamilton circuit in $G$. However, since every edge of $E$ is part of an $r$-vertex, an arc of it implicitly lies on every pseudo-hamilton and hamilton in $G$. On the other hand, *no* arc of an edge of $N$ lies on a hamilton circuit of $G'$. Thus, in testing a set, $S'$ of arcs in $G'$ for $H_i$-admissibility, no arcs that can't lie on a hamilton circuit of $G$ are in $S'$. The same is not true of a corresponding set of arcs in $G$. It follows that the probability that the first set of arcs forms an $H_i'$-admissible permutation in $G'$ should be as great or greater than that of a corresponding set $S$ in $G$. A corresponding argument can be made for $D$ and $D'$. Using $G'$, we use Algorithm G to almost always obtain a hamilton circuit. An alternate algorithm, $G_{no\ r\text{-vertices}}$, doesn't require a contracted graph $G'$. In this algorithm, let $a$ be a pseudo-arc vertex of degree 2. Then the probability of including an edge $[a\ b]$ of $G$ in a pseudo-hamilton circuit is 1. If $G$ contains a hamilton circuit, we can almost always obtain a finite sequence of $H_i$-admissible permutations, $s_i$, such that $H_{C'}$ is a hamilton circuit in $G$. We further prove that the ability to create a sequence of form $H_0, H_1, \ldots, H_i, \ldots$ is a necessary and sufficient condition for the existence of a hamilton circuit in $G$. The algorithm used is called Algorithm G. An analogous algorithm for digraphs, i.e., graphs in which each edge is given an orientation, is called Algorithm D. In contrast to Algorithm D where backtracking is necessary, if we fail to obtain an $H_i$-admissible permutation in G, we use a rotation from a pseudo-arc vertex. The running time for all three algorithms is $O(n^{1.5}(\log n)^4)$. As we increase the running time of the algorithm polynomially, say by multiplying it by *n*, we decrease the expected value of failure exponentially. From this observation, we make the following conjecture:

**Conjecture 1.1** *Let G be a graph containing a hamilton circuit. Then Algorithm G or Algorithm $G_{no\ r-vertices}$ always obtains a hamilton circuit in at most $O(n^5)$ running time.*

**Conjecture 1.2** *Let D be a digraph containing a hamilton cycle. Then Algorithm D always obtains a hamilton cycle in at most $O(n^5)$ running time.*

These are difficult conjectures to prove. If counter-examples exist, perhaps they might help improve the algorithms. It might be necessary to test *all* arcs out of each vertex. In the latter case, if a counter-example exists, we should be able to obtain it in $O(n^5)$ running time.

## 1.2 Theorems

We call the type of random graph constructed by Bollabás a *Boll graph*.

An analogous theorem for directed random graphs, theorem Frieze-Boll, was proven by Frieze in [13]. Henceforth, we assume that the set of vertices of each graph or directed graph is *V*. We define the randomness of our choice of edges as Boll does in [4]:

**Theorem 1.1.** *Let G be a random graph (random directed graph) containing a pseudo-hamilton circuit (cycle) H. Then the probability that a pseudo-3-cycle is H-admissible is*

$$\frac{n-3}{2(n-2)}$$

**Proof.** Let *H* be a pseudo-hamilton circuit represented by equi-distantly-spaced points traversing the circle in a clock-wise manner. Now construct a random chord *(1,H(j))* which represents an oriented edge or arc of *G*. Here $1 < j \leq n$. Thus, we cannot let *j* = 2, since (1,2) lies on the circle. Thus, the probability, $p_1$, that *H(j)* is chosen as the terminal point of the arc is $\frac{1}{n-2}$. Next, randomly construct an arc, *(j, H(k))*. It is easily shown by constructing hσ that *S = {(1, H(j)), (j, H(k))}* defines an *H*-admissible pseudo-3-cycle, *s = (1 j k),* if and only if the two arcs intersect in the circle *H*. The probability, *p,* that they intersect is

(Pr(We randomly choose (1, H(j)).)( Pr((j, H(k)) intersects H(j) in H.)

We are assuming that no arc chosen is an arc of the directed hamilton circuit *H*. Thus, if *m = n-2 and j' = j-2,*

$$p_2 = (\frac{1}{n-2})(\sum_{j=3}^{j=n}\frac{n-j}{n-2}) = \frac{(n-2)(n-3)}{2(n-2)^2} = \frac{n-3}{2(n-2)}$$

The following theorem of W. Hoeffding is given in [16]:

**Theorem 1.2.** *Let B(a,p) denote the binomial random variable with parameters a and p with*

$$BS(b,c;a,p) = Pr(b \leq B(a,p) \leq c)$$

*Then*

(i)          $BS(0,(1-\alpha)ap; a, p) \leq exp(-\alpha^2 ap/2)$

(ii)         $BS((1+\alpha)ap, \infty; a, p) \leq exp(-\alpha^2 ap/2)$

where $0 < \alpha < 1$.

**Theorem 1.3** *Let v be a randomly chosen vertex of a Boll random graph, $G_{m^*}$. Then the probability that v has precisely two edges of $G_{m^*}$ incident to it is at most*

$$\frac{\log(cn(\log n)^2)}{2n}$$

**Proof..** We first define hypergeometric probability.

Consider a collection of $N = N_1 + N_2$ similar objects, $N_1$ of them belonging to one of two dichotomous classes (say red chips), $N_2$ of them belonging to the second class (blue chips). A collection of $r$ objects is selected from these $N$ objects at random and without replacement. Given that

$$X \varepsilon N, \ x \leq r, \ x \leq N_1, \ r-x \leq N_2,$$

find the probability that exactly $x$ of these $r$ objects is chosen. $Pr(X = x)$ is given by the formula

$$Pr(X = x) = \frac{\binom{N_1}{x}\binom{N_2}{r-x}}{\binom{N}{r}}$$

where $x$ objects are red and $r - x$ objects are blue. Let $v$ be a randomly chosen vertex of $G_m$. We wish to obtain the probability that exactly two edges in $G_m$ are incident to it. Let

$N = N_1 + N_2$ where $N$ is the sum of the degrees of all vertices in $K_n$ and $N_1$ is the degree of $v$ ***in*** $K_n$. $r$ equals twice the minimum number of edges in $G_m$ for which $G_m$ is 2-connected as $n \to \infty$. Thus,

$$N = 2\binom{n}{2} = n(n-1)$$

$N_1 =$ the number of edges in $K_n$ incident to $v = n - 1$

$N_2 = N - N_1 = n(n-1) - (n-1) = (n-1)(n-1)^2$

$r =$ the sum of the degrees of the vertices in $G_{m^*}$

$$= \text{at most } [n\log(cn(\log n)^2)]$$

$$x = 2$$

In the definition of *r,* we assume that *c* is a positive number.

Then

$$Pr(X = 2) = \frac{\binom{n-1}{2}\binom{n^2 - 2n + 1}{[n\log(cn(\log n)^2)] - 2}}{\binom{n^2 - n}{[n\log(cn(\log n)^2)]}}$$

From W. Feller [10], using the approximation of the hypergeometric distribution to the binomial distribution when $N \to \infty$, let $p = \frac{n-1}{N} = \frac{1}{n}$ yielding

$$Pr(X = 2) \to$$

$$B(2;N,p) = \binom{[n\log(cn(\log n)^2)]}{2}\left(\frac{1}{n}\right)^2\left(1 - \frac{1}{n}\right)^{[n\log(cn(\log n)^2)] - 2}$$

$$\approx [.5\log^2(cn(\log n)^2)]\exp(-\log(cn(\log n)^2))$$

$$\approx \frac{[\log(cn(\log n)^2)]^2}{2cn(\log n)^2}$$

concluding the proof.

**Corollary 1.3.** *The probability that G contains more than O(log n) vertices of degree 2 approaches 0 as $n \to \infty$.*

**Proof.** Using Hoeffding's Theorem, let $p = \frac{(\log(cn(\log n)^2))^2}{2cn(\log n)^2}$, $a = n\log(cn(\log n)^2)$, $c \to \infty$. Thus,

$$ap = \frac{(\log(cn(\log n)^2))^3}{2c(\log n)^2}.$$ We now simplify $ap$.

$log(cn(log\ n)^2) = log\ c + log\ n + 2log(log\ n)$. But

$$\frac{log\ c + log\ n + 2log(log\ n)}{log\ n} \to 1 \text{ as } n \to \infty$$

Thus for very large $n$, $ap \to \frac{log(cn(log\ n)^2)}{2c} < \frac{log\ n}{c}$. From (ii) of his theorem, the probability, $p''$, that $(1 + \alpha)ap$ vertices are of degree 2 satisfies

$$p'' < exp(\frac{-\alpha^2((log(cn(logn)^2)^3}{6c(logn)})),$$

approaching 0 as $n \to \infty$. But for small $\alpha > 0$,

$$(1 + \alpha)ap < \frac{log\ n}{c}, \text{ concluding the proof.}$$

**Theorem 1.4.** Let $D_{m^*}$ be a Frieze-Boll directed graph. Then, given a randomly chosen vertex, $v$, the probability that a unique arc emanates from $v$ is no greater than $\frac{log\ cn}{cn}$ as $n \to \infty$.

**Proof.** We again use hypergeometric probability. W.l.o.g., let $K_n^D$ be the complete directed graph containing all arcs between any two vertices in $V$. Then let

$N$ = the number of arcs in $K_n^D$

$= n(n - 1)$

$N_1 = N - N_1 = (n-1)(n-1) = (n-1)^2$,

$r$ = the number of arcs in $D_{m^*}$

= at most $n(log\ n + k)$ where $k \to \infty$.

$x = 1$.

*Note*. Since $D_{m^*}$ is a Frieze-Boll digraph, each vertex, , $v$ has the property that $d^-(v) \geq 1,\ d^+(v) \geq 1$. From hypergeometric probability,

$$Pr(X = 1) = \frac{\binom{n-1}{1}\binom{(n-1)^2}{[n(\log n + k) - 1]}}{\binom{n^2 - n}{1}}$$

Again using the approximation of the hypergeometric distribution to the binomial distribution as $N \to \infty$, we obtain

$$Pr(X = 1) \to B(1;N,p) = \binom{[n((\log n) + k)]}{1}\left(\frac{1}{n-1}\right)\left(1 - \frac{1}{n-1}\right)^{[n((\log n) + k)] - 1}$$

$$\to [\log n + \log c] \exp(-[\log n + \log c])$$

$$\approx \frac{\log(cn)}{cn}$$

**Corollary 1.4.** *Suppose that $m^* = n(\log cn)$ where $c \to \infty$ as $n \to \infty$. Then the probability that there exist more than $(\frac{1}{c} + 1)^2 (\log n)^2 \to (\log(n))^2$ unique arcs of $D_{m'}$ emanating from or terminating in a vertex approaches 0 as $n \to \infty$.*

**Proof.** The probability, $p$, that a randomly chosen vertex, $v$, of $D_{m^*}$ has a unique arc emanating from it approaches $\frac{\log cn}{cn}$. The number of arcs in $D_{m^*}$ is at most $a = n(\log n + k) = n(\log n + \log c) = n(\log cn)$. The number of vertices is $n$. Thus, $ap$ is at most $\frac{(\log cn)^2}{c}$. From Hoeffding's theorem 1.4,

$$BS(1 + \alpha)ap, \infty; a, p) \leq \exp(-\frac{\alpha^2 (\log cn)^2}{3c}) \to 0 \text{ as } n \to \infty.$$

The same probability is true for the case where a unique arc terminates in $v$. Suppose that the number of vertices is greater than $(\frac{1}{c} + 1)(\log n)^2)$ vertices. Then

$$\frac{(log cn)^2}{c} = \frac{(log c)^2 + 2(log c)(log n) + (log n)^2}{c} > (\frac{1+c}{c})(log n)^2$$

$$\frac{(log c)^2}{c} + \frac{2(log c)(log n)}{c} > (log n)^2$$

$$\frac{(log c)^2}{c log n} + \frac{2(log c)}{c} > log n$$

which is impossible since the left hand side approaches a constant while the right-hand side approaches $\infty$.

The probability that a unique arc terminates in a vertex, $v$, is the same as the probability that it emanates from $v$. Thus, the total number of such arcs is at most $(\frac{1}{c} + 1)^2 (log n)^2 \to (log n)^2$ as $n \to \infty$.

An H-admissible product of two disjoint (pseudo) 2-cycles (for short an H-admissible *POTDTC* ), say $s = \{(a,b), (c,d)\}$, occurs if and only if the vertices $a, b, c, d$ traverse $H$ in a clockwise manner in one of the following ways:

(i)      $a - c - b - d$,

(ii)      $a - d - b - c$.

It follows that if the vertices are placed on $H$ in positions $\frac{2pj}{n}$ for $j = 1, 2, \ldots, n$, then, w.l.o.g., $[a, H(b)]$ and $[c, H(d)]$ are properly intersecting chords of $H$.

Before going further, we mention the following results, of which the first three are found in Dickson [8]:

(i)      $\sum_{j=1}^{j=n} j = \frac{(n+1)(n)}{2}$

(ii)      $\sum_{j=1}^{j=n} j^2 = \frac{n(2n+1)(n+1)}{6}$

(iii)      $\sum_{j=1}^{j=n} j^3 = \left(\frac{n(n+1)}{2}\right)^2$

(iv) $\sum_{j=1}^{j=n} j^4 = \dfrac{n(n+1)(2n+1)(3n^2+3n-1)}{30}$

(v) $\sum_{j=1}^{j=n} j^5 = \dfrac{n^2(n+1)^2(2n^2+2n-1)}{12}$

**Theorem 1.5.** *Let H be a pseudo-hamilton circuit of a random graph or a random directed graph G. Assume that G contains n vertices and that $e_1$ and $e_2$ are randomly chosen edges of G neither of which is an arc of H. Then the probability that $e_1$ and $e_2$ properly intersect (have no endpoints in common) is* $\dfrac{n-3}{3(n-2)}$.

**Proof.** W.l.o.g., let $e_1 = (1,j)$. Consider the probability, $p$, that $(1, j)$ properly intersects $(r, s)$ where $2 \leq r \leq j-1$ while $j+1 \leq s \leq n$. $j$ can range over the domain $[3,n]$. It follows that given a specific value for $j$, the number of integral values of $r$ is $j-2$: the edge, $[r, H(r)]$, is not permissible by hypothesis; furthermore, the loop $[r, r]$ is not an edge of $K_n$. The number of possible values for $s$ is $n-j$, namely, all vertices not contained in $[1, j]$. Thus, given a fixed value of $j$, the number of successes equals $(n-j)(j-2)$. It follows that the total number of successes is

$$\sum_{j=3}^{j=n} (n-j)(j-2)$$

On the other hand, for a fixed value of $j$, the number of failures equals the number of possible values of $r$ ($j-2$) multiplied by the number of possibilities for $s$. $s$ cannot lie in the closed interval $[n-j,n]$. Furthermore, it must be distinct from $r$ and $H(r)$. Therefore, the number of possibilities for $s$ is $j-2$. Therefore, the number of possibilties for failure is $(j-2)^2$. It follows that, using all values of $j$, the number of possibilities for failure is

$$\sum_{j=3}^{j=n} (j-2)^2$$

Thus, the probability of intersection is

$$\frac{\sum_{j=3}^{j=n}(n-j)(j-2)}{\sum_{j=3}^{j=n}(n-j)(j-2) + (j-2)^2}$$

Now let $j' = j-2$, $m = n-2$. Then the probability of success simplifies to

$$\frac{\sum_{j'=1}^{j'=m} mj' - (j')^2}{m(\sum_{j'=1}^{j'=m} j')} = \frac{\frac{m^2(m+1)}{2} - \left(\frac{m(2m+1)(m+1)}{6}\right)}{\frac{m^2(m+1)}{2}} = \frac{m(m+1)(3m-2m-1)}{m(m+1)(3m)} = \frac{m-1}{3m}$$

$$= \frac{n-3}{3(n-2)}$$

**Theorem 1.6**. *Let G be a random graph with n vertices and $\delta(G) \geq 3$ or a random directed graph, D, containing n vertices where both $\delta^+(D)$ and $\delta^-(D) \geq 2$. Assume that H = (1 2 3 ... n) is a pseudo-hamilton circuit containing one pseudo-arc vertex, n. Then the probability that we can obtain at least two H-admissible permutations containing n is at least*

$$p = \frac{286n^6 - 4326n^5 + 23489n^4 - 80546n^3 + 190342n^2 - 112242n - 27624}{360n^6 - 5040n^5 + 29160n^4 - 89280n^3 + 15264n^2 - 138240n + 51840}.$$

**Proof.** As previously done, let $n$ be placed at twelve o'clock and each vertex i at $\frac{2pi}{n}$ for $i = 1, 2, ..., n.-1$. To consider the worst possible case, assume that the degree of $n$ and $1$ in G are both 3. Alternately, the in-degree and out-degree of $n$ and $1$ in D are, in each case, 2. W. l. o. g., assume that H consists of arcs or pseudo-arcs going in a clock-wise direction. Let the following arcs exist: $c = (i, 1)$, $d = (j, 1)$, $a = (n, k)$, $a' = (n, l)$, $b = (k-1, r)$, $b' = (l-1, s)$. Since $n$ is the only pseudo-arc vertex of *H,* any *H*-admissible permutation must contain it. We thus have the following possibilities for sets of arcs corresponding to *H*-admissible 3-cycles:

(1) {(n, k), (k-1, r), (r-1, 1)},.

(2) {(n, l), (l-1, s), (s-1, 1)},

(3) {(n, k), (k-1, i+1), (i, 1)},

(4) {(n, k), (k-1, j+1), (j, 1)},

(5) {(n, l), (l-1, i+1), (i, 1)},

(6) {(n, l), (l-1, j+1), (j, 1)}.

(7) No *H*-admissible 3-cycles are formed.

The following are generally pseudo-arcs: (1), (r-1, 1); (2), (s-1, 1); (3), (k-1, i+1); (4), (k-1, j+1); (5), (l-1, i+1); (6), (l-1, j+1).

Before going on, we present a set of formulas that will be useful in computing the number of possibilities for various events to occur:

(i) $\sum_{j=1}^{j=n} j = \frac{(n+1)(n)}{2}$

(ii) $\sum_{j=1}^{j=n} j^2 = \frac{n(2n+1)(n+1)}{6}$

(iii) $\sum_{j=1}^{j=n} j^3 = \left(\frac{n(n+1)}{2}\right)^2$

(iv) $\sum_{j=1}^{j=n} j^4 = \frac{n(n+1)(2n+1)(3n^2+3n-1)}{30}$

(v) $\sum_{j=1}^{j=n} j^5 = \frac{n^2(n+1)^2(2n^2+2n-1)}{12}$

The simplest way to obtain *p* is to obtain its complement, *p'*, i.e., the number of possibilities for obtaining at most one *H*-admissible 3-cycle. This is equivalent to obtaining the number of possibilities for at most one of the seven given cases to occur; alternately, the number of cases in which at most one pair of arcs intersect. We first obtain the maximum number of possibilities. We randomly choose two edges incident to n and two edges incident to 1. Each edge chosen cannot be a loop - *(1, 1)* or *(n, n)* - or be an arc of *H* - *(n, 1)* or *(1,2)*. We thus have *n-2* vertices available as the terminal end of the edge. Therefore, we first randomly choose *i* and *j*. There are $\binom{n-2}{2}$ ways of choosing such a pair. The smaller of the two numbers is *i*, the other, *j*. Similarly there are $\binom{n-2}{2}$ ways of choosing *k* and *l*. Once we have chosen the arcs *(i, 1), (j, 1), (n, k), (n. l)*, we must randomly choose *(k-1, r)* and *(l-1, s)*. We have *(n - 2)* possibilities for each choice. Thus, the total number of possible choices is

$$\binom{n-2}{2}\binom{n-2}{2}(n-2)^2 = \frac{(n-2)^4(n-3)^2}{4}$$

$$= \frac{n^6 - 14n^5 + 81n^4 - 248n^3 + 424n^2 - 384n + 144}{4}$$

$$= \frac{360n^6 - 5040n^5 + 29160n^4 - 89280n^3 + 152640n^2 - 138240n + 51840}{1440} \quad (1.1)$$

W. l. o. g., we can assume that *i < j* and *k < l*. If we assume that at most one *H*-admissible 3-cycle occurs, then either none or else precisely one of the seven cases is valid. Suppose that case (6) occurs. We then have *(n, l)* intersecting *(j, 1)* implying that *l < j*. But *k < l*. Therefore, *k* must also be less than *j*. This implies that *both (n, k) and (n, l)* intersect *(j, 1)*. Therefore, we would obtain at least *two H*-admissible 3-cycles. Thus, we may eliminate case (6). Now assume that case (3) occurs. Then *(n, k)* intersects (i, 1). Therefore, *k < i*. But *i < j*. implying that *(n, k)* also intersects *(j, 1)*. We may thus eliminate case (3). Next, assume that case (5) occurs. Then *(n, l)* intersects *(i, 1)*. Thus, *l < i*. But *i < j*. Thus, *l < j*. Thus, we would again have two intersections and, therefore, two H-admissible 3-cycles.

Thus, the only cases remaining are (1), (2), (4) and (7) which are illustrated in figures 1.1, 1.2, 1.4, 1.7, respectively.

.

First consider case (1) or case (2). Since there can be at most one *H*-admissible 3-cycle, if *(n, k)* intersects *(k-1, r)*, then $k \geq j$. Otherwise, we would have at least two intersections. Similarly, if *(n, l)* intersects *(l-1, s)*, then $l > k \geq j$. Thus, for precisely case (1) or case (2) to occur, $i < j \leq k < l$. We now consider in how many ways that can occur. First, suppose that $i < j < k < l$. Then the only construction of arcs in which none intersect is *(i, 1), (j, 1), (n, k), n(, l)*. Thus, we must obtain the number of ways in which we can obtain four points from among *n-2*, i. e., $\binom{n-2}{4}$ and multiply it by $(n-2)^2$. The latter is the number of ways in which we can construct *b = (k-1, r)* and *b' = (l-1, s)*. The product obtained is

$$\frac{n^6 - 18n^5 + 131n^4 - 494n^3 + 1020n^2 - 1096n + 480}{24} \quad (1.2)$$

Next, it may occur that *j=k*. Thus, we compute $\binom{n-2}{3}(n-2)^2$:

$$\frac{n^5 - 13n^4 + 66n^3 - 164n^2 + 200n - 96}{6} = \frac{4n^5 - 52n^4 + 264n^3 - 656n^2 + 800n - 384}{24} \quad (1.3)$$

Adding (1.2) to (1.3), we obtain

$$\frac{n^6 - 14n^5 + 79n^4 - 230n^3 + 364n^2 - 296n + 96}{24} \quad (1.4)$$

Given the number of possibilities that $i < j \leq k < l$, four cases are possible:

(1) *a* intersects *b*, but *a' doesn't* intersect *b'*;

(2) *a' * intersects *b'*, but *a doesn't* intersect *b*;

(3) *a doesn't* intersect *b*, and *a' doesn't* intersect *b'*;

(4) *a* intersects *b*, *and a' intersects b'*.

If we subtract the number of cases in which (4) occurs from (1.3), given $i < j \leq k < l$, we obtain the number of cases where we obtain at most one intersection of arcs (*a* and *b* or *a'* and *b'*) which yield an *H*-admissible 3-cycle. We now obtain the number of possibilities that (4) will occur. First, consider the number of values available for *l*. Since $1 < k$, $l = k + 1$ is the smallest possible value for *l*. The largest value that *s* can take is *n-1* since $(l-1, s)$ must properly intersect $(n, l)$. It follows that *n-2* is the largest value that *l* can take. Continuing, if $l = k + 1$, then *s* can take the values *k+2, k+3, ..., n-1*. Thus, the total number of possibilities is *n - k - 2*. Since *l* can vary from *k+1* to *n-2*, the total number of possibilities for $a = (n,l)$ to be intersected by $(l-1, s)$ is

$$\sum_{l=k+1}^{l=n-2} n-l-1$$

Now consider the number of possibilities that $a = (n, k)$ will be intersected by $b = (k-1, r)$. *r* ranges from *k + 1* to *n - 1*. Thus, the number of possibilities for *r* is *n - j - 1*. On the other hand, $l > k \geq j$. Therefore, *k* varies from *j* to *n - 3*. It follows that - assuming that *j* is fixed - the number of possibilities that both *a* and *b* intersect *and a'* and *b'* intersect is

$$\sum_{k=j}^{k=n-3} \sum_{l=k+1}^{l=n-2} (n-k-1)(n-l-1)$$

Continuing, *i* ranges from *3* to *j - 1*, while *j* varies from *4* to *n - 3*: $j \leq k < n - 2$,

Thus, the total number of possibilities is:

$$\sum_{j=4}^{j=n-3} \sum_{i=3}^{i=j-1} \sum_{k=j}^{k=n-3} \sum_{l=k+1}^{l=n-2} (n-k-1)(n-l-1) \quad (1.5)$$

We compute (1.4) one summation at a time. Thus, consider $\sum_{l=k+1}^{l=n-2} (n-l-1)$. Let U = n - l - 1. Our sum then becomes $\sum_{U=n-k-2}^{U=1} U = \sum_{U=1}^{U=n-k-2} U = \frac{(n-k-2)(n-k-1)}{2}$.

Multiplying by *n - k - 1*, we obtain

$$\sum_{k=j}^{k=n-3} \frac{(n-k-2)(n-k-1)^2}{2} \quad (1.6)$$

Let $U = n - k - 1$. Then $k = n - U - 1$. If $k = j$, then $U = n - j - 1$. If $k = n - 3$, then $U = 2$.

We thus obtain

$$\sum_{U=n-j-1}^{U=2} \frac{U^2(U-1)}{2} = \sum_{U=2}^{U=n-j-1} \frac{U^3 - U^2}{2} = \frac{(n-j-1)^2(n-j)^2 - 4}{8} - \frac{(n-j-1)(n-j)(2(n-j)-1) - 6}{12}$$

$$= \frac{3(n-j-1)^2(n-j)^2 - 12 - 2(n-j-1)(n-j)(2(n-j)-1) + 12}{24} \quad (1.7)$$

Thus, our next computations are:

$$\sum_{j=4}^{j=n-3} \sum_{i=3}^{i=j-1} \frac{3(n-j-1)^2(n-j)^2 - 2(n-j-1)(n-j)(2(n-j)-1)}{24}$$

$$= \sum_{j=4}^{j=n-3} \frac{(j-3)[3(n-j-1)^2 - 2(n-j-1)(n-j)(2(n-j)-1)]}{24} \quad (1.8)$$

Let $U = n - j$. Then $j = n - U$, $j - 3 = n - U - 3$. If $j = 4$, then $U = n - 4$. If $j = n - 3$, $U = 3$. Thus, (1.6) yields

$$\sum_{U=n-4}^{U=3} \frac{(n-U-3)[(3U^2(U-1)^2 - 2U(U-1)(2U-1)]}{24}$$

$$= \sum_{U=3}^{U=n-4} \frac{(n-U-3)[(3U^2(U-1)^2 - 2U(U-1)(2U-1)]}{24} \quad (1.9)$$

Multiplying out in (1.9) yields

$$\sum_{U=3}^{U=n-4} \frac{-3U^5 + (3n+1)U^4 + (22-10n)U^3 + (8n-22)U^2 + (6-2n)U}{24} \quad (1.10)$$

$$= \frac{-6n^6 + 126n^5 - 1095n^4 + 5040n^3 - 12957n^2 + 17640n + 9936}{288}$$

$$+ \frac{18n^6 - 273n^5 - 63n^4 - 2708n^3 - 1725n^2 + 3489n + 1254}{720}$$

$$+ \frac{-10n^5 + 162n^4 - 1038n^3 + 3286n^2 + 3324n - 15444}{96}$$

$$+ \frac{192n^4 - 696n^3 + 1046n^2 - 4966n + 9240}{144}$$

$$+ \frac{16n^4 - 212n^3 + 1046n^2 - 4966n + 9240}{144}$$

$$+ \frac{-2n^3 + 20n^2 + 222n - 792}{48}$$

$$= \frac{6n^6 - 66n^5 - 1091n^4 - 4926n^3 + 60882n^2 + 2518n - 18456}{1440} \quad (1.11)$$

As we noted earlier in (1.4), the number of cases in which $i < j \leq k < l$ when expressed as a fraction with 1440 in the denominator is

$$\frac{60n^6 - 840n^5 + 4740n^4 - 13800n^3 + 21840n^2 - 17760n + 5760}{1440} \quad (1.12)$$

It follows that the number of cases in which $i < j \leq k < l$ and at most one *H*-admissible 3-cycle occurs is (1.12) minus (1.11) yielding

$$\frac{54n^6 - 774n^5 + 5831n^4 - 8874n^3 - 39042n^2 - 20278n + 24216}{1440} \quad (1.13)$$

We now consider case (4). Here $i \leq k < j < l$. *a* intersects *d*, but *a'* doesn't intersect *b'*. Given fixed values for *i*, *k* and *j*, the number of cases where $a' = (n, l)$ such that $b' = (l - 1, s)$ *doesn't* intersect *a'* has *s* going through the values $l - 2, l - 3, ..., 1, n$. Thus, the number of such cases is $l - 1$. But *l* varies from $l = j$ to $l = n - 1$. Thus, our first sum is

$$\sum_{l=j}^{l=n-1} l - 1 \quad (1.14)$$

Since $i \leq k < j \leq l$, $k$ varies from $i$ to $j - 1$. Since $a = (n, k)$ doesn't intersect $b = (k-1, r)$, the cases in which this occurs are $r = k-2, k-3, ..., 1, n$. Thus, there are $k - 1$ cases in which $a$ doesn't intersect $b$.

This leads to

$$\sum_{k=i}^{k=j-1}\sum_{l=j}^{l=n-2}(k-1)(l-1) \quad (1.15)$$

$i$ varies from $3$ to $j - 1$, while $j$ varies from $4$ to $n-1$. Therefore, our final expression is

$$\sum_{j=4}^{n-1}\sum_{i=3}^{j-1}\sum_{k=i}^{j-1}\sum_{l=j}^{n-1}(k-1)(l-1) \quad (1.16)$$

First, consider

$$\sum_{l=j}^{l=n-1}(k-1)(l-1) \quad (1.17)$$

Let $U = l - 1$, $l = U - 1$. If $l = j$, then $U = j + 1$. If $l = n - 1$, $U = n$.

We thus obtain

$$\sum_{U=j+1}^{U=n}(k-1)U = (k-1)([\frac{n(n+1)}{2} - \frac{j(j+1)}{2}]) = (k-1)[\frac{n^2+n-j^2-j}{2}] = (k-1)[\frac{(n-j)(n+j+1)}{2}]$$

Next,

$$\sum_{k=i}^{k=j-1}(k-1)[\frac{(n+j)(n-j+1)}{2}]$$

Let $U = k - 1 \Rightarrow k = U - 1$. If $k = i$, $U = i + 1$. If $k = j - 1$, $U = j$. We thus obtain

$$\sum_{U=i}^{U=j}U[\frac{(n+j)(n-j+1)}{2}] = \frac{(n+j)(n-j+1)(j+i)(j-i+1)}{4} \quad (1.18)$$

Now consider

$$\sum_{i=3}^{i=j-1}\{\frac{(n+j)(n-j+1)}{4}\}(j-i+1)(j+i) \quad (1.19)$$

Let $U = j - i + 1$. Then $i = j - U + 1 \Rightarrow i + j = 2j - U + 1$. If $i = 3$, $U = j - 2$. If $i = j - 1$, $U = 2$.

Therefore, our next sum is

$$\sum_{U=2}^{U=j-2} \{\frac{(n+j)(n-j+1)}{4}\}U(2j-U+1) = \sum_{U=2}^{U=j-2} \{\frac{(n+j)(n-j+1)}{4}\}[(2j+1)U - U^2]$$

This yields

$$\{\frac{(n+j)(n-j+1)}{4}\}\{[(2j+1)\frac{(j-2)(j-1)}{2}] - \frac{(j-2)(j-1)(2j-3)}{6}\}$$

$$= \frac{(n+j)(n-j+1)(j-2)(j-1)(4j+6)}{24} \quad (1.20)$$

We now simplify (1.20) and apply our last summation:

$$\sum_{j=4}^{j=n-1} \frac{-2j^5 + 5j^4 + (2n^2 + 2n + 2)j^3 + (-3n^2 - 3n - 11)j^2 + (-5n^2 - 5n + 6)j + (6n^2 + 6n)}{12}.$$

Using formula (i) - (v), we obtain the sum of each term separately and then add up the sums.

$$-\frac{(n-1)(n)}{6}[\frac{2n^4 - 4n^3 + n^2 + n}{12} - \frac{276}{12}] = (n-1)(n)[\frac{-2n^4 + 4n^3 - n^2 - n + 276}{72}]$$

$$\frac{5(n-1)(n)}{12}[\frac{6n^3 - 9n^2 + n + 1}{30} - \frac{98}{30}] = (n-1)(n)[\frac{30n^3 - 45n^2 + 5n - 485}{360}]$$

$$\frac{(n^2 + n + 1)(n-1)(n)}{6}[\frac{n^2 - n}{4} - \frac{36}{4}] = (n-1)(n)[\frac{n^4 - 36n^2 - 37n - 36}{24}]$$

$$\frac{(-3n^2 - 3n - 11)(n-1)(n)}{12}[\frac{2n-1}{6} - \frac{14}{6}] = (n-1)(n)[\frac{-6n^3 + 39n^2 + 23n + 165}{72}]$$

$$\frac{(-5n^2 - 5n + 6)(n-1)(n)}{12}[\frac{1}{2} - \frac{6}{2}] = (n-1)(n)[\frac{25n^2 + 25n - 30}{24}]$$

$$(\frac{6n^2 + 6n}{12})(n-4) = \frac{6n^3 - 18n^2 - 24n}{12}$$

The sum of the expressions on the right-hand side of each equation is

$$(n^2 - n)\{\frac{-10n^4 + 20n^3 - 5n^2 - 5n + 1360}{360}$$
$$+ \frac{30n^3 - 45n^2 + 5n - 485}{360}$$
$$+ \frac{15n^4 - 540n^2 - 555n - 540}{360}$$
$$= \quad + \frac{-30n^3 + 195n^2 + 115n + 825}{360}$$
$$+ \frac{375n^2 + 375n - 450}{360}\}$$
$$+ \frac{180n^3 - 540n^2 - 720n}{360}$$

$$= (n^2 - n)[\frac{5n^4 + 20n^3 - 20n^2 - 65n + 730}{360}] + \frac{180n^3 - 540n^2 - 720n}{360}$$

$$= \frac{5n^6 + 15n^5 - 40n^4 + 35n^3 + 335n^2 - 1430n}{360}$$

$$= \frac{20n^6 + 60n^5 - 160n^4 + 140n^3 + 1340n^2 - 5720n}{1440} \quad (1.21)$$

Subtracting (1.21) and (1.13) from (1.1), we obtain

$$\frac{286n^6 - 4326n^5 + 23489n^4 - 80546n^3 + 190342n^2 - 112242n + 27624}{1440} \quad (1.22)$$

It follows that the probability of obtaining at most one H-admissible 3-cycle is at most

$$\frac{286n^6 - 4326n^5 + 23489n^4 - 80546n^3 + 190342n^2 - 112242n + 27624}{360n^6 - 5040n^5 + 29160n^4 - 89280n^3 + 152640n^2 - 138240n + 51840} \quad (1.23)$$

**Corollary 1.6** *Let D be a random directed graph with $\delta^+(D) \geq 2$, $\delta^-(D) \geq 2$, with H a pseudo-hamilton cycle of D. Then the probability of obtaining at least two H-admissible permutations containing an arbitrary pseudo-arc vertex, v, is the same as in the previous case.*

**Proof.** All we have to do is assume that v has two arcs emanating from it and H(a) has at least two arcs terminating in it.

Before beginning a sketch of the algorithm, we discuss the probability of obtaining an $H_i$-admissible pseudo-3-cycle in $R_3$, a regular 3-out graph obtained from the directed graph, $D_3$. Suppose that $a$ is a pseudo-arc vertex of the pseudo-hamilton circuit, $H_i$, while $H_{i''}$ is the graph consisting of $H_i$ together with all arcs of $R_3$ symmetric to arcs of $H_i$. (For the sake of argument, we assume that each edge of $R_3$ is a pair of symmetric arcs.) Now randomly chose an edge incident to $a$, say $[a, H_i(b)]$, lying in $R_3 - H_{i''}$, and an edge incident to b, say $[b, H_i(c)]$, in $R_3 - H_{i''}$. The question then arises: May we assume that these two edges are actually chosen independently of each other? The answer is yes: Each of these edges was obtained from randomly chosen arcs of $D_3$. It is possible that the arcs in $D_3$ might be $(a, H_i(b))$ and $(H_i(c), b)$. The crucial thing is that the probability that the two arcs intersect approaches $\frac{1}{2}$ as $n \to \infty$. The corresponding edges form an *H*-admissible permutation if and only if they intersect. Thus, as $n \to \infty$, we may assume that randomly chosen edges of the form $[a, H_i(b)]$, $[b, H_i(c)]$ have a probability of $\frac{1}{2}$ of defining an admissible 3-cycle. Similarly, edges of the form $[a, H_i(b)], [c, H_i(d)]$ have a probability of $\frac{1}{3}$ of defining an admissible *POTDTC*. We may thus randomly choose edges from $R_3$ when applying Algorithm G to it. These definitions may be used interchangeably. If *(a, b)* is a directed arc in a directed graph, *(b, a)* is an arc symmetric to *(a, b)*. A random edge chosen in Algorithm G of the next section is always assumed to be incident to a fixed vertex, say v. In that sense, it behaves like an arc emanating from v. Let $h = (a_1\ a_2\ ...\ a_n)$ denote an n-cycle of $S_n$. Suppose we apply an *H*-admissible 3-cycle, *s,* to *h* to obtain *h'* = *hs* where *s* = *(a b c).* Then

$$h' = (a\ h(b)\ ...\ c\ h(a)\ ...\ b h(c)\ ...)\ .$$

Since we have omitted only subpaths belonging to h,

$$a' = [a, H(b), ..., cH(a), ..., bH(c), ...]$$

is an abbreviation of $H'$. Let $H$ be the pseudo-hamilton circuit corresponding to the n-cycle h in $S_n$, while $\sigma$ is the pseudo-3-cycle corresponding to s. Thus, the term H-admissible is used interchangeably with h-admissible. In particular, the pseudo-hamilton circuit, H, corresponding to h and represented by

$$A = [a, H(a), ..., b, H(b), ..., c, H(b), ...]$$

can be replaced by the abbreviation

$$A' = [a, H(b), ..., c, H(a), ..., b, H(c), ...]$$

representing $H'$. As an example, if

$$h = (1\ 2\ 3\ 4\ 5\ 6\ 7\ 8\ 9\ 10\ 11\ 12),\ s = (1\ 4\ 8),$$

$$h' = hs = (1\ 5\ 6\ 7\ 8\ 2\ 3\ 4\ 9\ 10\ 11\ 12)$$

$$= (1\ h(4)\ ...\ 8\ h(1)\ ...\ 4\ h(8)\ ...)$$

Thus, $h'$ can be represented by the abbreviation

$$a' = (1\ h(4) ...\ 8\ h(1)\ ...\ 4\ h(8)\ ...)$$

Correspondingly, since the remaining points of $H'$ all occur in the order in which they occurred in $H$, the pseudo-hamilton circuit, $H'$, is completely determined by the abbreviation

$$A' = (1\ H(4)\ ...\ 8\ H(1)\ ...\ 4\ H(8)\ ...).$$

In particular, s maps $aH(a)$ into $aH(c)$, $bH(b)$ into $bH(c)$), $cH(c)$ into $cH(a)$ where $a = 1, b = 4, c = 8$. Essentially, we are partitioning H into three subpaths that are joined together to form the pseudo-hamilton circuit $H'$. In this example, s is the permutation associated with the abbreviation $A$. Now let $s'' = (2\ 6)(3\ 7)$ be applied to h to obtain

$$h'' = hs'' = (1\ 2\ 3\ 4\ 5\ 6\ 7\ 8\ 9\ 10\ 11\ 12)$$

$$= (1\ 2\ h(6)\ 8\ h(4)\ 6\ h(2)\ 4\ h(8)\ ...)$$

Here we are partitioning h into four subpaths joined together to form $h''$. In general, the format for an abbreviation using an h-admissible *POTDT*, $s'' = (a\ c)(b\ d)$, is

$$hs'' = (a\ h(c)\ ...\ d\ h(b)\ ...\ ch(a)\ ...\ b\ h(d)\ ...)\ .$$

Before going on, we define a rotation. Let

$$S = [v_i, v_{i+1}, v_{i+2}, ..., v_{i+k}]$$

be a subpath of a pseudo-hamilton circuit, $H$, where $v_i$ is a pseudo-arc vertex of $H$. Assume that $[v_i, v_{i+k}]$ is an edge lying in $G - H$. Then

$$S' = [v_i, v_{i+k}, v_{i+k-1}, ..., v_{i+1}]$$

is a subpath of a pseudo-hamilton circuit, $H'$, where $S'$ replaces $S$ of $H$ to yield $H'$. This procedure is called a *rotation* with respect to $v_i$ and $v_{i+k}$. Bollobás, Fenner and Frieze use rotations in [5]. Before, going on, we note that each rotation defines an $H$-admissible permutation. Let $H_R$ define pseudo-hamilton-cycle after the application of the rotation $R$ to $H$. Then there exists an $H$-admissible permutation $r$ such that $Hr = H_R$. Without loss of generality, let

$$H = (1\ 2\ 3\ ...\ n).$$

If $R$ is defined by the arc $(1\ 2m)$, then $r = (1\ 2m-1\ 2m-3\ ...\ 3)(2\ 2m\ 2m-2\ ...\ 4)$ On the other hand, $(1\ 2m+1)$ yields $r = (1\ 2m\ 2m-2\ ...\ 2\ 2m+1\ 2m-1\ ...\ 3)$. In particular, if $m = 1$, then r is $(1\ 3)(2\ 4)$, or $(1\ 2\ 3)$, respectively. Thus, when the *SCORE* values of all $H$-admissible permutations are 0, rotations may exist whose *SCORE* values are 1 or 2. Theorem 1.7 deals with these cases. If $S = [v_1, v_2, ..., v_r]$ is a subpath of $G$ containing interior vertices $v_2, ..., v_{r-1}$ each of which is of degree 2, then $v_1 v_2 ... v_r$ is an *r*-vertex of the *contracted graph*, $G'$. When constructing $G'$, we also delete the edge $[v_1\ v_2]$ (if it exists). Our reason for doing this is that $[v_1\ v_r]$ lies on no hamilton circuit of $G$. In a rotation, the order of the vertices of an *r*-vertex is changed. Thus, if an *r*-vertex is $v_1 v_2 ... v_r$, after a rotation it becomes $v_r v_{r-1} ... v_1$. In constructing the contracted graph of $D$, let $P = [v_1, v_2, \ldots, v_r]$ be a path in $D$ where each $v_i$ $(1 < i < r)$ has total degree two. If $r = 2$, then the in-degree of $v_2$ is one. When constructing $D'$, we must delete all arcs incident to $v_2, v_3, \ldots, v_{r-1}$, as well as all arcs emanating from $v_1$ or terminating in $v_r$ to form the $r$-vertex $v_1 v_2 \ldots v_r$. Furthermore, if $v_1 v_2 \ldots v_r$ is an $r$-vertex in $G'$ and $v_r v_{r+1} \ldots v_s$ is an $s+1$-vertex, we must combine them into $v_1 v_2 \ldots v_r v_{r+1} \ldots v_s$, an $(r + s)$-vertex, that then becomes a vertex in $G'$. Before including this $(r+s)$-vertex in $G'$, we must delete *all* edges of $G$ incident to $v_r$. A reasonable question is: Are $G'$ and $D'$ random? In the usual sense, no. Consider the sets $E$ and $N$ where $E$ is the set of edges of degree 2 removed from $G$ during the construction of $G'$, while the elements of $N$ are the remaining edges of $G$ removed from $G$ while constructing $G'$. The elements of $E$ have the property that each edge lies on *every* hamilton circuit in $G$. On the other hand, $N$ contains edges of $G$ each of which lies on *no* hamilton circuit in $G$. However, since every edge of $E$ is part of an $r$-vertex, an arc of it implicitly lies on every pseudo-hamilton and hamilton in $G$. On the other hand, *no* arc of an edge of $N$ lies on a hamilton circuit of $G'$. Thus, in testing a set, $S'$ of arcs in $G'$ for $H_i$-admissibility, no arcs that can't lie on a hamilton circuit of $G$ are in $S'$. The same is not true of a corresponding set of arcs in $G$. It follows that the probability that the first set of arcs forms an $H_i'$-admissible permutation in $G'$ should be as great or greater than that of a corresponding set $S$ in $G$. A corresponding argument can be made for $D$ and $D'$. Next, we discuss the relationship between $G$, $G - H_i$, and the circle $H_i$ $(0 \leq i \leq n-1)$. First, let $i = 0$. We obtain the random graph, $G$, by randomly choosing $m$ edges and simultaneously placing them in $H_0$. Thus, each pair of edges has probability $\dfrac{1}{\binom{n}{2}} = \dfrac{2}{n(n-1)}$ of being chosen, while each vertex has probability $\dfrac{1}{n}$. The important fact is that these probabilities are independent of the probabilities of each of the other edges. We now note that $G - H_i$ is generally not a random graph. Without loss of generality, let $H_i = (1\ 2\ \ldots\ n)$. Suppose we choose the arc $(1\ j)$.

Now suppose that $G$ is a random graph containing the set of vertices $V = \{1, 2, \ldots, n\}$. Let $j$ be an arbitrary vertex other than 1. Then the probability that an arbitrary edge incident to $j$ has its other endpoint in the interval $[[1, 2, \ldots, n]$ is $\frac{j-1}{n-1}$. Now assume that we have placed the vertices equidistant from each other along the circle $H_i = (1, 2, \ldots, n\}$. Suppose we randomly choose the "arc" $e = (1\ j)$. Furthermore, if an arc of $G$ lies on $H_i$, we recognize this fact before choosing an arc $f = (j-1\ k)$ emanating from $j-1$. Consider $(j-1\ j)$. If it is an edge of $G$, the number of possibilities for $k \leq j$ is $j - 2$. Thus, the probability that $k \leq j$ is at most $\frac{j-2}{n-2}$. On the other hand, if $j - 1$ is a pseudo-arc vertex, then $(j-1\ j)$ doesn't lie in $G$. It follows that the probability $k \leq j$ is again at most $\frac{j-2}{n-2}$. Thus, using Algorithms G, $G_{no\ r-vertices}$, or D, the probability that $f$ intersects $e$ as a chord in the interior of $H_i$ is greater than the probability that this will occur when arc $f$ is chosen randomly from the random graph $G$. A similar situation occurs when we test for $H_i$-admissible *POTDTC*'s. Furthermore, as the algorithms proceeds, the set of arcs, pseudo-arcs and pseudo-arc vertices changes. Assume that every edge is of $G$ is represented by a symmetric pair of arcs. Then, if $v_1 v_2 \ldots v_r$ is an *r*-vertex of $G'$, the same is true of $v_r v_{r-1} \ldots v_1$. Our only caveat is that *we can't have $v_1 v_2 \ldots v_r$ in the same pseudo-hamilton cycle as $v_r v_{r-1} \ldots v_1$ or vice-versa*: $v_1 v_2 \ldots v_r v_{r-1} \ldots v_1$ is a sequence of non-disjoint 2-cycles. Thus, they can't both occur in an *n*-cycle. Given $v_1 v_2 \ldots v_r$ in a pseudo-hamilton cycle, the only way we can obtain $v_r v_{r-1} \ldots v_1$ is if $v_1 v_2 \ldots v_r$ lies on that part of a rotation in which the orientations of the arcs that define it are reversed. The new pseudo-hamilton circuit obtained from replacing subpaths of the form $S$ by *r*-vertices is $H'$. We note that, by construction, the number of pseudo-arc vertices of $H'$ is never greater than the number in $H$. Furthermore, after deleting all edges of $G_{m*}$ from vertices of degree two, all of the remaining arcs of $G_{m*}$ occur in $G'_{m*}$.

Let $t$ be number of vertices of degree 2 in $G_{m*}$. We remove at most $2t$ edges from $G_{m*}$ to form $V_C$, the set of vertices of the contracted graph, $G'_{m*}$. As we will prove in theorems 1.10 and 1.11, given that a graph or digraph contains a hamilton circuit or cycle, $H_C$, there always exists a sequence of $H_i$-admissible pseudo 3-cycles and $POTDTC's$ that yields $H_C$. To obtain an $H$-admissible permutation *(a b c)* requires that *[a, H(b)]* and *[b, H(c)]* properly intersect in a circle along which the vertices of *H* are equally spaced. On the other hand, *H*-admissibility of a *POTDT, (a c)(b d)*, requires that *[a, c]* and *[b, d]* properly intersect. We search each iteration in depth up to *(log n)²* permutations. Thus, using abbreviations when testing for *H*-admissibility, we need only consider at most $24(\log n)^3$ points for the first iteration, $48(\log n)^3$ points for the second one, ... , $48r(\log n)$ points for the *r*-th one. Thus, in *r* iterations, we randomly go through $24(\log n)^3$ points. Suppose a recalculation of $H_i$ occurs every $n^{\frac{1}{2}}$ iterations. Then (as will be shown in more detail in section 1.4), it requires $O(n(\log n)^3)$ r. t. to go through $n^{\frac{1}{2}}$ iterations. On the other hand, it requires $O(i)$ running time (r. t.) to construct a rotation in an abbreviation containing *i* points. It follows that it requires at most $O(n(\log n))$ running time to use a rotation in each of $n^{\frac{1}{2}}$ iterations. It follows that in $O(n^{1.5}(\log n))$ iterations, it requires at most $O(n^{1.5}(\log n)^4)$ r. t. to complete these operations. The latter will also be the r. t. for both Algorithm G and Algorithm D in the next section. The algorithm essentially consists of sequentially obtaining a sequence $H_0, A_1, \ldots A_{[n^5]}, H_{[n^5]}, A_{[n^5]+1}, \ldots, A_{2[n^5]}, H_{2[n^5]}, \ldots n$ which, applying Algorithm G to graphs, the number of pseudo-arc vertices is a monotonically decreasing function. Using Algorithm D for directed graphs, we are required to backtrack after certain iterations. However, using a large number of iterations, the number of successful iterations becomes considerably greater than the number of failures. Thus, we here also replace pseudo-arc vertices by arc vertices. In general, as $n \to \infty$, we approach a hamilton circuit in the former case, and a hamilton cycle in the latter. Definitions and examples of the following data structures comes from Knuth [17]: Given *m* entries, using a balanced, binary search tree, we can, respectively, locate, insert, or delete any element, or rebalance the tree in $O(\log n)$ r. t. In this paper, a LIFO, double-ended queue – henceforth called a *queue* – is a linear list in which all insertions are made at the beginning of the list, while deletions may be made at either end of the list. A randomly chosen edge in Algorithm G of the next section is always assumed to be incident to a fixed vertex, say $v$. In that sense, it behaves like an arc. In general, when choosing an edge incident to an *r*-vertex, all we need do is to determine whether an edge exists which is terminating in or emanating from the *r*-vertex. For instance, if $v_1 v_2 ... v_r$ has an edge of $H_i'$ incident to $v_1$, we must choose an edge incident to $v_r$. On the other hand, if no edge incident to either vertex lies on $H_i'$, then we can choose an edge incident to either vertex. However, if our *r*-vertex is written, $v_1 v_2 ... v_r$ and we choose an edge incident to $v_1$, we must rewrite our *r*-vertex as - $v_1 v_2 ... v_r$ in the *abbreviation* representing $H_i'$. This indicates that the subpath represented by the *r*-vertex goes from $v_r$ to $v_1$. Now consider the case when we're working on a directed graph, $D$. Since we don't use rotations in Algorithm D, we can only accept as arcs in $D'$ those arcs of $D$ which *terminate* in $v_1$ and *emanate* from $v_r$. As we shall see in the next section, in both Algorithm G and Algorithm D, we construct abbreviations each of which represents $[n^{.5}]$ iterations. Furthermore,

(i) If, using Algorithm G, an iteration yields no $H_i$-admissible permutation, then the last action of the iteration is the construction of a rotation, say from

$$S = \{a_1, a_2, ... , a_i, ...\},$$

to

$$S' = \{a_1, a_i, a_{i-1}, ... , a_2, a_{i+1}, a_{i+2}, ...\}.$$

(ii) We change the signs of all vertices in the rotation except the first one. This indicates that we are traversing $H_i'$ in a counter-clockwise manner. Thus, $S'$ should be written

$$\{a_1, a_i, a_{i-1}, ... , a_{i+1}, a_{i+2}, ...\}.$$

In doing so, we must keep in mind that if $a_i$ is a pseudo-arc vertex in S, $[a_j,a_{j+1}]$ is a pseudo-arc. Thus, after the rotation, $[a_{j+1},a_j] = [a_j,a_{j+1}]$ is also a pseudo-arc. It follows that $a_{j+1}$ is a pseudo-arc vertex in S', a subpath of $H'_{i+1}$. Similarly, $[a_j,a_{j-1}]$ is a pseudo-arc if and only if $[a_{j-1},a_j]$ was one in S. Thus, $a_j$ is not necessarily a pseudo-arc of $H_{i+1}$. When we have obtained a hamilton circuit, H', in a contracted graph, G' or D', since the r-vertices contain a subpath of vertices, we respectively obtain a corresponding hamilton circuit (in G) or a hamilton cycle in D.

## 1.3 Algorithms G and D

In this section, we give algorithms for obtaining a hamilton circuit covering a number of cases. In particular, our paper answers two conjectures of Frieze in the affirmative:

(1)    As $n \to \infty$, $D_{2-in,2-out}$ almost always contains a hamilton circuit.

(2)    As $n \to \infty$, $R_3$, the regular 3-out graph, almost always contains a hamilton circuit.

Before going on, we note that both Algorithm G and Algorithm D act upon their respective contracted graphs G' and D'. Consider $G'_{m*}$, the contracted graph of the Boll graph, $G_{m*}$. It is a random graph containing n' vertices each of which is incident to three or more vertices. On the other hand, $R_3$ consists of edges obtained from three random arcs emanating from each vertex. Both algorithms depend upon the intersection of arcs emanating from fixed vertices. Thus, if we can almost always obtain a hamilton circuit in $R_3$, then the same must be true for G'. Similarly, D' has the property that two or more random arcs terminate in each vertex; also, two or more arcs terminate in each vertex. Thus, if Algorithm D is valid for $D_{2-in,2-out}$, then it must be true for the Frieze-Boll digraph, $D_{m*}$. We therefore assume that Algorithm G is applied to $R_3$, while Algorithm D is applied to $D_{2-in,2-out}$. In [12], Fenner and Frieze proved that $R_3$ is 3-connected, while $D_{2-in,2-out}$ is strongly connected. Thus, in both cases, there always exists a simple path connecting any pair of vertices v and w. In general, G and D are both represented as a balanced, binary search tree whose branches are numbered 1 through $n$ together with respective counters that register the number of edges incident to each vertex. Each edge of G is represented as a pair of symmetric arc. Those, two entries are required for each edge of its edges. For simplicity, $D$ denotes $D_{2-in,2-out}$ while $G$ denotes $R_3$. We generally use the word "arc(s)" when discussing both edges and arcs. Since the only edges employed in the algorithm are those incident to a fixed vertex, they are essentially used as arcs. In the algorithm, starting with a randomly chosen initial pseudo-hamilton circuit, $H_0$, we successively construct new ones using $H_i$-admissible permutations, $s_i$ $(i=0,1,2,...)$, to respectively obtain $H_{i+1}$ $(i=1,2,3,...)$. Using theorem 1.6 with respect to $G$ and $D$, we $G'$ replaces $G'$, $D'$ replaces $D$. In general, $G'$ has minimum degree 3, while $\delta^+(D') \geq 2$, $\delta^-(D') \geq 2$. Henceforth, for simplicity, let $n$ be the number of vertices in both the original graph and its contracted graph. For all graphs and directed graphs other than $D'$, we construct $h_0$ by randomly choosing vertices from the balanced, binary search tree obtained from $S = \{1,2,...,n\}$, deleting entries from $S$ when they are chosen, after which we rebalance the search tree. Let 1 be the initial vertex of $ORD(h_0)$ and define $ORD(1) = 1$. If $v_2$ is the second vertex chosen, let $ORD(v_2) = 2$, etc. ... As we construct $h_0$, we place each successive vertex, $v_i$, along with its ordinal number, $ORD(v_i)$, in a balanced binary search tree, $H_0$, in which the search key is $ORD(v_i)$. *Going in a clockwise manner starting at 1 around $h_0$, $ORD(a) < ORD(b)$ implies that a occurs before b.* We, simultaneously, construct a balanced, binary search tree, $ORD^{-1}(H_0)$, in which the numbers from 1 through $n$ occur in sequential order, in which each *integer* is followed by its ordinal number with respect to $h_0$. Here the search key is the numerical value of each number from 1 through $n$. This allows us to access the ordinal value of any given vertex on $h_0$ in at most $O(\log n)$ running time. As we construct $h_0$, we check each new arc in the pseudo-hamilton circuit, $H_0$, to see if it is a pseudo-arc or an arc of $G'$. If it is a pseudo-arc, $(v,v')$, we place $(v,v')$ on a balanced, binary search tree called *PSEUDO* that has the following properties: Each pseudo-arc vertex is placed on *PSEUDO* in two different branches: The first branch, I, contains all initial vertices in increasing order of magnitude, the second, T, all terminal vertices in increasing order of magnitude. On each of these two branches, we construct branches numbered 2, 3, 4, …, that denote the *degree* of the initial (terminal) vertex being placed on I (T). With *D*, we only require placing the *initial* vertices of pseudo-arcs in increasing order of magnitude on I. As in the previous case, we have sub-branches that denote the out-degree of the vertex being placed on I. When we use a rotation in *G'*, we reverse the order of edges, Thus, if *(v,v')* is a pseudo-arc before it is contained in a rotation, *after* the rotation *v'* is a pseudo-arc vertex. The reason for this is that we only choose arcs *emanating* from vertices. Iterations of our algorithm replace $H_0$ by successive pseudo-hamilton circuits, $H_i$ *(i = 1,2,3, ...)* where $H_i$ replaces $H_{i-1}$. Following each iteration in *G'*, if we obtain a new pseudo-arc (or arcs), *(v,v')*, from the $H_i$-admissible permutation chosen during the iteration, we place it in *PSEUDO*. We then use its initial vertex, *v'*, (if there is no sign in front of it) or *v''* (if "-" precedes it) in the next iteration. Before going on, let $SCORE = e_s - a_s$ where $e_s$ is the number of edges in *G' - $H_i$* associated with an $H_i$-admissible permutation, *s,* and $a_s$ is the number of arc vertices in *s*. Define *BACKTRACK* as a queue with the followintg property: If $H_i s_i = H_{i+1}$ where $s_i$ moves at most four points, we place $s_i^{-1}$ in *BACKTRACK*. In Algorithm G, if the *SCORE* values of all permutations are 0,we check the last entry of *BACKTRACK* to make sure we don't choose an $H_i$-admissible permutation that yields $H_{i-1}$. In Algorithm D, when we fail, we use it to obtain $H_{i-1}$. Continuing, if more than one $H_i$-admissible permutation is obtained, we choose the one with the greatest *SCORE* value. We then choose a pseudo-arc vertex of greatest degree from *PSEUDO*. During each iteration, we use up to *log n* arcs emanating from each of at most three pseudo-arcs vertices to construct admissible pseudo 3-cycles or *POTDTC's*. In obtaining the degree of a vertex in *G'*, we don't count arcs *(a b)* where *b* is an *r*-vertex having *opposite* orientation to an *r*-vertex on $H_i$. If the $H_i$-admissible permutation chosen during an iteration has no pseudo-arcs, we randomly pick a pseudo-arc vertex from *PSEUDO* that has largest degree. We can also apply a modified version of Algorithm G to graphs containing vertices whose degree is 2. In that case, if $a$ is a pseudo-arc vertex of degree 2, we first try to obtain an $H_i$-admissible permutation containing an arc $(a\ b)$ in $G$. If we cannot do so, we construct the rotation defined by $(b\ a)$. We then place $R(a\ b)$ in $BACKTRACK$. We note that the probability of obtaining the edge $[a\ b]$ in a hamilton circuit is 1. Also, suppose that one of the following occurs: (1) We obtain an $H_i$-admissible permutation of greatest $SCORE$ value containing a pseudo-arc vertex of degree 2, (2) We obtain an $H_i$-admissible permutation of greatest $SCORE$ value containing an arc terminating in a vertex of degree 2. If either (1) or (2) occurs, we choose such a permutation regardless of the fact that other $H_i$-admissible permutations of the same $SCORE$ value yield a pseudo-arc vertex of greater degree. It may be that we don't need to change any of the orientations to obtain a hamilton circuit. If we have to do so, we can often save time by using rotations. We thus use the following rule: If we obtain an iteration all of whose $SCORE$ values are zero, we choose an $H_i$-admissible permutation, $s_i$, whose pseudo-arc vertex, $a$, has the largest degree. If a rotation out of $a$ has a $SCORE$ value greater than 0, we apply it to $H_i$ to obtain $H_{i+1}$. Otherwise, we apply $s_i$ to $H_i$ to obtain $H_{i+1}$. One other thing is worth mentioning. If only one pseudo-arc vertex is left, we have choose pseudo 3-cycles with one exception – if we obtain a 2-cycle, $(a\ b)$, whose $SCORE$ value is one, we use the distance along $H_i$ of the $ORD$ values from $ORD(a)$ to $ORD(b)$. We then choose the smaller segment $S = \min(ORD(a)\ ORD(b)), (ORD(b)\ ORD(a))\}$. We then randomly choose up to $\log n$ vertices in the interior of $S$, and test up to $\log n$ arcs out of each vertex to test for $H_i$-admissible POTDTC's. Any found would have a $SCORE$ value of at least zero. We do the same test in both Algorithms G and D. In Algorithm D, if we obtain *no* $H_i$-admissible permutation during an iteration, we choose its last entry, $s'$, in $BACKTRACK$ yielding $H_i\ s' = H_{i-1}$. We then continue the algorithm. Assume that the initial pseudo-arc vertex, $a$, of an iteration is of degree 2. If we fail, we use the rotation defined by *(b a)* mentioned earlier. Otherwise, we choose the last entry in BACKTRACK, $s' = s^{-1}$, to obtain $H_i s' = H_{i-1}$. We then continue the algorithm. Assume that $j$ is the number of successes and $i$ the number of failures in $ITER = i + j$ iterations of the algorithm. Given the probability

$$p = \frac{286n^6 - 4326n^5 + 23489n^4 - 80546n^3 + 190342n^2 - 112242n - 27624}{360n^6 - 5040n^5 + 29160n^4 - 89280n^3 + 15264n^2 - 13824n + 51840},$$

$$p' = \frac{74n^6 - 714n^5 + 5674n^4 - 8734n^3 - 175078n^2 - 25998n + 79464}{360n^6 - 5040n^5 + 29160n^4 - 89280n^3 + 15264n^2 - 138240n + 51840},$$

$$p - p' = \frac{212n^6 - 3612n^5 + 17818n^4 - 71812n^3 + 365430n^2 - 86244n - 107088}{360n^6 - 5040n^5 + 29160n^4 - 89280n^3 + 15264n^2 - 138240n + 51840}.$$

Thus, from theorems 1.2 and 1.3 and the Law of Large Numbers, for very large values of $n$, the number of times we succeed minus the number of times we fail is at least $.558 j$. In general, if we have two or more $H_i$-admissible permutations, we do not use one that is the inverse of $s_{i-1}$ where $H_{i-1} \sigma_{i-1} = H_i$. For a large number of iterations, the net number of successful iterations in the first phase of the algorithm is at least

$$(.588)(2n(\log n)) = 1.176n(\log n)$$

Using results from the classical Occupancy Problem, this number is large enough so that we can almost always go through each vertex of $G'$ at least once. This implies that we almost always will obtain a hamilton path. In the second phase of the algorithm, we apply another $2n(\log n)$ iterations to almost always obtain a hamilton circuit in $G'$. For $i \geq 0$, the iterations of our algorithm replace pseudo-arcs on $H_i$ by edges in $G'$. Let $a$ be a pseudo-arc vertex of $H_0$. Then there are always at least six possibilities for choices of edges with respect to $a$ as shown in theorem 1.6:

(1)    Both $[a, H_0(b)]$ and $[b, H_0(c)]$ belong to $G' - H_0$;

(2)　　Both $[a,H_0(b')]$ and $[b',H_0(c')]$ belong to $G' - H_0$;

(3)　　Both $[a,H_0(b)]$ and $[H_0(a),c]$ belong to $G' - H_0$;

(4)　　Both $[a,H_0(b')]$ and $[H_9(a),c']$ belong to $G' - H_0$;

(5)　　Both $[a,H_0(b)]$ and $[H_0(a),c]$ belong to $G' - H_0$;

(6)　　Both $[a,H_0(b')]$ and $[H_0(a),c']$ belong to $G' - H_0$.

A permutation $(a\ b\ c)$ is an $H_0$-admissible pseudo-3-cycle if and only if $ORD(a), ORD(b), ORD(c)$ occur in a clockwise manner as we traverse $H_0$. Given that $G'$ is 3-connected, consider the probability of success in an iteration in $G'$. We pick the worst possible cases. We assume that $H_0$ has precisely one pseudo-arc vertex, $a$. Thus, $H_0(a)$ and $b$ are both arc vertices. Furthermore, assume that $b$ and $b'$ each has two edges incident to it lying on $H_0$, while, since $a$ is a pseudo-arc vertex, $a$ and $H_0(a)$ each has precisely one edge incident to it lying on $H_0(a)$. We now go into greater detail on the construction of $H_0$-admissible permutations. Randomly select up to $log\ n$ edges incident to $a$ and another $log\ n$ incident to $H_0(a)$. (In any case, since every vertex of $G'$ is at least of degree three, we can always test at the minimum the six cases given above.) We systematically test pairs of edges – one incident to $a$, the other incident to $b$, $b'$, or $H_0(a)$ – for $H_0(a)$-admissibility as pseudo-3-cycles. If we obtain an $h_0$-admissible pseudo-3-cycle, say $s_0$, then $h_1 = h_0 s_0$ which implies that $H_1 = H_0 \sigma_0$. Suppose we cannot obtain an $h_0$-admissible pseudo-3-cycle and we have two or more pseudo-arc-vertices on $H_0$, say $v$ and $w$. Then we randomly choose up to $log\ n$ edges incident to $v$ and another set of up to $log\ n$ edges incident to $w$. Using pairs of edges – one edge incident to $a$, the other, incident to $v$ – we test for $h_0$-admissible pairs of POTDT. (We again note that $h_0$-admissibility occurs if and only if the two edges intersect as chords in the interior of the circle upon which $H_0$ is defined.) If $s$ is an $h_0$-admissible POTDT, then $h_1 = h_0 s_0$ implying that $H_1 = H_0\sigma_0$. If we obtain no $H_0$-admissible permutations in an iteration, define its *SCORE* value to be zero. Before going on, we prove the following theorem:

**Theorem 1.7** *Let $a$ be a pseudo-arc vertex of $H_i$, a pseudo-hamilton circuit of a contracted graph, $G'$. Suppose that one of the following holds:*

  *(a) $H_i(a)$ is a pseudo-arc vertex.*

  *(b) If $(a,b)$ is an arc of $G'$, $C = (H_i(a)\ H_i(b))$ is also an arc of $G'$.*

*Then the rotation generated by $(a\ b)$ has a SCORE value greater than 0.*

**Proof. (a)** Let $[a\ a_1\ a_2\ ...\ a_r\ b]$ be a subpath of $H_i$. A rotation with respect to $(a,b)$ yields $[a\ b\ a_r\ a_{r-1}\ ...\ a_1\ a_{r+1}]$. If $[a_i\ a_{i+1}]$ is a pseudo-arc of $H_i$, then $[a_{i+1}\ a_i]$ is a pseudo-arc of $H_i^R$, the pseudo-hamilton circuit obtained from $H_i$ after the rotation. $H_i(a) = a_1$, $H_i(b) = a_{r+1}$. If $B = (a_1\ a_{r+1}) = (H_i(a)\ H_i(b))$ is an arc of $G'$, then two pseudo-arc vertices, $a$ and $H_i(a)$, have become arc vertices. If $B$ is a pseudo-arc, then one pseudo-arc vertex has become an arc vertex. In either case, the *SCORE* value of the rotation is at least one.

  **(b)** Since $C$ is an arc, $H_i(a)$ is an arc vertex. On the other hand, $a$ has been changed from a pseudo-arc vertex to an arc vertex. Thus, the *SCORE* value of the rotation is at least one.

After $[n^5]$ iterations of Algorithm G, we construct $H_1$ out of *ABBREV*. We continue the algorithm using $H_1$ in place of $H_0$. The algorithm continues in this manner until we've gone through $H_{[2n^5(\log(cn))]}$. We then use at most another $2n(\log(n))$ iterations to obtain a hamilton circuit.

*Comment 1.1* If $(a\ b\ c)$ is $H_0$-admissible, we form a balanced, binary search tree called $A_i$ where

$$A_1 = [ORD(a), H_0(ORD(b)), ..., ORD(c), H_0(ORD(a)), ..., ORD(b), H_0(ORD(c))]$$

We now note that $H_0(a\ b\ c) = H_1$, $H_0 = H_1(a\ b\ c)^{-1} = H_1(a\ c\ b)$. When we backtrack in Algorithm D, it is possible that more than one of the edges in $S = \{[a, H_1(c)], [c, H_1(b)], [b, H_1(a)]\}$ is a pseudo-arc. Before selecting an $H_i$-admissible permutation in Algorithm G, we check the unique entry of *BACKTRACK* to see if the permutation selected is that entry. If so, we choose another permutation, or (if necessary) construct a rotation to obtain a new pseudo-arc vertex. Alternately, in Algorithm D, under such a circumstance, we select the last entry in *BACKTRACK*, say $(acb)$, to backtrack to $H_{i-1} = H_i(acb)$. The use of abbreviations shortens the running time of the algorithm. As we proceed in the construction of the abbreviations $A_i$ $(i=1,2,...,[n^5])$, we use the following rules:

Let $v$ be a random element of $V$.

(1) If the successor to $ORD(v)$ doesn't explicitly occur in $A_i$, then

$$ORD(H_i(v)) = ORD(H_0(v)).$$

(2) Otherwise, $ORD(H_i(v)) =$ the successor of $ORD(v)$ on $A_i$.

To clarify the construction of $A_1$ and $A_2$, we use the following example:

Let

$$H_0 = (1\ 14\ 8\ 4\ 3\ 12\ 7\ 13\ 10\ 6\ 11\ 5\ 15\ 9\ 2)$$

$$= (ORD(1)ORD(2)\ ...\ ORD(11)ORD(12)ORD(13)ORD(14)ORD(15))$$

From the latter, we note that the orders of the elements of the permutation are a subset of the natural numbers. It follows that in a first rotation, those ordinal numbers which do not appear in $A_1$ and which are assumed to have negative signs in front of them must be a subset of $\{n, n-1, n-2, ..., i, i-1, ...\}$. Thus, we come to the following rules for rotations:

(1) During a rotation with respect to $v$ and $v_i$, $[v, v_i]$ becomes an arc of $H_0$, while $[v_{i-1}, v_{i+1}]$ may change from a pseudo-arc to an arc or vice-versa. All other arcs or pseudo-arcs formed have the same designation (pseudo-arc or arc) as they had previously. More simply, if $(v_c, v_{c+1})$ was a pseudo-arc before the rotation, then $(v_{c+1}, v_c)$ is a pseudo-arc after the rotation. Similarly, if $(v_c, v_{c+1})$ was an arc before the rotation, then $(v_{c+1}, v_c)$ is an arc after the rotation. The reasoning is straightforward here: Both arcs come from the same edge that either lies on $H_0$ or doesn't lay there.

(2) As mentioned earlier, if an edge of $H_0$, say $[v_c, v_{c+1}]$, lies in *PSEUDO* and has elements which do not explicitly occur in an abbreviation, then if a "-" precedes each of the two elements, $v_{c+1}$ precedes $v_c$ and is a pseudo-arc vertex. If no sign precedes each of them, then $v_c$ is a pseudo-arc vertex.

Continuing with our example,

$$ORD(1)=1, ORD(2) = 14, ORD(3) = 8, ORD(4)\ 4,$$

$$ORD(5) = 3, ORD(6) = 12, ORD(7) = 7, ORD(8) = 13,$$

$$ORD(9) = 10, ORD(10) = 6, ORD(11) = 11, ORD(12) = 5,$$

$$ORD(13) = 15, ORD(14) = 9, ORD(15) = 2.\ ORD(13) = 15, ORD(14) = 9, ORD(15) = 2$$

Suppose $s_1 = (1\ 4\ 7)$. Since $ORD(1) = 1, ORD(4) = 4, ORD(7) = 7,$ if *1, 4, 7* traverse $H_0$ in a clockwise manner, $s_1$ is $H_0$-admissible. Using

$$H_0(ORD(1)) = ORD(2),$$

$$H_0(ORD(4)) = ORD(5),$$

$$H_0(ORD(7)) = ORD(8),$$

$$A_1 = \{ORD(a)H_0(ORD(b))\ ...\ ORD(c)H_0(ORD(a))\ ...\ ORD(b)H_0(ORD(c))\ ...\}$$

$$= \{ORD(1)ORD(5)\ ...\ ORD(7)ORD(2)\ ...\ ORD(4)ORD(8)\ ...\}.$$

Henceforth, we simplify notation by using ordinal numbers in abbreviations. Assume now that *7* is a pseudo-arc vertex and that we want to construct a rotation out of 7. Let *[7, 10]* be an edge of *G'*. $ORD(7) = 7, ORD(10) = 9$. Then the rotation with respect to $ORD(7)$ and $ORD(9)$ transforms $A_1$ into $A_2$ in the following way:

$$A_2 = \{1\ 5\ ...\ \underline{7\ 9}\ 8\ 4\ ---\ \underline{2\ 10}\ ...\}$$

The dashes between 4 and 2 indicate that the ordinal numbers are consecutively decreasing in value. Thus, [7 2 3 4 8 9 10] in $A_1$ becomes [7 9 8 4 3 2 10]. The important thing to note is rotations *never increase* the number of pseudo-arc vertices in $H_1$. First, 7 is no longer a pseudo-arc vertex since it is followed by 9 where $[ORD(7),ORD(9)] = [7,10]$ is an edge of $G'$. On the other hand, $[ORD(2),ORD(10)]$ is generally not an edge of $G'$, although if rotations are done often enough it may - in a particular case - belong to $G'$. Now suppose $[ORD(4),ORD(8)]$ belongs to $G'$. Then $[ORD(8),ORD(4)] = [ORD(4),ORD(8)]$ also belongs to $G'$, while if $ORD(3)$ is a pseudo-arc vertex, then $[ORD(3),ORD(4)] = [ORD(4),ORD(3)]$ doesn't belong to the graph. It follows that if $ORD(3)$ is a pseudo-arc vertex of $A_1$, then $ORD(4)$ is also a pseudo-arc vertex of $A_2$. Now we want to construct $A_3$. Assume that $ORD(2) = 14$ is a pseudo-arc vertex. W.l.o.g., let $s_2 = (14\ 9\ 12)$.

$$14 = ORD(2),\ 9 = ORD(14),\ 12 = ORD(6)$$

We place the respective ordinal numbers in $A_2$ underlined and in italics to see if they occur in a clockwise manner (going from left to right and starting at 1 again if necessary).

$$A_2 = (1\ 5\ \underline{6}\ ...\ 7\ 9\ 8\ 4\ ---\ \underline{2}\ 10\ ...\ \underline{14}\ ...).$$

Going from left to right, we obtain *6, 2, 14*. The numbers occur in a clockwise manner in the circle defined by $H_2$. Thus, $s_2$ is $H_2$-admissible. Before we can construct $A_3$ to represent $H_3$, we must obtain $H_2(ORD(2))$, $H_2(ORD(14))$, $H_2(ORD(6))$. The successor of $ORD(2)$ in $A_2$ is $ORD(10)$. On the other hand, the successor of $ORD(14)$ doesn't explicitly occur in $A_2$. Therefore, it is its successor in $H_0$, namely, $ORD(15)$. Finally, consider $ORD(6)$. The successor of $ORD(6)$ doesn't explicitly occur in $A_2$. Therefore, its successor is its successor in $H_0$, namely, $ORD(7)$. Thus,

$$6 \rightarrow 2 \rightarrow 10, 2 \rightarrow 14 \rightarrow 15, 14 \rightarrow 6 \rightarrow 7$$

yielding

$$A_3 = (1\ 5\ \underline{6\ 10}\ ...\ \underline{14}\ 7\ 9\ 8\ 4\ ---\ \underline{2\ 15})$$

Now let 12 = *ORD(6)* and 3 = *ORD(5)* be pseudo-arc vertices. Let

$$s_3 = (ORD(5)\ ORD(13))(ORD(6)\ ORD(2)) = (3\ 5)(12\ 14)$$

be a product of two disjoint pseudo-2-cycles (*POTDTC*) which we wish to test for $H_3$-admissibility. From

$$A_3 = (1\ \underline{5}\ \underline{6}\ 10\ ...\ \underline{13}\ 14\ 7\ 9\ 8\ 4\ \text{---}\ \underline{2}\ 15)$$

*[5 13]* intersects *[6 2]* in the circle $H_3$. Therefore, $s_3$ is $H_3$-admissible. In this case, $A_4$ is of the form

$$(a\ H_3(b)\ ...\ d\ H_3(c)\ ...\ b\ H_3(a)\ ...\ cH_3(d)\ ...).$$

Continuing,

$$H_3(ORD(5)) = ORD(6),\ H_3(ORD(6)) = ORD(10),$$

$$H_3(ORD(13)) = ORD(6),\ H_3(ORD(2)) = ORD(15).$$

We thus obtain

$$A_4 = (1\ 5\ 14\ 7\ 9\ 8\ 4\ \text{---}\ 2\ 10\ ...\ 13\ 6\ 15)$$

Before we apply a new permutation in an iteration, we check *BACKTRACK* to see if the new permutation, $s_i$, is at the end of the queue. If it is, we construct a rotation out of the initial pseudo-arc vertex to obtain $H_{i+1}$. If we are using Algorithm D, we use $s_i$ to backtrack to $H_{i+1} = H_{i-1}$. After we apply a new permutation, *(a b c)* or *(a c)(b d)* to an abbreviation, we place its inverse *((a c b) or (a c)(b d))* in *BACKTRACK*. In general, if *[v, v\*]* changes from a *pseudo-arc* on $H_i$ to an *arc*, *[v,$H_{i+1}$(v'')]*, on $H_{i+1}$, we delete *[v, v\*]* from *PSEUDO*. On the other hand, if *[v',v\*\*]* is a *new pseudo-arc* on $H_{i+1}$, we place *[v', v\*\*]* in *PSEUDO*. In this manner, we construct *PSEUDO*. Going back to our example, if we fail in an iteration applied to $A_4$ and *ORD(13)* is a pseudo-arc vertex, we construct a new rotation using an edge of *G'* not on $H_4$, say *[ORD(13),ORD(7)]*. After the rotation, *[ORD(6),ORD(9)]* generally becomes a pseudo-arc, while *[ORD(13),ORD(7)]* becomes an arc of $H_5$. If *[ORD(6),ORD(9)]* is an arc of $H_5$, we obtain a pseudo-arc from *PSEUDO*. If *PSEUDO* contains no pseudo-arc, then $H_5$ is a hamilton circuit. In general, using Algorithm G, the number of arcs in *PSEUDO* is a monotonically decreasing function that approaches 1 as we go through all of the vertices in $G'$. In the digraph $D'$, we generally must backtrack. In some cases, when we backtrack, we increase the number of pseudo-arcs by 1 or 2. When we reach the abbreviation $A_{[n^5]}$, using $A_{[n^5]}$ and $H_0$, we construct $H_{[n^5]}$ and ($ORD(H_{[n^5]})$). We then delete $H_0$ and $A_{[n^5]}$; we next use $H_{[n^5]}$ to construct abbreviations $A_{[n^5]+1}$, $A_{[n^5]+2}$, ... , $A_{2[n^5]}$. Using $H_{[n^5]}$ and $A_{2[n^5]}$, we construct $H_{2[n^5]}$, and then delete $H_{[n^5]}$ and $A_{2[n^5]}$. This procedure continues throughout the algorithm. From theorem 1.6, the probability of an iteration yielding at least two $H_i$-admissible permutations is at least .7 when $n = 20$. As $n$ increases, the probability monotonically increases, approaching $\frac{143}{180}$ as $n \to \infty$.

Assume that $n$ is very large. If we use *2n(log n)* iterations - where each iteration chooses up to *2(log n)²* edges randomly chosen to obtain an admissible permutation - using the classical Occupancy Problem, we prove that we almost always successfully pass through every vertex in $V$. If there are *fewer* than *log n* edges incident to a vertex, $a$, we may use each edge in $G' - H_i$ incident to $a$. If there are at least two pseudo-arc vertices in *PSEUDO*, we randomly choose two of them, say $a$ and $b$, with which we construct $H_i$-admissible pseudo-POTDT's. We may use up to *log n* edges incident to each pair of pseudo-arc vertices in the constructions. Using *2n(log n)* iterations, the probability that we will be able to obtain a hamilton path approaches 1 as $n \to \infty$. After obtaining a hamilton path, say $H'_P$, which contains only one pseudo-arc vertex, a hamilton circuit, $H'_C$, is obtained using another *2n(log n)* iterations: Since $G'$ is a random graph, at some point we obtain an $H_i$-admissible 3-cycle, say *(p q r)*, such that each edge in

$$S = \{[p, H_i(q)], [q, H_i(r)], [r, H_i(p)]\}$$

lies in $G' - H'_i$. Let $A_{h'p+j}$ be the abbreviation in which we obtain $H'_C$, while $H_{h'p}$ is hamilton path associated with $A_{h'p+j}$ Using $H_{h'p}$ and the sign in front of an *r*-vertex, we can determine the correct orientation of the subpath represented by an *r*-vertex in the hamilton circuit, $H_C$, of *G'*. We call the algorithm just described, Algorithm G. The problem in applying the algorithm to a directed graph is that we can't use rotations in it. Thus, we actually have to backtrack. When backtracking, we always advance the index of $H_i$ implying $H_i(acb) = H_{i+1} = H_{i-1}$. $(ac)(bd)$ is its own inverse. If we fail, we place it in *BACKTRACK*. We always assume that if *(a b c)* is $H_i$-admissible, *(a c b)* is $H_{i+1}$-admissible. However, if $H_{i+1}$ has fewer pseudo-arc vertices than $H_i$,

$$S_{i+1} = \{[a, H_{i+1}(c)], [c, H_{i+1}(b)], [b, H_{i+!}(a)]\}$$

contains fewer than two arcs. *Although we use the inverse of s $_{i+1}$ ((a c b) or (a c)(b d)) in all cases, we should keep this fact in mind*. On the other hand, if $H_i$ and $H_{i+1}$ have the same number of pseudo-arc vertices, then $S_{i+1}$ always contains at least two arcs. This is why the probability must be greater than $\frac{1}{2}$ that we have at least two admissible permutations in an iteration: if we have only one, say *(a c b)*, it may well be that we have to use it to backtrack when we're working with a directed graph. In the case of a directed graph, *D*, its contracted graph, *D'*, has at least two arcs entering and two leaving each vertex in *V*. We have no trouble defining the orientation of any subpath $S_\beta$ represented by an *r*-vertex.. It has fixed orientation throughout the algorithm. Thus, if an *r*-vertex is $v_1v_2...v_r$, its initial vertex is $v_1$ and its terminal vertex is $v_r$. Since we don't use rotations, *PSEUDO* consists of pseudo-arc vertices – not pseudo-arcs. Otherwise, the algorithm is the same as Algorithm G. We call the algorithm for directed graphs Algorithm D. As mentioned earlier, the probability that *G* has at least two $H_i$-admissible permutations approaches $\frac{143}{180}$ as $n \to \infty$. Correspondingly, the probability for failure approaches $\frac{37}{180}$. It follows that the net number of successful iterations is at least $(\frac{53}{90})$ *ITER* . Here *ITER* is the number of iterations. If *n* is very large, we again use *2n(log n)* iterations in order to successfully go through each vertex in *D'*. We then require another *2n(log n)* iterations to obtain a hamilton circuit.

*Comment 1.2* If $n \geq 20$, and *G* is a random graph containing a hamilton circuit, then the net number of successes (successes minus failures) in *G'* approaches *(.588)(2n(log n)( ITER )*. Thus, we can almost always obtain a hamilton circuit even though *n* is finite. As we shall illustrate in detail later, as we multiply the number of iterations by *n,* the expected value of failure decreases exponentially. In fact, for comparatively small *r,* the expected value obtained by multiplying *the total number of graphs containing n vertices* by the probability of failure approaches zero. Since the expected value of failure is less than 1, this indicates that if *G* contains a hamilton circuit, Algorithm G almost always obtains one.

The running time of both algorithms is $O(n^{1.5}(\log n)^4)$.

**1.4 The Probability of Success.**

In any directed graph, *D*, considered here,

$$\delta^+(D'-H_i) \geq 1, \delta^-(D'-H_i) \geq 1$$

*i = 0,1,2, ...* . Furthermore, if *a* is a pseudo-arc vertex, $\delta^+(a) \geq 2, \delta^-(a) \geq 1$ in *D' - $H_i$*, while $\delta^+(H_i(a)) \geq 1, \delta^-(H_i(a)) \geq 2$. Alternately, if *G* is a graph, $\delta(G' - H_i) \geq 1$ *(i = 0,1, ...)*. Furthermore, if *a* is a pseudo-vertex of $H_i$, then $\delta(a) \geq 2$ in *G' - $H_i$*. Therefore, the probability of success in each iteration when searching in depth in up to $(\log n)^2$ trials is

$$p = \frac{286n^6 - 4326n^5 + 23489n^4 - 80546n^3 + 190342n^2 - 112242n - 27624}{360n^6 - 5040n^5 + 29160n^4 - 89280n^3 + 15264n^2 - 138240n + 51840}.$$

If *n* is very large, it follows by the Law of Large Numbers that the number of successes in *2n(log n)* iterations approaches at least $\frac{286}{180}n(\log n) = \frac{143}{90}n(\log n) > 1.588n(\log n)$, while the number of failures is at most $.413n(\log n)$. It follows that the net number of successes approaches at least $1.175\ 1.175n(\log n)$. Before going on, we describe the Classical Occupancy Problem. In order to obtain the probability of success in Algorithms G and D, we require a formula for the probability of success of a particular case of the occupancy problem. In particular the problem of going through all $n$ vertices of $G'$ or $D'$ using $r$ randomly chosen arcs is equivalent to that of randomly distributing $r$ balls in $n$ boxes. The proof we use and formula obtained is given in Feller [10], Chapter IV. II. Suppose we randomly distribute $r$ balls in $n$ cells. What is the probability that each cell will be occupied? We first note that each arrangement has probability $n^{-r}$ of occurring. Let $A_k$ be the event that cell number $k$ is empty $(k = 1, 2, \ldots n)$. In this instance, all $r$ balls are placed in the remaining $n-1$ cells. This can be done in $(n-1)^r$ different ways. Similarly, there are $(n-2)^r$ arrangements leaving two pre-assigned cells empty, etc. . Accordingly,

$$p_i = (1 - \frac{1}{n})^r, \quad p_{ij} = (1 - \frac{2}{n})^r, \quad p_{ijk} = (1 - \frac{3}{n})^r, \ldots$$

There are $\binom{n}{v}$ ways of selecting $v$ cells from $n$ cells. Thus, for every value of $v \leq n$, the probability of exactly $v$ cells being empty is

$$S_v = \binom{n}{v}(1 - \frac{v}{n})^r$$

Before going on, we note that $S_v$ is the sum of all probabilities containing the same number of subscripts, i.e., $S_v = \sum p_{i_1 i_2 \ldots i_v}$. In Chapter I, Section I, of [10], the following theorem is proven:

**Theorem 1.9.** *Let $S_r$ be the sum of all probabilities containing $r$ subscripts. The probability $P_1$ of the realization of at least one among the events $A_i$ $(i = 1, 2, \ldots, N)$ is given by*

$$P_1 = S_1 - S_2 + S_3 - S_4 - \ldots + (-1)^N S_N$$

It follows that $p_0(r,n)$, **the probability that all cells are occupied** is $1 - S_1 + S_2 - \ldots$ or

$$p_0(r,n) = \sum_{v=0}^{v=n} (-1)^v (1-\frac{v}{n})^r$$

Let $(n)_v = n(n-1)(n-2)...(n-v+1)$. Then $(n-v)^v < (n)_v < n^v$ implying that

$$n^v(\frac{n-v}{n})^r \frac{(n-v)^v}{n^v} < (n)_v (\frac{n-v}{n})^r < n^v(\frac{n-v}{n})^r$$

$$n^v(1-\frac{v}{n})^{v+r} < v! S_v < n^v(1-\frac{v}{n})^r \quad *$$

The following double inequality is from Chapter 2, Section 8 of [10].

**Theorem 1.8.** *Let* $0 < t < 1$. *Then* $e^{-\frac{t}{1-t}} < 1-t < e^{-t}$

Let $t = \frac{v}{n}$. First, consider the inequality

$$v! S_v < n^v(1-t)^r$$

$(1-t)^r < e^{-rt}$. Thus,

$$v! S_v < n^v e^{-rt} = n^v e^{-\frac{rv}{n}}$$

On the other hand,

$$v! S_v > n^v(1-t)^{v+r} > n^v e^{-\frac{t(v+r)}{1-t}} = n^v e^{-\frac{v(v+r)}{n-v}}$$

Thus,

$$(ne^{-\frac{v+r}{n-v}})^v < v! S_v < (ne^{-\frac{r}{n}})^v.$$

$v! S_v$ is tightly bounded by the double inequality given in theorem 1.8.

$$n^v e^{-\frac{(v+r)v}{n-v}} < v! S_v < n^v e^{-\frac{rv}{n}}$$

Let $\lambda = ne^{-\frac{r}{n}}$ remain bounded as $n \to \infty$. Consider

$$\frac{e^{-(\frac{v^2+rv}{n-v})}}{e^{-\frac{rv}{n}}} = e^{-(\frac{v^2+rv}{n-v})+\frac{rv}{n}} = e^{-\frac{v^2(n+r)}{n(n-v)}}$$

Let $r = 1.175n(\log n)$ in $\lambda$. As $n \to \infty$, we obtain

$$e^{-\frac{v^2(n-1.175n(\log n))}{n(n-v)}} = e^{-\frac{v^2(1-1.175\log n)}{n-v}} \to e^0 = 1.$$

Thus, assuming $v$ is held fixed, the ratio of the first and last expression in the double inequality approaches 1. This implies that $\frac{v!S_v}{n^v e^-} \to 1$ as $n \to \infty$, implying that as $n \to \infty$,

$$S_v \to \frac{\lambda^v}{v!} = \frac{\left(\frac{1}{n^{.175}}\right)^v}{v!}. \text{ Thus, as } n \to \infty,$$

$$P_0(n, 1.175n\log(n)) = \sum_{v=0}^{v=n}(-1)^v \binom{n}{v} S_v \to \sum_{v=0}^{v=n}(-1)^v \frac{\left(\frac{1}{n^{.175}}\right)^v}{v!} \to e^{-\frac{1}{n^{.175}}} \to e^0 = 1.$$

It follows that as $n \to \infty$, if $v$ is an arbitrary vertex in $G'$, the probability that in $1.175\,n(\log n)$ successful iterations, we pass through each vertex of $G'$ approaches 1. Thus, the probability of obtaining a hamilton path in $G'$ (and therefore in G) approaches 1 as $n \to \infty$. Assume now that $H_P$ is a hamilton path in $G'$. We now apply another $2n(\log n)$ iterations each containing only to $H_{p+i}$-admissible 3-cycles to obtain one consisting only of arcs in $G' - H_{p+i}$.

The probability of such an arc existing in a 3-cycle is at least $\frac{1}{n}$ implying

Pr (The third arc obtained constructing an $H_{p+i}$ - admissible

$$\text{3-cycle is } not \text{ in } G' - H_{p+i}.) = 1 - \frac{1}{n}.$$

By the Law of Large Numbers, the number of successful iterations is at least 1.175n(log n). Therefore, the probability of *not* obtaining a hamilton circuit is at most

$$(1-\frac{1}{n})^{1.175n(\log n)} < e^{-1.175(\log n)} = \frac{1}{n^{1.175}}$$

Thus, the probability of obtaining an $H_{p+i}$ - admissible 3-cycle containing only arcs in $G' - H_{p+i}$ is $1 - \frac{1}{n^{1.175}}$ which approaches $1$ as $n \to \infty$. It follows that the probability of success of the whole algorithm is at least $e^{-\frac{1}{n^{1.175}}}(1 - \frac{1}{n^{1.175}}) \to 1$

## 1.5 Further Results

Henceforth, given a graph or directed graph, $G$, the edges of the pseudo-hamilton circuit, $H_0'$, lie in $K_n - G'$. Also, any of the following will be denoted by as an $H_i$-admissible permutation: $H_i$-admissible 3-cycle, $H_i$-admissible pseudo-3-cycle, $H_i$-admissible *POTD* 2-cycles, $H_i$-admissible *POTD* pseudo-2-cycles. For simplicity, an $H_i$-admissible *POTD* pseudo-2-cycles may have one or two pseudo-2-cycles

**Theorem 1.9.** *If $D$ is a random directed graph, as $n \to \infty$, there almost always exists a hamilton circuit, $H_C'$ in $D'$ obtainable from an arbitrary pseudo-hamilton circuit, $H_0'$, by sequentially using only $H_i'$-admissible permutations without backtracking.*

**Proof**. Algorithm D proves the theorem except for backtracking. But, as shown in section 1.1,

$$h_C' = (h_0') \prod_{i=1}^{i=c+p} \sigma_i$$

where $p$ stands for the number of permutations which were used in backtracking. But each of these permutations changed $H_i'$ to $H_{i-1}'$. Thus, each of these permutations was the inverse of the preceding one. It follows that, starting at $H_0'$ and removing all of the permutations which were erased by backtracking, we are left with a product of 3-cycles or *POTDTC* each of which is $H_{\alpha_i}'$-admissible ($i = 0,1,2,...,N$) where $N$ is less than $4n[\log n]$. This product yields $H_C'$. From the statement that precedes theorem 1.9, we may assume that $H_0'$ contains only pseudo-arc vertices. Before proving theorem 1.10, we prove the following lemmas:

**Lemma 1.10.1.** *Assume that we obtain $H_i$ by applying an $H_{i-1}$-admissible permutation to $H_{i-1}$ ($i = 1,2, ...$) Let $\sigma_i = (H_i)^{-1} H_C$ where $H_C$ is a hamilton circuit of a graph, G. Then the pseudo-arc vertices (the initial vertices of the pseudo-arcs) in PSEUDO are the same as the points moved by $\sigma_i$.*

**Proof.** All of the edges of $H_0$ lie in $K_n$ - G. Thus, the only edges lying on $H_i$ ($i = 1,2, ...$) are edges of $H_C$. It follows that for any value of $i$, if $\sigma_i = (H_i)^{-1} H_C$, any identity element of the permutation $\sigma_i$ corresponds to an edge, say $e$, of $H_C$ which lies on $H_i$. It follows that $e$ is an edge of $H_i$ whose initial vertex is an arc vertex. On the other hand, all edges of $H_C$ which have initial vertices moved by $\sigma_i$ are pseudo-arc vertices of $H_i$. It follows that

$|PSEUDO| = |\sigma_i|$.

**Lemma 1.10.2.** *Let G be a graph containing a hamilton circuit $H_C$. Let $H_i$ ($i = 0,1,2, ...$) be a pseudo-hamilton circuit of G obtained from successively applying admissible permutations to $H_0$, $H_1$, ... , $H_{i-1}$. Assume $\sigma_i = (H_i)^{-1} H_C$. Suppose all of the points of a disjoint cycle of $\sigma_i$, say C, do not traverse $H_i$ in a counter-clockwise manner, $|C| > 2$, and, either,*

   *(i) going in a clockwise manner, three consecutive points of C, say a, b, c, traverse $H_i$ in a clockwise manner*

*or*

   *(ii) three consecutive points of C, say c, a ,b, traverse $H_i$ in a clockwise manner where ( c a ), ( a b ) are arcs of C.*

*Then (a b c) is an $H_i$-admissible permutation.*

**Proof.** If $a, b, c$ traverse $H_i$ in a clockwise manner, $[a, H_i(b)]$ and $[b, H_i(c)]$ intersect in the circle containing the evenly- spaced vertices of $H_i$. It follows from theorem 1.1 that (a b c) is an $H_i$- admissible permutation. In case (ii), unless $C = (a\ b\ c)$, we assume that $(b\ c)$ rather than $(c\ a)$ plays the role of pseudo-arc in the $H_i$-admissible pseudo 3-cycle $(a\ b\ c)$.

**Corollary 1.10.2**. *Let $|s|$ denote the number of points moved by a permutation, $s$, of $S_n$. Then if*

$H_i = H_{i-1}(a\ b\ c)$ *and* $H_i \sigma_i = H_{i\ C}$, $|\sigma_i| \leq |\sigma_{i-1}| - 2$.

**Proof**. From lemma 1.10.1, the vertices of $\sigma_i$ are precisely the elements of *PSEUDO*. Since $(a b c)$ is $H_i$-admissible, one of the following pairs of arcs lies on C:

$$\{(a\ b), (b\ c)\}, \{(b\ c), (c\ a)\}, \{(c\ a), (a\ b)\}$$

It follows that each of the vertices $a, b, c$ is a pseudo-arc vertex of $H_i$. Application of $H_{i-1}$ to $(a\ b\ c)$ to form $H_{i-1}$ deletes at least two arcs from $C$ as well as two pseudo-arcs or pseudo-arc vertices from $PSEUDO_{i-1}$. If one of the three arcs $(c\ a), (a\ b), (b\ c)$ doesn't belong to $\sigma_{i-1}$, then $H_i \sigma_i = H_C$ where $|\sigma_i| = |\sigma_{i-1}| - 2$. If all three arcs $(a\ b), (b\ c)$ and $(c\ a)$ belong to $\sigma_{i-1}$, then $|\sigma_i| = |\sigma_{i-1}| - 3$.

**Lemma 1.10.3.** *Let C be a disjoint cycle of $\sigma_i$ such that neither condition (i) nor condition (ii) of Lemma 1.8.2 holds for any three consecutive points. Then if $(a\ c)$ is an arbitrary arc of C, there exists at least one arc $(b\ d)$ of a cycle of $\sigma_i$ such that $(a\ c)(b\ d)$ is an $H_i$-admissible permutation.*

**Proof.** If neither condition *(i)* nor condition *(ii)* holds for $C$, then we cannot obtain an $H_i$-admissible permutation of the form $(a\ b\ c)$ using at least two adjacent arcs of $C$. (In what follows, we use the ordinal values in $ORD(H_i)$ of $a, b, c$. Since there is always a way of expressing the three ordinal values in the range *[1, n]*, we assume that in each example this is the case. Given the previous statement, $x > y$ means that $ORD(x) > ORD(y)$. If all sets of three consecutive points of C (for example, $c, a, b$) either traverse $H_i'$ in a counter-clockwise manner with $(i) c > a > b$ or $(ii) c > a, a < b, b > c$, then $(a\ b\ c)$ cannot be $H_i$-admissible. As an example, let $H_i' = (1\ 2\ 3\ 4\ 5\ 6\ 7\ 8\ 9\ 10)$, while $c = 8, a = 5, b = 3$. Then $(a\ b\ c) = (5\ 3\ 8)$ which is not $H_i$-admissible. Since no set of three consecutive points of C form an $H_i$-admissible permutation, every arc, $(a\ c)$, of C must interlace with an arc, *(b d)*, belonging to some cycle (possibly *C*) of $\sigma_i$. For suppose that that is not the case.

$H_i$ applied to the pseudo-2-cycle (or 2-cycle), *(a c)*, yields the product of disjoint cycles

$$(a\, H_i(c)\,...\,)(c\, H_i(a)\,...\,) = P_1 P_2.$$

Here the points of $P_1$ and $P_2$ include all of the elements of $\{1, 2, ..., n\}$. But $H_i \sigma_i = H_C$. Since there exists no arc of $\sigma_i$ interlacing with *(a c)*, $P_1$ contains no point, *p*, which lies on an arc of $\sigma_i$ connecting *p* to some point outside of $P_1$. Thus, $P_1$ is a disjoint cycle of $H_C$ containing fewer than *n* points which is impossible since, by definition, $H_C$ is a Hamilton circuit (an *n*-cycle).

**Example 1.1.** Let $H_i = (1\ 2\ 3\ 4\ 5\ 6\ 7\ 8\ 9\ 10)$ with *(a c)* = (3 7). Then

$H_i(3\ 7) = (1\ 2\ 3\ 8\ 9\ 10)(4\ 5\ 6\ 7)$. Now let *b* = 4 and *d* = 9. Then

$H_i(3\ 7)(4\ 9) = H^* = (1\ 2\ \underline{3}\ 8\ \underline{9}\ 5\ 6\ \underline{7}\ \underline{4})$. Note the occurrence of 3, 7, 4, and 9 in $H^*$.

3 – 9 – 7 – 4 : the points of *(a c)* interlace with those of *(b d)*.

*Note.* A cycle that has an even number of points is an odd permutation, while one with an odd number of points is an even permutation. *(odd)(odd)* = *even*, *(odd)(even)* = *odd*, *(even)(odd)* = *odd*. Thus, a cycle containing *n* points when multiplied by *(a c)* can never yield another cycle containing *n* points.

**Lemma 1.10.4.** *If $H_{i-1}$ is applied to an $H_{i-1}$-admissible permutation, (a c)(b d), in order to obtain $H_i$ then $|\sigma_i| \leq |\sigma_{i-1}| - 2$.*

**Proof.** A product of two disjoint pseudo-2-cycles obtained from $\sigma_{i-1}$ contains four pseudo-arc vertices of $H_{i-1}$ and at least two arcs of $\sigma_{i-1}$. Thus, $|\sigma_i|$ is at least two less than $|\sigma_{i-1}|$.

**Corollary 1.10.4.** *The number of pseudo-arcs or vertices in $PSEUDO_i$ is at least two fewer than the number in PSEUDO.*

**Proof.** From lemma 1.10.1, the number of pseudo-arcs or vertices in $PSEUDO_i$ equals the number of points moved by $\sigma_i$. It follows from lemma 1.10.4 that the corollary is valid.

**Theorem 1.10** *Let G be a graph containing a hamilton circuit, $H_C$, while $H_0$ is a pseudo-hamilton circuit each of whose edges lies in $K_n$ - G. Then there always exists a set, S, of successive $H_i$-admissible permutations, $\sigma_i$ (i=0,1,2,3,...) that obtains $H_C$.*

**Proof.** Lemmas 1.8.2 and 1.8.3 prove the theorem when $H_0$ lies in $K_n$ - G.

**Corollary 1.10** *Let G and $H_C$ be defined as in theorem 1.10. Then the number, N, of sets, $S_i$, (i = 1,2, ...) of $H_i$-admissible permutations obtaining $H_C$ satisfies the inequalities*

$$\frac{2\left[\frac{n}{2}\right]^3 + 3\left[\frac{n}{2}\right]^2 + \left[\frac{n}{2}\right]}{6} < N < \frac{2\left[\frac{n}{3}\right]^3 + 3\left[\frac{n}{3}\right]^2 + \left[\frac{n}{3}\right]}{6}$$

**Proof.** $H_0\sigma_0 = H_C$. The number of points moved by $\sigma_0$ is $n$. We first note that from Feller [?], the expected number of cycles in a randomly chosen permutation is approximately log n. Thus, the expected number of points on a cycle is $\left[\frac{n}{\log n}\right]$. Let A be an arbitrary arc on a cycle of $\sigma$. Then

(a) If A doesn't lie on a 2-cycle, any arc, B, adjacent to it has a probability of $\frac{1}{2}$ of intersecting it.

(b) If B isn't adjacent to A, its probability of its intersecting A is $\frac{1}{3}$.

(c) If A is an arc of a 2-cycle, *TWO*, then any arc of $\sigma$ not on *TWO* has a probability of $\frac{1}{3}$ of intersected it.

Each $H_i$-admissible permutation obtained from arcs of $\sigma$ almost always moves two points; rarely, it moves three points; it almost never moves four points. Given (a), (b) and (c), the expected number of $H_i$-admissible pseudo-3-cycles is $\frac{1}{2}i$, the expected number of $H_i$-admissible pseudo POTDTC's is $\frac{1}{3}(\binom{i}{2} - i)$. Let $N$ be the expected number of sets $S_i$.

$$\sum_{i=3}^{\left[\frac{n}{3}\right]} \frac{1}{2}i + \frac{1}{3}(\binom{i}{2} - i) = \frac{1}{6}(\sum_{i=3}^{i=\left[\frac{n}{3}\right]} i^2) = \frac{2\left[\frac{n}{3}\right]^3 + 3\left[\frac{n}{3}\right]^2 + \left[\frac{n}{3}\right]}{6}.$$

On the other hand, the minimum number of points of σ moved by an $H_i$-admissible permutation is two. Therefore, an upper bound on the numbers of sets is

$$\frac{2\left[\frac{n}{2}\right]^3 + 3\left[\frac{n}{2}\right]^2 + \left[\frac{n}{2}\right]}{6}.$$

**Theorem 1.11** *Let G be a graph containing a hamilton circuit, $H_C$, while $H_0$ is an arbitrary pseudo-hamilton circuit. Then we can always use successive $H_i$-admissible permutations to obtain $H_C$ provided we may assume at most once that each arc vertex of $H_0$ is a pseudo-arc vertex when we construct the $H_i$-admissible permutations.*

**Proof.** The proof is the same as that of theorem G except that, in some cases, we must use arc vertices instead of pseudo-arc vertices at most once to eventually obtain $H_C$.

*Comment 1.3.* The limitation of theorem 1.10 with regards to $H_0$ does not necessarily invalidate the use of Algorithm G. Even if $H_0$ contains some arc vertices, using the contracted graph, $G'$, we can generally pass through the original arc vertices of $H_0'$ as the algorithm proceeds. The best procedure would be to start by using the pseudo-arc vertices of $H_0'$ in constructing admissible permutations. If necessary when constructing admissible permutations, we consider *the arc* vertices as pseudo-arc vertices the first time we use them. Using randomly chosen rotations in Algorithm G as well as arc vertices (if necessary) generally gives us the opportunity to go through each pseudo-arc vertex of $H_0'$ in $2n(\log n)$ iterations of the algorithm. The second phase of the algorithm (to obtain a hamilton circuit in $G'$) generally requires another $2n(\log n)$ iterations. We note that neither theorem 1.9 nor theorem 1.10 requires the use of backtracking or rotations.

**Corollary 1.11.1** *Let $H_0'$ be a pseudo-hamilton circuit of the contracted graph, $G'$, of $G$, where all of the edges of $H_0'$ lie in $K_n - G'$. Then a necessary condition for $G$ to contain a hamilton circuit is that each vertex, $v$, of $G'$ lies on an $H_0'$-admissible permutation of $G'$.*

**Proof.** From theorem 1.10, if $H_0'$ lies in $K_n - G'$, then no arcs of the hamilton circuit, $H_C'$, lie on $H_0'$. Thus, we may use any vertex, $v$, of $G'$ to start the algorithm given in theorem 1.10, i.e., $v$ must lie on an $H_0'$-admissible permutation.

**Theorem 1.12.** *Let $G$ be a random graph of degree $n$ that contains a hamilton circuit. Then if Algorithm G is allowed to run through $12n^3(\log n)$ iterations, the probability of obtaining a hamilton circuit in $G$ is at least*

$$(1 - \frac{2(1 - e^{-\frac{1}{n^{6n^2-1}}})}{n^{6n^2}})(e^{-\frac{1}{n^{6n^2-1}}})(1 - \frac{1}{n^{6n^2}})$$

*The running time of the algorithm is $O(n^{3.5}(\log n)^4)$.*

**Proof.** $G$ is a finite random graph having $n$ vertices. Applying algorithm G, we first use $6n^3(\log n)$ iterations to successfully pass through every vertex in $V$. This yields a hamilton path. We then use another $6n^3(\log n)$ iterations to obtain a hamilton circuit. Successfully passing through a vertex means that an arc, $(a, H_i(b))$, obtained either by the use of an $H_i$-admissible permutation, or a rotation using the vertices $a$ and $H_i(b)$, belongs to $H_{i+1}$. Thus, an iteration always yields at least one successful arc. It follows that $6n^3(\log n)$ iterations yields at least $6n^3(\log n)$ successful arcs. Our next step is to obtain a lower bound for the probability that we go through every vertex of $V$ using $6n^3(\log n)$ iterations. To do this, we use theorem 6.D, Section 6.2 of Barbour, Holst, Janson [3]. In what follows, we rename theorem 6.D in [3], theorem 1.13. Theorem 1.13 is used to obtain the total variation distance between the probability distribution of the random variable, $W_*$, of the occupancy problem and that of the Poisson distribution, $\lambda_*$. An outline of the occupancy problem now follows:

Let r balls be thrown independently of each other into n boxes with probability $p_k$ of hitting the $k-th$ box. Set the random variable, $X_k$, equal to the number of balls hitting the $k$-th box. For the $X's$, the multinomial distribution holds:

$$Pr[X_1 = j_1, X_2 = j_2, ..., X_n = j_n] = \frac{r!}{j_1!...j_n!} p_1^{j_1}...p_n^{j_n}, \sum_{i=1}^{n} j_i = r.$$

Consider the number of boxes hit with at most $m \geq 0$ balls. This random variable may be written as

$$W_* = \sum_{k=1}^{k=n} I_k = \sum_{k=1}^{k=n} I[X_k \leq m]$$

We now make the following definitions:

$\pi_k$ is the probability that the $k$-th box will be hit by at most $m$ balls. If the $k$-th box is hit by at most m balls, we add 1 to the number of boxes hit by at most $m$ balls. Thus, if $E(W_*)$ is the expected number of boxes hit by at most m balls, $E(W_*) = \lambda_* = \sum_{k=1}^{k=n} \pi_k$ where

$$\pi_k = Pr[X_k \leq m] = \sum_{j=0}^{j=m} \binom{r}{j} p_k^j (1-p_k)^{r-j}$$

Let $U_*$ be a Poisson random variable with $E(U_*) = \lambda_*$. In theorem 1.13 [3], we approximate the total variation distance, $d_{TV}$, between the point probabilities of the distributions $L(W_*)$ and Poisson($\lambda_*$).

**Theorem 1.13** For $W_*$ the number of boxes with at most $m$ balls,

$$d_{TV}(L(W_*), Po(\lambda_*)) \leq (1-e^{-\lambda_*})\{ \max \pi_k + \frac{r}{\lambda_*}(\frac{log(r)+mlog(log(r))+5m}{r-log(r)-mlog(log(r))-4m}\lambda_* + \frac{4}{r})^2 \}$$

*where r > log(r) + mlog(log(r)) + 4m.*

We are interested in the case when $m = 0$, i.e., each of the boxes is hit by at least one ball.

Furthermore, since it is equally likely that a randomly thrown ball hits each of the boxes,

$$p_1 = p_2 = ... = p_n = \frac{1}{n}.$$

and

$$\pi_1 = \pi_2 = ... = \pi_n = \pi = Pr[X_k = 0] = \sum_{j=0}^{j=0}(\frac{1}{n})^0(1-\frac{1}{n})^{6n^3 \log(n)}$$

Before going on, let $p^* = \frac{\log(r)+m\log(\log(r))+4m}{r}$. Note that in the proof of theorem 1.13 given in [3], we are given the inequality

$$\sum_{\substack{k=1 \\ p_k > p}}^{k=n} \frac{p_k}{1-p_k} \pi_k < \frac{4}{r}$$

If $r = 6n^3 \log(n)$, $p = \frac{1}{n} > p^*$. We therefore obtain

$$\sum_{k=1}^{k=n} \frac{p_k}{1-p_k} \pi_k = \sum_{k=1}^{k=n} \frac{\frac{1}{n}}{1-\frac{1}{n}}(1-\frac{1}{n})^{6n^3 \log(n)} = \frac{\lambda_*}{n-1}$$

Thus, given $m = 0$,

$$d_{TV}(L(W_*),Po(\lambda_*)) \leq (1 - e^{-\frac{1}{n^{6n^2}-1}})\{\frac{1}{n^{6n^2}} + \frac{r}{\lambda_*}(\frac{\log(r)}{r}\lambda_* + \frac{\lambda_*}{n-1})^2\}$$

*implying that*

$$d_{TV}(L(W_*),Po(\lambda_*)) \leq (1 - e^{-\frac{1}{6n^2-1}})\{\frac{1}{n^{6n^2}} + \frac{\lambda_*}{r}(\log(r)+1)^2\}$$

$$\leq (1-e^{-\frac{1}{6n^2-1}})\{\frac{1}{n^{6n^2}} + \frac{1}{6n^{6n^2+2}\log(n)}(\log(6n^3 \log(n))+1)^2\}$$

$$< \frac{2(1-e^{-\frac{1}{n^{6n^2}-1}})}{n^{6n^2}}$$

Thus, the error obtained in approximating the probability that $6n^3(\log n)$ random arcs will successfully go through each vertex of $V$ using the Poisson distribution with $\lambda_* = n(1-\frac{1}{n})^{6n^3 (\log n)} < \frac{1}{n^{6n^2-1}}$ is less than $\frac{2(1-e^{-\frac{1}{n^{6n^2-1}}})}{n^{6n^2}}$. It follows that the probability of going through all vertices of $V$ and obtaining a hamilton path is at least

$$(1-\frac{2(1-e^{-\frac{1}{n^{6n^2-1}}})}{n^{6n^2}})(e^{-\frac{1}{n^{6n^2-1}}}).$$

Given that we have obtained a hamilton path, $H_p$, we now find the probability of obtaining a hamilton circuit using another $6n^3 (\log n)$ successful arcs. Define $(a, H_{P+j}(b))$ to be the third arc of an $H_{P+j}$-admissible 3-cycle. The probability that it is an edge of $G'$ is at least $\frac{1}{n}$. Now assume $(a, H_{P+j}(b))$ is a successful arc. The question arises: What is the probability that $(b, H_{P+j}(a))$ is an arc or pseudo-arc of $G'$? Again, the probability is at least $\frac{1}{n}$. Thus, during an iteration, the probability of not obtaining a successful arc is at most $(1-\frac{1}{n})$. It follows that the probability that we can't obtain at least one arc during $6n^3 (\log n)$ iterations is at least

$$(1-\frac{1}{n})^{6n^3 (\log n)} < e^{-6n^2 (\log n)} = \frac{1}{n^{6n^2 (\log n)}}$$

It follows that the probability of obtaining a hamilton circuit from a hamilton path is at least $1 - \frac{1}{n^{6n^2}}$. Thus, the probability of success of Algorithm G is at least

$$(1 - \frac{2(1-e^{-\frac{1}{n^{6n^2-1}}})}{n^{6n^2}})(e^{-\frac{1}{n^{6n^2-1}}})(1-\frac{1}{n^{6n^2}})$$

*Comment 1.4.* When we multiply our number of iterations by $n$, we multiply the running time of the algorithm by n. From the proof of theorem 1.10, as we increase the running time polynomially, we decrease the expected value for failure exponentially. .

**Theorem 1.14**. *Let D be a random directed graph containing $n \geq 30$ vertices. Furthermore, assume that D contains a hamilton cycle. Then Algorithm D obtains a hamilton cycle in D with probability at least*

$$(1 - \frac{1}{n^{.32n^2}})(1 - \frac{2(1 - e^{-\frac{1}{n^{2.562n^2-1}}})}{n^{2.562n^2}})(1 - \frac{1}{n^{1.281n^2}})$$

The algorithm is essentially the same as in theorem 1.10 except that we don't use rotations but we do use backtracking. The probability for a successful iteration is

$$p = \frac{286n^6 - 4326n^5 + 23489n^4 - 80546n^3 + 190342n^2 - 112242n - 27624}{360n^6 - 5040n^5 + 29160n^4 - 89280n^3 + 15264n^2 - 13824n + 51840};$$

the probability for failure is

$$p' = \frac{74n^6 - 714n^5 + 5674n^4 - 8734n^3 - 175078n^2 - 25998n + 79464}{360n^6 - 5040n^5 + 29160n^4 - 89280n^3 + 15264n^2 - 138240n + 51840}.$$

Thus, the net probability for success after backtracking is taken into account, $p_{net}$, is

$$p_{net}(n) = p - p' = \frac{212n^6 - 3612n^5 + 17818n^4 - 71812n^3 + 365430n^2 - 86244n - 107088}{360n^6 - 5040n^5 + 29160n^4 - 89280n^3 + 15264n^2 - 138240n + 51840}$$

For $n \geq 30$, $p(n) > .7135599$, while $p'(n) < .28644$ implying that $p_{net}(n) = p(n) - p'(n) > .427$. We now note that every net success yields an $H_i$-admissible permutation that successfully passes through at least two vertices. We now use Hoeffding's Theorem (theorem 1.2) to obtain the probability of obtaining at least

$$(.5)(.427)(6n^3(\log n)) = 1.281n^3(\log n)$$

net successes. Let $\alpha = .5$, $a = 6n^3(\log n)$, $p = .427$. Then

$$BP((0, (.5)(.427)(6n^3(\log n))); 6n^3(\log n), .427) \leq e^{-.125(2.562n^3(\log n))} < \frac{1}{n^{.32n^3}}$$

Thus, the probability of obtaining at least $1.281n^3(\log n)$ $H_i$-admissible permutations is at least $1 - \frac{1}{n^{.32n^3}}$. Since each permutation passes through at least two vertices successfully, the total number

71of vertices passed through successfully is at least $(2)(1.281n^3(\log n)) = 2.562n^3(\log n)$. We next use theorem 1.13 to obtain the maximum error in approximating the point-wise total variation distance between the $W_i$ and Poisson distributions. In this case, $m = 0$, $p = \frac{1}{n}$, $\pi = (1 - \frac{1}{n})^{2.562n^3(\log n)}$,

$\lambda_* = n\pi = n(1-\frac{1}{n})^{2.562n^3(\log n)}$. Thus,

$$d_{TV}(L(W_*), Po(\lambda_*)) \leq \frac{2(1-e^{-\frac{1}{2.562n^2-1}})}{n^{2.562n^2}}$$

In this case, the Poisson probability when $m = 0$ is $e^{-\frac{1}{n^{2.562n^2-1}}}$.

Thus, the probability of obtaining a hamilton path is at least

$$(1 - \frac{2(1-e^{-\frac{1}{2.562n^2-1}})}{n^{2.562n^2}})(e^{-\frac{1}{n^{2.562n^2-1}}})$$

Finally, the probability of not being able to obtain an $H_i$-admissible 3-cycle having all three arcs in $G'$ in $1.281n^3(\log n)$ successful iterations is at most

$$(1 - \frac{1}{n})^{1.281n^3(\log n)} < e^{-1.281n^2(\log n)} = \frac{1}{n^{1.281n^2}}.$$

Therefore, the probability of obtaining a hamilton cycle using $12n^3(\log n)$ iterations is at least

$$(1 - \frac{1}{n^{.32n^2}})(1 - \frac{2(1-e^{-\frac{1}{n^{2.562n^2-1}}})}{n^{2.562n^2}})(1 - \frac{1}{n^{1.281n^2}})$$

The running time of the algorithm is $O(n^{3.5}(\log n)^4)$.

Before going on to conjecture 1.1, we prove the following theorem.

**Theorem 1.15.** *Let $G'$ be the contracted graph of an arbitrary graph $G$ where $G$ contains a hamilton circuit. Suppose $H'$ is a pseudo-hamilton circuit of $G'$. Assume that the following are true:*

*(1) The vertices of $H'$ are equally spaced along a circle C.*

*(2) a is a pseudo − arc vertex of $H'$ .*

(3) $e_1 = [H'(a), H'(d)]$ and $e_2 = [d, H'(e)]$ are edges of $G' - H'$ which do not intersect in C.

(4) $[a,x]$, an edge of $G' - H'$, determines a rotation, $r$, of $H'$ such that $H'(a)$ is a pseudo-arc vertex of $H^* = H'r$.

Then if – going in a clockwise direction – x lies between $H'(e)$ and $d$, $e_1$ and $e_2$ intersect in $H'$ and determine an $H'$-admissible pseudo-3-cycle.

**Proof.** W.l.o.g., let $H' = (1\ 2\ 3\ ...\ \underline{21}\ ...\ \underline{30}\ ...\ \underline{39}\ \underline{40}\ 41\ ...\ n)$.

Assume that $a = 20$, $H'(a) = 21$, $H'(d) = 40$, $d = 39$, $H'(e) = 30$, $x = 35$. Then

$$H^* = (1\ 2\ ...\ 20\ 35\ 34\ 33\ 32\ 31\ \underline{30}...\ \underline{21}\ 36\ 37\ 38\ \underline{39}\ \underline{40}\ ...\ n).$$

We note that $r$ reversed the order of $H'(a) = 21$ and $H'(e) = 30$. Thus, 30, 21, 39, 40 interlace the vertices of $e_1$ and $e_2$, i.e., $e_1$ and $e_2$ intersect in the circle defined by $H^*$. It follows that since $H'(a)$ is a pseudo-arc vertex by hypothesis, the latter two edges define a $H^*$-admissible pseudo-3-cycle.

**Corollary 1.15** Let $e_1$ and $e_2$ be defined as in theorem 1.15. Define

$e_3 = [H'(d), g] = [f, g]$. Assume that $r = [a, x]$ defines a rotation. Let

(a) $x < d$, $e_1$ and $e_2$ intersect,

or

(b) $x > f$, $e_1$ and $e_2$ intersect,

Then $e_1$ and $e_2$ define an $H^*$-admissible pseudo-3-cycle.

**Proof.** [a] On any circle $H^*$ in which --- going in a clockwise direction --- $d$ is the predecessor of $H'(d)$ and $e_1$ and $e_2$ intersect, they define an $H^*$-admissible pseudo-3-cycle.

[b] If $x > f$ and $e_1$ and $e_3$ intersect, $r$ changes $f$ from a successor of $H'(d)$ going clockwise along $H'$ to a predecessor of $H'(d)$ going clockwise along $H^*$. Thus, the two edges define an $H^*$-admissible pseudo-3-cycle.

**Conjecture 1.1.** *Let G be an arbitrary graph that contains a hamilton circuit. Let G' be its contracted graph. Then, using Algorithm G, we can always obtain a hamilton circuit of G' in polynomial time.*

*Comment 1.5.* Theorem 1.15 and its corollary show that it is possible by using a rotation to convert two edges that don't define an admissible permutation into a pair defining an admissible permutation. Finally, we note that in theorems 1.12 and 1.14, the probability of failure decreases exponentially as the running time increases polynomially. An analagous conjecture (Conjecture 1.2) is hypothisized for Algorithm D.

In general, it is useful to cut down on the number of arc vertices in $H_0$. To do this, we first define a $cH_0$-admissible permutation. An $H_0$-admissible permutation most of whose edges lie in $K_n - G'$ is called a $cH_0$-admissible permutation.

The restriction on the number of arc vertices allowed in $H_0'$ is not that difficult to deal with provided that the number of edges in $G'$ is considerably smaller than the number in its complement, $K_n - G'$. We fist note that Theorem G only requires that $G'$ (and therefore $G$) contain a hamilton circuit. Thus, we don't have to randomly construct $H_0'$. Any good heuristic or the algorithm given in this paper may be used to construct $H_0'$. If we wish to obtain a $cH_i' = H_0$ containing a minimum number of arcs in $G'$, we could use quantum computing to eliminate edges in $G'$ from $cH_0'$, replacing them with pseudo-arcs of $G'$ until we obtain an $H_0'$ containing close to the minimum possible number of arcs vertices in any possible $H_0'$.

**Example 1.2** Let $G$ be the set of all edges, $[a,b]$, in $K_{32}$ such that one vertex is even and the other is odd. Assume that all vertices not in italics and underlined are arc vertices of G. Define $cPSEUDO_i$ as the set of arc vertices of $G$ lying in $H_i$.

Let

$$H_C = (1\ 2\ 3\ 4\ 5\ 6\ 7\ 8\ 9\ 10\ 11\ 12\ 13\ 14\ 15\ 16\ 17\ 18\ 19\ 20$$

<p style="text-align:center">21 22 23 24 25 26 27 28 29 30 31 32)</p>

be a hamilton circuit in *G*. Define

$$cH_0 = (1\ 22\ \underline{28}\ 29\ \underline{21}\ 12\ \underline{8}\ 15\ 25\ 27\ \underline{19}\ 6\ \underline{4}\ 23\ \underline{7}\ 24$$

$$\underline{31}\ \underline{18}\ 11\ 13\ \underline{5}\ \underline{2}\ \underline{9}\ 32\ 30\ 20\ \underline{26}\ \underline{17}\ 14$$

$$\underline{16}\ \underline{3}\ \underline{10})$$

$cPSEUDO_0 = \{1, 28, 21, 8, 19, 4, 7, 24,$

$\qquad\qquad 31, 18, 5, 2, 9, 26, 17, 16, 3, 10\}$

If $s_0 = (1\ 28\ 6)$ then

$$cH_1 = (1\ 29\ \underline{21}\ 12\ \underline{8}\ 15\ 25\ 27\ \underline{19}\ 6\ 22\ 28\ \underline{4}\ 23\ \underline{7}\ \underline{24}\ \underline{31}\ \underline{18}\ 11\ 13$$

$$\underline{5}\ \underline{2}\ \underline{9}\ 32\ 30\ 20\ \underline{26}\ \underline{17}\ 14\ \underline{16}\ \underline{3}\ \underline{10})$$

$cPSEUDO_1 = \{21, 8, 19, 4, 7, 24, 31, 18, 5, 2, 9, 26, 17, 16, 3, 10\}$

If $s_1 = (21\ 24)(4\ 3)$ then

$$cH_2 = (1\ 29\ 21\ \underline{31}\ \underline{18}\ 11\ 13\ \underline{5}\ \underline{2}\ \underline{9}\ 32\ 30\ 20\ \underline{26}\ \underline{17}\ 14\ \underline{16}\ 3\ 23\ \underline{7}$$

$$24\ 12\ \underline{8}\ 15\ 25\ 27\ \underline{19}\ 6\ 22\ 28\ 4\ \underline{10})$$

$cPSEUDO_2 = \{31, 18, 5, 2, 9, 26, 17, 16, 7, 8, 19, 10\}$

If $s_2 = (5\ 26)(9\ 8)$, then

$$cH_3 = (1\ 29\ 21\ \underline{31}\ \underline{18}\ 11\ 13\ 5\ \underline{17}\ 14\ \underline{16}\ 3\ 23\ \underline{7}\ 24$$

$$12\ 8\ 32\ 30\ 20\ 26\ \underline{2}\ 9\ 15\ 25\ 27\ \underline{19}\ 6\ 22\ 28\ 4\ \underline{10})$$

$cPSEUDO_3 = \{31, 18, 17, 16, 7, 2, 19, 10\}$

If $s_3 = (31\ 2\ 4)$, then

$$cH_4 = (1\ 29\ 21\ 31\ 9\ 15\ 25\ 27\ \underline{19}\ 6\ 22\ 28\ 4\ \underline{18}\ 11\ 13\ 5\ \underline{17}\ 14\ \underline{16}\ 3$$

$$23\ \underline{7}\ 24\ 12\ 8\ 32\ 30\ 20\ 26\ 2\ \underline{10})$$

$cPSEUDO_4 = \{19, 18, 17, 16, 7, 10\}$

If $s_4 = (18\ 14\ 10)$, then

$$cH_5 = (1\ 29\ 21\ 31\ 9\ 15\ 25\ 27\ \underline{19}\ 6\ 22\ 28\ 4\ 18\ \underline{16}\ 3\ 23\ 7\ 11$$

$$13\ 5\ \underline{17}\ 14\ 24\ 12\ 8\ 32\ 30\ 20\ 26\ 2\ \underline{10})$$

$cPSEUDO_5 = \{19, 16, 17, 10\}$

If $s_5 = (19\ 6\ 2)$, then

$$cH_6 = (1\ 29\ 21\ 31\ 9\ 15\ 25\ 27\ 19\ 3\ 23\ 7\ 11\ 13\ 5\ \underline{17}\ 14\ 24\ 12\ 8\ 32$$

$$30\ 20\ 26\ 2\ 6\ 22\ 28\ 4\ 18\ 16\ \underline{10})$$

$cPSEUDO_6 = \{17, 10\}$

We can define $cH_6$ as $H_0$. We cannot obtain a pseudo-hamilton circuit with fewer than two pseudo-arcs: There always must be a change from an even number to an odd number and from an odd number to an even number.

## 1.6 General $H_0$-Admissible Permutations

In general, given a random graph, $G'$, or a digraph, $D'$, let $H_C$ be a hamilton circuit in $G'$ or a hamilton cycle in $D'$. Define $H_C^{-1}G'$ as the graph obtained from $G'$ in which $H_C^{-1}$ has been applied to every arc of $G'$. The probability that a cycle of length $2r + 1$ in $H_C^{-1}G'$ none of whose arcs lie on $H_C$ is $H_c$-admissible is $\frac{1}{r+1}$. A paper of D. Walkup [22] gives a description of *all* possible $H_C$-admissible permutations. In particular, it gives a recursion formula for finding possible cycle lengths of an *H*-admissible permutation of a fixed length. It also gives the number of $H$-admissible permutations moving a fixed number of points, $m$, allowing us to obtain the probability that an even permutation of length $m$ is $H$-admissible. Thus, This is useful in sections 3a and 3b of Phase 3 of the algorithm in chapter 3 where we need to put together distinct cycles to test for $H_i$-admissible n-cycles. Thus, it seems likely that a random graph, $G$, where $H_c^{-1}G'$ has numerous cycles, has a large number of $H_C$-admissible permutations. It follows that (at least probabilistically) such a graph $G$ must have numerous hamilton circuits. The same thing applies to a random graph $D$..

# 1.7 A Heuristic for the Symmetric Traveling Salesman Problem

In this section, we discuss a heuristic for the symmetric traveling salesman problem using Algorithm G.

(1) Randomly construct an *n*-cycle, say $h_0'$. Let the corresponding pseudo-hamilton circuit be $H_0'$.

(2) Apply a simple heuristic (say, Lin-Kernigan [19]) to $H_0'$ to obtain an approximation to a smallest sum of weights of an *n*-cycle, $H_0$.

We now make the following definitions:

The *defining arcs of an admissible permutation* are two arcs associated with the permutation that intersect. A set of arcs is *good* if the sum of their weights is less than the sum of the weights of arcs on $H_i$ that have the same respective initial vertices. Then the following hold: If the sum of the weights of the arcs associated with an $H_i$-admissible permutation is less than the sum of the weights of the arcs with corresponding initial vertices, then the permutation is good. Furthermore, if a rotation determined by the arc *[x, y]* has the property that

$$w[x,y] + w[H_i(x),H_i(y)] < w[x,H_i(x)] + w[y,H_i(y)],$$

then the rotation is good.

(3) We next make rotations through each of the vertices of $H_i$ until we can no longer make a good rotation.

(4) An iteration now consists of one of the following:

    (a) a good admissible 3-cycle which is not followed by a rotation.

    (b) an admissible 3-cycle such that the sum of the weights of a set of defining arcs which added to the sum of the weights of a rotation defines a good set of arcs.

(5) If $H_i$ is a weighted *n*-cycle before 4(a) and 4(b) have been applied to it, while $H_{i+1}$ is the result after such applications, then

$$W(H_i) - W([a, H_i(a)]) > W(H_{i+1}) - W([H_i(a), H_{i+1}(H_i(a))]).$$

*Note.* We require at least $O(i)$ r. t. to test a sequence of an admissible permutation followed by one or two rotations on $A_i$ to see if a good set of arcs is defined. We then must delete arcs after the testing is done.

As we proceed in (4) and (5), if $W(H_{i_1}) < W(H_0)$, we place $W(H_{i_1})$ and $H_{i_1}$ in a queue. If $W(H_{i_2}) < W(H_{i_1})$, we delete $W(H_{i_1})$ and $H_{i_1}$ from the queue and replace them with $W(H_{i_2})$ and $H_{i_2}$. The algorithm concludes if, starting at a minimum weight $W^* = W(H_{i_j})$, no $H_{i_{j'}}$ lies between $H_{i_j}$ and $H_{i_{j+[n\log n]}}$ where $W(H_{i_{j'}}) < W^*$.

## 1.8 Notes

In [8], Erdös and Renyi give conditions under which a random graph is 2-connected. In [20], Palásti gives conditions for a directed graph to be strongly connected. In [17], Kömlos and Szemerédi prove that if a random graph, $G$, is constructed by randomly choosing edges from the complete graph, $K_n$, until every vertex has a minimum degree of at least two, then as $n \to \infty$, G contains a hamilton circuit with probability approaching 1. A similar theorem was proven by Bollobás in [3]. In [12], Frieze proved an analogous theorem for random directed graphs.

Let $N = \binom{n}{2}$. In [8], Erdös and Renyi give the probability

$$\frac{1}{\binom{N}{k_n}}$$

for randomly choosing a graph, $G$, with n vertices and $k_n$ edges. Let $E(n, k_n)$ be the event "The random graph, $G$, containing $n$ vertices and $k_n$ edges is 2-connected". As $n \to \infty$, they give the following limit distribution for the event $E(n, k_n)$:

Given $k_n = \frac{1}{2}n(\log n + \log(\log n) + c_n)$,

$$= 0 \text{ if } c_n \to -\infty,$$

$$\lim Pr(E(n,k_n)) = \exp(e^{-2c}) \text{ if } c_n \to c,$$

$$= 1 \text{ if } c_n \to \infty.$$

It follows that if $G$ is 2-connected, then each vertex is at least of degree 2. In [20], Palásti proved the following: Let $D_{n,N}$ be a random directed graph, where each edge is equally likely to be chosen from among all edges in $K_n^D$, the complete directed graph on all vertices (including all loops). Then if $N = N_c$ where $N_c = [n(\log n) + c]$ and $c$ is an arbitrary number, the probability that $D_{n,N}$ is strongly connected has the probability $\lim_{n \to \infty} P_{n,N_c} = e^{-2e^{-c}}$. It follows that as $c$ approaches a very large positive number, the above limit approaches 1. In [17], Kómlos and Szemerédi proved the following theorem:

**Theorem KS.** *Let the edges in $K_n$ be numbered $e_1, e_2, \ldots, e_N$ where $N = \binom{n}{2}$. Suppose that each edge in $K_n$ has been chosen in the following manner: The first edge has been chosen randomly with probability $\frac{1}{\binom{n}{2}}$, the second edge has been chosen randomly from the remaining edges with probability $\frac{1}{\binom{n}{2} - 1}$, etc. We stop this process at the first instant when every valence is at least 2, say when m edges have been chosen. Then*

$$\lim_{n \to \infty} Pr(G_m \text{ is hamiltonian.}) = \lim_{n \to \infty} Pr(\delta(G_m) \geq 2) = 1$$

Theorem KS also appears to have been proven by Ajtai, Kömlos and Szeremédi in [1].

In [4], Bollobás proved the following:

**Theorem 1.1a**. *Let $\{e_1, e_2, \ldots, e_N\}$ ($N = \binom{n}{2}$) be a random permutation of $K_n$. If $G_m = \{V, \{e_1, e_2, \ldots, e_{m*}\}\}$ and $m* = \{\min(m: \delta(G_m) \geq 2\}$, then $\lim_{n \to \infty} Pr(G_{m*}$ is hamiltonian.$) = 1$.*

We call the type of random graph constructed by Bollabás a *Boll graph*.

An analogous theorem for directed random graphs, theorem Frieze-ABKS, was proven by Frieze in [13]. Henceforth, we assume that the set of vertices of each graph or directed graph is V. We define the randomness of our choice of edges as Boll does in [4]:

Let

$$p(E) = \{e_i \mid i = 1, 2, \ldots, \frac{n(n-1)}{2}\}.$$

be a random permutation of the edges in $K_n$. With $m < \frac{n(n-1)}{2}$, define $G_m$ to be the random graph with edges $E = \{e_i \mid i = 1, 2, \ldots, m\}$. Define $k(G)$ to be the *vertex connectivity* of G.

Let $\mathbb{N}$ be the set of natural numbers. In [6], Bollobás proved the following theorem:

**Theorem 1.2a.** *Let Q be a monotonically increasing property of graphs and t a function defined by*

$$t(Q) = t(Q;G) = \min (m: G \text{ has } Q).$$

*Then, given $d \in \mathbb{N}$, as $n \to \infty$,*

$$t(\delta(G) \geq d) = t(k(G) \geq d).$$

It follows that if a Boll graph, G, has each vertex of degree at least 2, then G is almost always 2-connected. In [22], Wormald proved that as $n \to \infty$ almost all random graphs on n vertices that are of degree r are r-connected. In [15], Frieze, Jerrum, Molloy and Wormald proved that almost all random regular graphs on n vertices that are of degree 3 have hamilton circuits as $n \to \infty$. Let $D_m$ be a random, directed graph in which $m$ arcs have been randomly chosen to emanate from each vertex. If we change each arc into an unoriented edge, then the resulting graph, $R_m$, is a regular m out-degree graph. In [11], Fenner and Frieze proved that $R_m$ is m-connected and that $D_{\text{i-in, o-out}}$ is strongly connected. These results are necessary to apply Algorithms G and D, respectively, to $R_3$ and $D_{2-in,2-out}$. Let

$$p(E) = \{e_i \mid i = 1, 2, \ldots, \frac{n(n-1)}{2}\}$$

be a random permutation of the edges in $K_n$. With $m < \frac{n(n-1)}{2}$, $G_m$ is the random graph with edges

$$E = \{e_i \mid i = 1, 2, \ldots, m\}.$$

From theorem 1.2a, it follows that if a Boll graph, $G$, has each vertex of degree at least 2, then $G$ is almost always 2-connected. Assume that $G$ is a random graph satisfying the hypotheses of Theorem ABKS. Then as $n \to \infty$, Algorithm G yields a hamilton circuit with probability approaching 1. Similarly, let $D_{m*}$ be a random directed graph that satisfies the hypotheses of Theorem Frieze-ABKS. Then as $n \to \infty$, Algorithm D obtains a hamilton cycle in $D_{m*}$ with probability approaching 1. $V$ and $E = \{e_i \mid i = 1, 2, \ldots, m\}$ define the respective vertices and arcs of a random directed graph, $D_p$, in which each arc is chosen with a fixed probability, p. In [2], Angluin and Valiant describe an $O(n(\log n))$ algorithm, A, such that

$$\lim_{n \to \infty} Pr(A \text{ obtains a hamilton circuit in } D.) = 1 - n^{-\alpha}.$$

Here $p = \frac{c(\log n)}{n}$ and c is dependent on α. With

$$m = \frac{1}{2}(n(\log n) + n \log(\log n) + c(n)),$$

$G_m$ is a graph chosen randomly from among all graphs containing m edges. The latter was the definition of a random graph used by Erdös and Renyi in [9]. Bollobás, Fenner and Frieze used the same definition of a random graph in [5]. In [5], they gave an algorithm for obtaining a hamilton circuit in a random graph in $G_m$ with running time $O(n^{3+o(1)})$. In [20], McDiarmid proved that if $D_m$ is a random directed graph with $m = n(\log n + c)$, then

$$\lim_{n \to \infty} Pr(D_m \text{ is hamiltonian.}) = 1$$

In [13], Frieze gave a sharp threshold algorithm with running time $O(n^{1.5})$ for obtaining a hamilton circuit as $n \to \infty$. In [14], Frieze and Luczak proved that as $n \to \infty$, $R_5$ almost always has a hamilton circuit. In [7], Cooper and Frieze proved that as $n \to \infty$, $D_{3-in,3-out}$ almost always has a hamilton circuit. Example 1.3 follows. The initial *SCORE* values of some permutations were incorrect. The only effect of this was to generate more iterations than necessary.

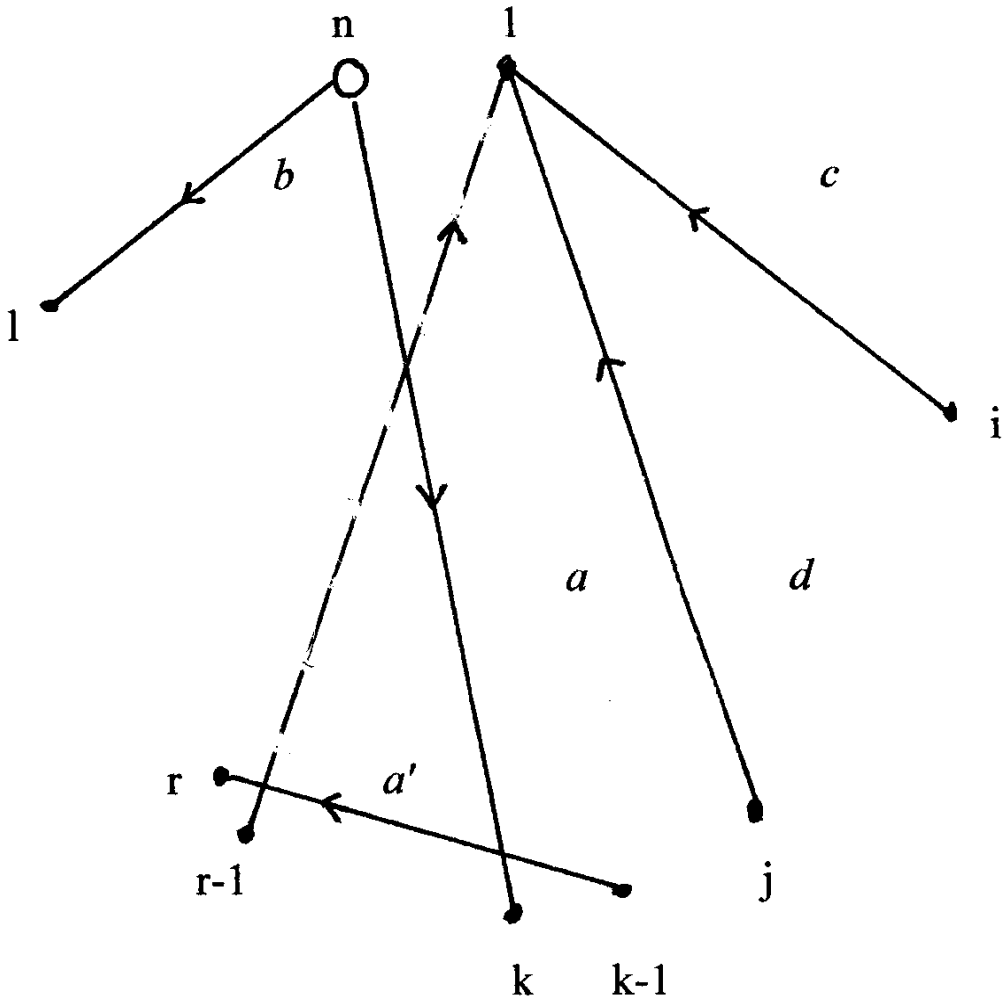

*a* intersects *a'*

Fig. 1.1

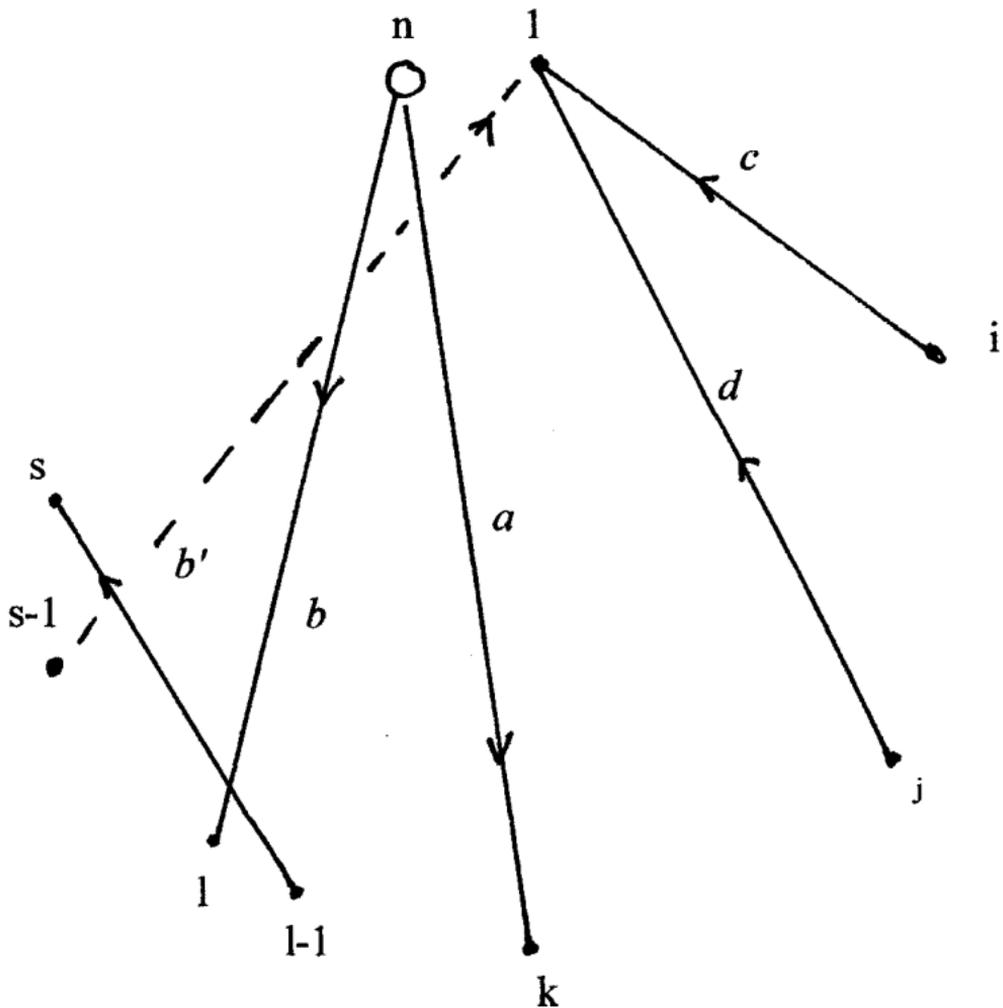

b intersects b'

Fig. 1.2

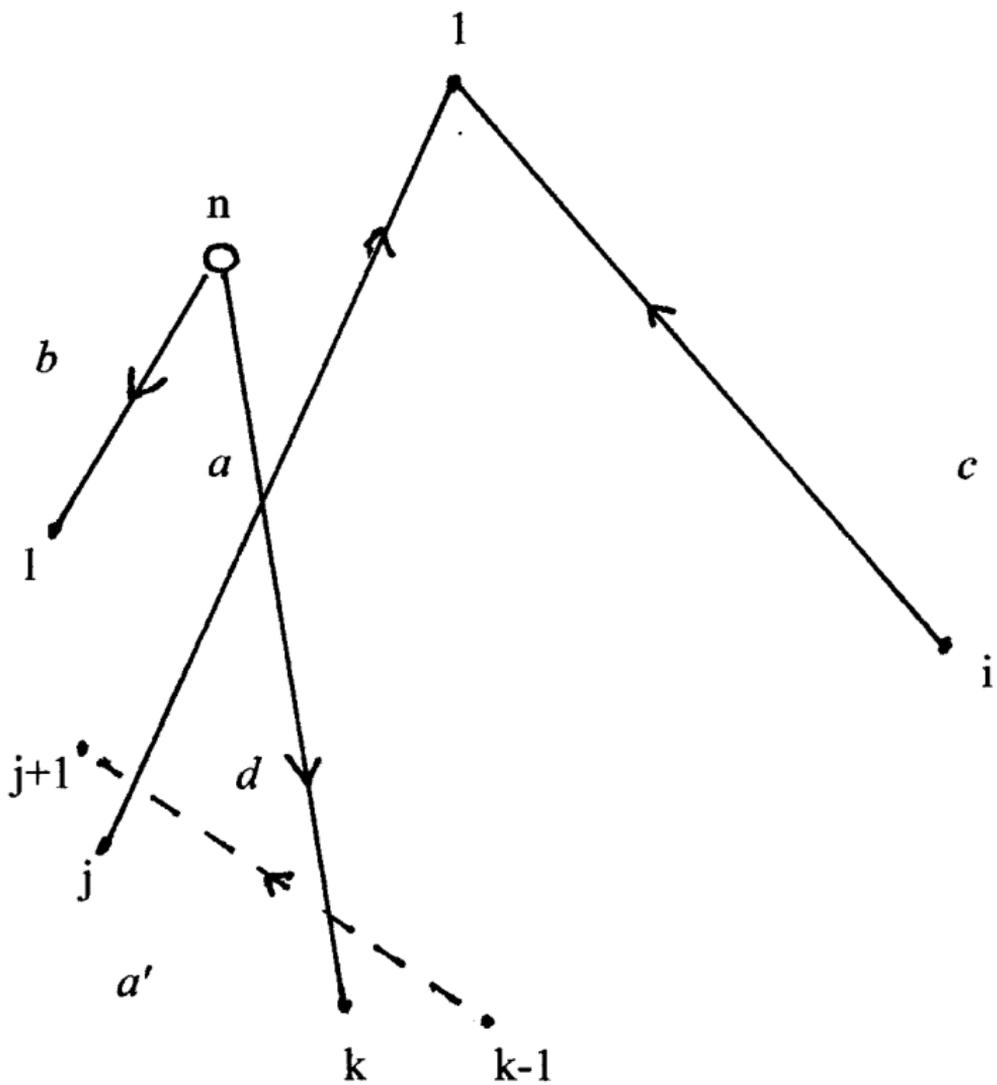

*a* intersects *d*

Fig. 1.4

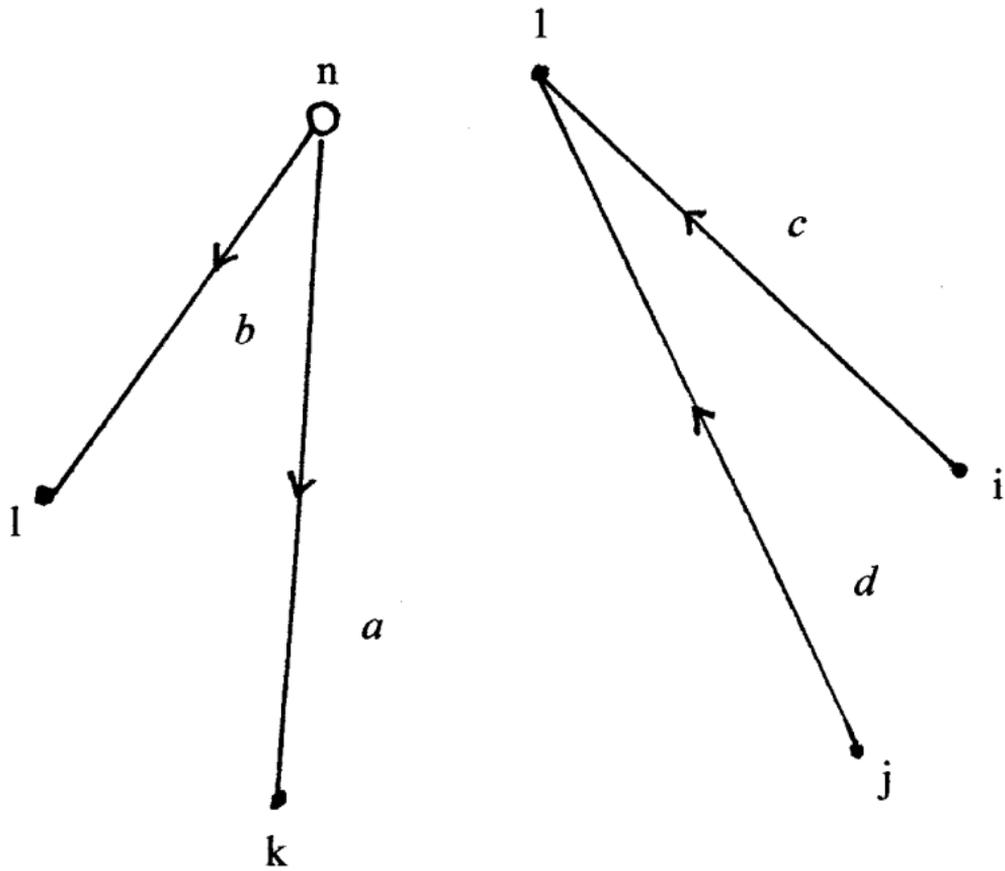

No H-admissible 3-cycles exist

Fig. 1.7

**Example 1.3** The arcs of graph $G'$ are

| | |
|---|---|
| 1: 7, 20, 9-6-4 | 14: 18, 19, 11-12-2 |
| 2-12-11: 16, 14 | 16: 3, 8, 11-12-2 |
| 3: 24, 16, 17, 14, 21-15-10 | 17: 25, 8, 7, 22 |
| 4-6-9: 22, 1 | 18: 5, 14, 8, 13 |
| 5: 18, 8, 13, 16, 24 | 19: 14, 23, 4-6-9 |
| | 20: 23, 1, 24, 13, 2-12-11 |
| 7: 1, 13, 25, 17 | 21-15-10: 23, 25 |
| 8: 16, 17, 5, 18 | 22: 17, 2-12-11, 9-6-4 |
| 9-6-4: 13, 19 | 23: 19, 20, 10-15-21 |
| 10-15-21: 24, 3 | 24: 3, 20, 5, 21-15-10 |
| 11-12-2: 22, 20 | 25: 17, 7, 10-15-21 |
| 13: 7, 4-6-9, 18, 20 | |

$H_0$ = (1  2-12-11  3  4-6-9  5  7  8  10-15-21  13  14  16  17  18  19  20  22  23  24  25)

PSEUDO = {all vertices}

ORD                                  $ORD^{-1}$

$ORD(1) = 1$                         $ORD^{-1}(1) = 1$

$ORD(2) = 2 - 12 - 11$               $ORD^{-1}(2 - 12 - 11) = 2$

$ORD(3) = 3$                         $ORD^{-1}(3) = 3$

$ORD(4) = 4 - 6 - 9$                 $ORD^{-1}(4 - 6 - 9) = 4$

$ORD(5) = 5$                $ORD^{-1}(5) = 5$

$ORD(6) = 7$                $ORD^{-1}(7) = 6$

$ORD(7) = 8$                $ORD^{-1}(8) = 7$

$ORD(8) = 10-15-21$         $ORD^{-1}(10-15-21) = 8$

$ORD(9) = 13$               $ORD^{-1}(13) = 9$

$ORD(10) = 14$              $ORD^{-1}(14) = 10$

$ORD(11) = 16$              $ORD^{-1}(16) = 11$

$ORD(12) = 17$              $ORD^{-1}(17) = 12$

$ORD(13) = 18$              $ORD^{-1}(18) = 13$

$ORD(14) = 19$              $ORD^{-1}(19) = 14$

$ORD(15) = 20$              $ORD^{-1}(20) = 15$

$ORD(16) = 22$              $ORD^{-1}(22) = 16$

$ORD(17) = 23$              $ORD^{-1}(23) = 17$

$ORD(18) = 24$              $ORD^{-1}(24) = 18$

$ORD(19) = 25$              $ORD^{-1}(25) = 19$

Henceforth, ordinary natural numbers represent vertices, while underlined natural numbers represent the ordinal values of vertices in permutations.

$H_0$-admissible 3-cycles: (1) (3 10 14), (2) (3 10 13). $SCORE(i) = 2, (i = 1,2)$.

   (1) (3 10 14): (3 11), (10 14), (13 4): (3 16), (14 19), (18 4-6-9).

   (2) (3 10 13): (3 11), (10 13), (12 4): (3 16), (14 18), (17 4-6-9).

*SCORE* (2)= 2. (2) yields the arcs (3 16), (14 18) and the pseudo-arc (17 4-6-9).

Applying these to $H_0$, we obtain **in ordinal values**

$$H_{01} = ABBREV = (\underline{1 \ \ldots \ \ 3 \ \ 11 \ 12 \ \ 4 \ \ \ldots \ \ 10 \ \ 13 \ \ \ldots})$$

$PSEUDO = \{1, 2\text{-}12\text{-}11, 4\text{-}6\text{-}9, 5, 7, 8, 9, 10\text{-}15\text{-}21, 13, 16, 17, 18, 20, 23, 24,$

$25\}$

The entries in $PSEUDO$ are always the natural number values of vertices.

We obtain the following $H_{01}$-admissible pseudo-3-cycles:

(1) (<u>12  6  8</u>): (<u>12  7</u>), (<u>6  9</u>), (<u>8  4</u>): (17  8), (7  13), (10-15-21  4-6-9).

(2) (<u>12  15</u>)(<u>18  14</u>): (<u>12  16</u>), (<u>15  4</u>), (<u>18  15</u>), (<u>14  19</u>): (17  22), (20  17), (13  20), (19  14).

$SCORE(i) = 2, (i = 1,2)$. Choose (2).

We obtain

$$H_{02} = ABBREV = (\underline{1 \ \ \ldots \ \ 3 \ \ 11 \ \ 12 \ \ 16 \ \ \ldots \ \ 18 \ \ 15 \ \ 4 \ \ \ldots \ \ 10 \ \ 13 \ \ 14 \ \ 19})$$

$PSEUDO = \{1, 2\text{-}12\text{-}11, 4\text{-}6\text{-}9, 10\text{-}15\text{-}21, 22; 13, 23, 25; 7, 8, 18, 20; 5\}.$

Choose 5 from $PSEUDO$.

We obtain the following $H_{02}$-admissible permutations:

(1) (<u>5  6  8</u>): (<u>5  7</u>), (<u>6  9</u>), (<u>8  6</u>): (5  8), (7  13), (10-15-21  7).

(2) (<u>5  6  19</u>): (<u>5  7</u>), (<u>6  1</u>), (<u>19  6</u>): (5  8), (7  1), (25  7).

(3) (<u>5  6  14</u>): (<u>5  7</u>), (<u>6  19</u>), (<u>14  6</u>): (5  8), (7  25), (19  7).

(4) (<u>5  8  17</u>): (<u>5  9</u>), (<u>8  18</u>), (<u>17  6</u>): (5  13), (10-15-21  24), (23  7).

(5) (<u>5  8  2</u>): (<u>5  9</u>), (<u>8  3</u>), (<u>2  6</u>): (5  13), (10-15-21  3), (2-12-11  7).

$SCORE(2) = 3$. Choose (<u>5  6  19</u>): (<u>5  7</u>), (<u>6  1</u>), (<u>19  6</u>).

$$H_{03} = ABBREV = (\underline{1 \ 2 \ 3 \ 11 \ 12 \ 16 \ 17 \ 18 \ 15 \ 4 \ 5 \ 7 \ 8 \ 9 \ 10 \ 13 \ 14 \ 19 \ 6}).$$

$PSEUDO = \{1, 2\text{-}12\text{-}11, 4\text{-}6\text{-}9, 10\text{-}15\text{-}21, 22;\ 13, 19, 23;\ 8, 18, 20\}$

$H_1 = (\underline{1\ 2\ 3\ 4\ 5\ 6\ 7\ 8\ 9\ 10\ 11\ 12\ 13\ 14\ 15\ 16\ 17\ 18\ 19})$

| $ORD$ | $ORD^{-1}$ |
|---|---|
| $ORD(1) = 1$ | $ORD^{-1}(1) = 1$ |
| $ORD(2) = 2-12-11$ | $ORD^{-1}(2-12-11) = 2$ |
| $ORD(3) = 3$ | $ORD^{-1}(3) = 3$ |
| $ORD(4) = 16$ | $ORD^{-1}(4-6-9) = 10$ |
| $ORD(5) = 17$ | $ORD^{-1}(5) = 11$ |
| $ORD(6) = 22$ | $ORD^{-1}(7) = 19$ |
| $ORD(7) = 23$ | $ORD^{-1}(8) = 12$ |
| $ORD(8) = 24$ | $ORD^{-1}(10-15-21) = 13$ |
| $ORD(9) = 20$ | $ORD^{-1}(13) = 14$ |
| $ORD(10) = 4-6-9$ | $ORD^{-1}(14) = 15$ |
| $ORD(11) = 5$ | $ORD^{-1}(16) = 4$ |
| $ORD(12) = 8$ | $ORD^{-1}(17) = 5$ |
| $ORD(13) = 10-15-21$ | $ORD^{-1}(18) = 16$ |
| $ORD(14) = 13$ | $ORD^{-1}(19) = 17$ |
| $ORD(15) = 14$ | $ORD^{-1}(20) = 9$ |
| $ORD(16) = 18$ | $ORD^{-1}(22) = 6$ |
| $ORD(17) = 19$ | $ORD^{-1}(23) = 7$ |
| $ORD(18) = 25$ | $ORD^{-1}(24) = 8$ |
| $ORD(19) = 7$ | $ORD^{-1}(25) = 18$ |

$H_1 = (1\ \ 2\text{-}12\text{-}11\ \ 3\ \ 16\ \ 17\ \ 22\ \ 23\ \ 24\ \ 20\ \ 4\text{-}6\text{-}9\ \ 5\ \ 8\ \ 10\text{-}15\text{-}21\ \ 13\ \ 14\ \ 18\ \ 19\ \ 25\ \ 7)$.

Choose 8.

(1) (<u>12 3 7</u>): (<u>12 4</u>), (<u>3 8</u>), (<u>7 13</u>): (8 16), (3 24), (23 10-15-21).

(2) (<u>12 3 4</u>): (<u>12 4</u>), (<u>3 5</u>), (<u>4 13</u>): (8 16), (3 17), (16 10-15-21).

(3) (<u>12 15 16</u>): (<u>12 16</u>), (<u>15 17</u>), (<u>16 13</u>): (8 18), (14 19), (18 10-15-21).

(4) (<u>12 15 1</u>): (<u>12 16</u>), (<u>15 2</u>), (<u>1 13</u>): (8 18), (14 2-12-11), (1 10-15-21).

(5) (12 3)(1 8): (12 4), (3 13), (1 9), (8 2): (8 16), (3 10-15-21), (1 20), (24 2-12-11).

SCORE(1) = 2 : (12 4), (3 8), (7 13).

$$H_{11} = ABBREV = (\ldots 3\ 8\ \ldots\ 12\ 4\ \ldots\ 7\ 13\ \ldots)$$

$$PSEUDO = \{2\text{-}12\text{-}11,\ 4\text{-}6\text{-}9,\ 10\text{-}15\text{-}21,\ 22;\ 1,\ 13,\ 18,\ 19,\ 20,\ 25;\ 16\}$$

Choose 1.

(1) (1 8 10): (1 9), (8 11), (10 2): (1 20), (24 5), (4-6-9 2-12-11).

(2) (4 2)(16 11): (4 3), (2 15), (16 12), (11 17): (16 3), (2 5), (18 8), (5 19).

(3) (4 2)(16 10): (4 3), (2 5), (16 11), (10 17): (16 3), (2 17), (18 5), (4-6-9 19).

SCORE(3) = 2. (4 2)(16 10): (4 3), (2 0 5), (16 11), (10 17).

$$H_{12} = ABBREV = (1\ 2\ 5\ 6\ 7\ 13\ 14\ 15\ 16\ 11\ 12\ 4\ 3\ 8\ 9\ 10\ 17\ 18\ 19)$$

$$= (1\ 2\text{-}12\text{-}11\ 17\ 22\ 23\ 10\text{-}15\text{-}21\ 13\ 14\ 18\ 5\ 8\ 16\ 3\ 24\ 20\ 4\text{-}6\text{-}9\ 19\ 25\ 7\}$$

$$PSEUDO = \{2\text{-}12\text{-}11,\ 4\text{-}6\text{-}9,\ 10\text{-}15\text{-}21,\ 22;\ 1,\ 13,\ 19,\ 20\}$$

Choose 13.

The following $H_{12}$ permutations have SCORE values of 2:

(1) (14 9 6): (14 10), (9 7), (6 15): (13 4-6-9), (20 23), (22 14).

(2) (14 9 1): (14 10), (9 2), (1 15): (13 4-6-9), (20 2-12-11), (1 14).

(3) (14 18 2): (14 19), (18 5), (2 15): (13 7), (25 17), (2-12-11 14).

Unfortunately, I chose a permutation that had a SCORE value of 1: (14 9 19).

(14 9 19): (14 10), (9 1), (19 15).

$$H_{13} = ABBREV = (1\ 2\ 5\ 6\ 7\ 13\ 14\ 10\ 17\ 18\ 19\ 15\ 16\ 11\ 12\ 4\ 3\ 8\ 9)$$

$$= \{2\text{-}12\text{-}11,\ 4\text{-}6\text{-}9,\ 10\text{-}15\text{-}21,\ 22;\ 1,\ 7,\ 19\}.$$

Choose 1 and 19.

(1) (1 18)(17 19): (1 19), (18 2), (17 15), (19 18): (1 7), (25 2-12-11), (19 14), (7 25).

$SCORE(1) = 2$.

$$H_{14} = ABBREV = (1\ 19\ 18\ 2\ 5\ 6\ 7\ 13\ 14\ 10\ 17\ 15\ 16\ 11\ 12\ 4\ 3\ 8\ 9)$$

$$= (1\ 7\ 25\ 2\text{-}12\text{-}11\ 17\ 22\ 23\ 10\text{-}15\text{-}11\ 13\ 4\text{-}6\text{-}9\ 19\ 14\ 18\ 5\ 8\ 16\ 3\ 24\ 20)$$

$$PSEUDO = \{2\text{-}12\text{-}11,\ 4\text{-}6\text{-}9,\ 10\text{-}15\text{-}11,\ 22;\ 25\}$$

Choose 25.

The following permutations have a *SCORE* value of 1:

(1) (18 2 12): (18 5), (2 4), (12 2): (25 17), (2-12-11 2), (8 2-12-11).

(2) (18 2 17): (18 5), (2 15), (17 2): (25 17), (2-12-11 14), (19 2-12-11).

(3) (18 2 9): (18 5), (2 1), (9 2): (25 17), (2-12-11 1), (20 2-12-11).

(4) (18 2 6): (18 5), (2 7), (6 2): (25 17), (2-12-11 23), (22 2-12-11).

Choose (1) since the degree of its pseudo-arc vertex, 8, is 4.

(18 2 12): (18 5), (2 4), (12 2) yields

$$H_2 = H_{04} = ABBREV = (1\ 19\ 18\ 5\ 6\ 7\ 13\ 14\ 10\ 17\ 15\ 16\ 11\ 12\ 2\ 4\ 3\ 8\ 9)$$

$$= (1\ 7\ 25\ 17\ 22\ 23\ 10\text{-}15\text{-}21\ 13\ 4\text{-}6\text{-}9\ 19\ 14\ 18\ 5\ 8\ 2\text{-}12\text{-}11\ 16\ 3\ 24\ 20)$$

$$PSEUDO = \{4\text{-}6\text{-}9,\ 10\text{-}15\text{-}11,\ 22;\ 8\}$$

| *ORD* | *ORD⁻¹* |
|---|---|
| ORD(1) = 1 | ORD⁻¹(1) = 1 |
| ORD(2) = 7 | ORD⁻¹(2 – 12 – 11) = 15 |
| ORD(3) = 25 | ORD⁻¹(3) = 17 |
| ORD(4) = 17 | ORD⁻¹(4 – 6 – 9) = 9 |
| ORD(5) = 22 | ORD⁻¹(5) = 13 |
| ORD(6) = 23 | ORD⁻¹(7) = 2 |
| ORD(7) = 10 – 15 – 21 | ORD⁻¹(8) = 14 |
| ORD(8) = 13 | ORD⁻¹(10 – 15 – 21) = 7 |
| ORD(9) = 4 – 6 – 9 | ORD⁻¹(13) = 8 |
| ORD(10) = 19 | ORD⁻¹(14) = 11 |
| ORD(11) = 14 | ORD⁻¹(16) = 16 |
| ORD(12) = 18 | ORD⁻¹(17) = 4 |
| ORD(13) = 5 | ORD⁻¹(18) = 12 |
| ORD(14) = 8 | ORD⁻¹(19) = 10 |
| ORD(15) = 2 – 12 – 11 | ORD⁻¹(20) = 19 |
| ORD(16) = 16 | ORD⁻¹(22) = 5 |
| ORD(17) = 3 | ORD⁻¹(23) = 6 |
| ORD(18) = 24 | ORD⁻¹(24) = 18 |
| ORD(19) = 20 | ORD⁻¹(25) = 3 |

Choose 8.

(1) (14 3 6): (14 4), (3 7), (6 15): (8 17), (25 10-15-21), (23 2-12-11).

(2) (14 15 10): (14 16), (15 1), (10 15): (8 16), (2-12-11 14), (19 2-12-11).

(3) (7 17)(9 4): (7 18), (17 8), (9 5), (4 10): (10-15-21 24), (3 13), (4-6-9 22), (17 19).

(4) (7 17)(9 19): (7 18), (17 8), (9 1), (19 10): (10-15-21 24), (3 13), (4-6-9 1).

(5) (7 16)(9 4): (7 17), (16 8), (9 5), (4 10): (10-15-21 3), (16 13), 4-6-9 22), (17 19).

(6) (7 16)(9 19): (7 17), (16 8), (9 1), (19 10): (10-15-21 3), (16 3), (4-6-9 1), (20 19).

(7) (<u>14 3 5</u>): (<u>14 4</u>), (<u>3 6</u>), (<u>5 15</u>): (8 17), (25 23), (22 2-12-11).

(8) (<u>14 15 19</u>): (<u>14 16</u>), (<u>15 1</u>), (<u>19 15</u>): (8 16), (2-12-11 1), (20 2-12-11).

(9) (<u>14 15 5</u>): (<u>14 16</u>), (<u>15 6</u>), (<u>5 15</u>): (8 16), (2-12-11 23), (22 2-12-11).

$SCORE(i) = 0$, $(i = 1,2,...,9)$.

$$H_{21} = ABBREV = (1 \ldots 4\ 10 \ldots 17\ 8\ 9\ 5 \ldots 7\ 18 \ldots)$$

Choose (3). It contains a pseudo-arc vertex, 3, that is the initial vertex of the arc

(3 10-15-21). The latter arc can be used to construct a rotation.

(<u>17 7</u>)(<u>9 4</u>): (<u>7 18</u>), (<u>17 8</u>), (<u>9 5</u>), (4 10) yields

$$H_{22} = ABBREV = (1 \ldots 4\ 10 \ldots 17\ 8\ 9\ 5 \ldots 7\ 18\ 19)$$

Using the rotation defined by (3 10-15-21) = (<u>17 7</u>), we obtain

$$H_{23} = ABBREV = (1 \ldots 4\ 10 \ldots 17\ 7\ \text{- - -}\ 5\ 9\ 8\ 18\ 19).$$

= (1 7 25 17 19 14 18 5 8 2-12-11 16 3 21-15-10 23 22 9-6-4 13 24 20)

We give the usual representation of vertices to show the changes in orientation of 10-15-21 and 4-6-9.

$$PSEUDO = \{13, 23;\ 8, 17\}$$

Choose 17.

(1) (<u>4 13 9</u>): (<u>4 14</u>), (<u>13 8</u>), (<u>9 10</u>): (17 8), (5 13), (9-4 19).

(2) (<u>4 13 15</u>): (<u>4 14</u>), (<u>13 16</u>), (<u>15 10</u>): (17 8), (5 3), (2-12-11 19).

(3) (<u>4 2 3</u>): (<u>4 3</u>), (<u>2 4</u>), (<u>3 10</u>): (17 25), (7 17), (25 19).

(4) (<u>4 13 8</u>): (<u>4 14</u>), (<u>13 18</u>), (<u>8 10</u>): (17 8), (5 24), (13 9).

(5) (<u>4 6 18</u>): (<u>4 5</u>), (<u>6 19</u>), (<u>18 10</u>): (17 13), (23 20), (24 19).

(6) (<u>4 6</u>)(<u>14 3</u>): (<u>4 5</u>), (<u>6 10</u>), (<u>14 4</u>), (<u>3 15</u>): (17 22), (23 19), (8 17), (25 2-12-11).

*SCORE(6) = 2*. Choose (6). (<u>4 6</u>)(<u>14 3</u>): (<u>4 5</u>), (<u>6 10</u>), (<u>14 4</u>), (<u>3 15</u>).

$$H_{24} = (\underline{1 \ldots 3 \ 15 \ldots 17 \ 7 \ 6 \ 10 \ldots 14 \ 4 \ 5 \ 9 \ 8 \ 18 \ 19})$$

$$PSEUDO = \{13, 25\}$$

The only $H_{24}$-admissible permutation is (<u>8 11 3</u>). Its *SCORE* value is 0.

(<u>8 11 13</u>): (<u>8 12</u>), (<u>11 14</u>), (<u>13 18</u>): (13 18), (14 8), (5 24).

$$H_{25} = (\underline{1 \ldots 3 \ 15 \ldots 17 \ 7 \ 6 \ 10 \ 11 \ 14 \ 4 \ 5 \ 9 \ 8 \ 12 \ 12 \ 18 \ 19})$$

$$PSEUDO = \{14, 25\}.$$

Neither 14 nor 30 contained an arc whose terminal vertex was an *r*-vertex on $H_{24}$.

Choose 14 and search for $H_{24}$- admissible permutations.

(1) (<u>11 8 1</u>): (<u>11 12</u>), (<u>8 2</u>), (<u>1 14</u>): (14 18), (13 7), (1 8).

(2) (<u>11 8</u>)(<u>3 14</u>): (<u>11 12</u>), (<u>8 14</u>), (<u>3 4</u>), (<u>14 15</u>): (14 18), (13 8), (8 2-12-11), (25 17).

(3) (<u>11 8 4</u>): (<u>11 12</u>), (<u>8 5</u>), (<u>4 14</u>): (14 18), (13 22), (17 8).

(4) (<u>11 8 16</u>): (<u>11 12</u>), (<u>8 17</u>), (<u>16 14</u>): (14 18), (13 3), (16 8).

*SCORE(i) = 0 (i = 1,2,3,4)*. Since (8 2-12-11) is a pseudo-arc obtained from (2), we choose (2).

$H_{26}$ = ABBREV = (<u>1 … 6 9 8 14 … 17 7 6 10 … 13 18 19</u>)

= (1 7 25 17 22 9-6-4 13 8 2-12-11 16 3 21-15-10 23 19 14 18 5 24 20)

We now obtain a rotation using the arc (8 17) that satisfies theorem 1.7.

$H_{27}$ = ABBREV = (<u>14 4 --- 18 13 --- 10 6 7 17 --- 15 5 9 8</u>)

= (<u>1 --- 18 13 --- 10 6 7 17 --- 15 5 9 8 14 4 --- 2</u>)

= (1 20 24 5 18 14 19 23 10-15-21 3 16 11-12-2 22 9-6-4 13 8 17 25 7)

$$PSEUDO = \{13\}).$$

We now obtain

(1) (8̲ 1̲3̲ 1̲6̲): (8̲ 1̲2̲), (1̲3̲ 1̲5̲), (1̲6̲ 1̲4̲): (13 18), (5 11-12-2), (16 8).

This yields $H_{28}$ = ABBREV = (1 ---- 18 13 --- 10 6 7 17 --- 15 5 9 8 14 4 --- 2)

= (1 20 24 5 11-12-2 22 22 9-6-4 13 18 14 19 23 10-15-21 3 16 8 17 25 7),

a hamilton circuit.

We next give an example of a heuristic variation of Algorithm G, Algorithm $G_{heuristic}$. It differs from Algorithm G only in one precedure: If all of the SCORE values of an iteration are 0, or if we fail to obtain an $H_i$-admissible permutation, we use theorem 1.7 to try to obtain a rotation with a positive SCORE value. Given a pseudo-arc vertex, $a$, of greatest degree, we test all arcs $(a\ b_j)$ out of $a$ to obtain one defining a rotation with a positive SCORE value. If we are unable to obtain one, we choose an arc $(a\ b_j)$ where deg($H_i(b_j)$) is greatest. If we haven't obtained a hamilton circuit, we continue the algorithm using the procedures of Algorithm G. The greatest value of this algorithm is its simplicity. Given an arc $(a\ b_j)$ it is simple to check whether or not $(H_i(a)\ H_i(b))$ is an arc of $G'$.

**Example 1.4** For simplicity, we start with $H_{27}$ of example 1.1. Our pseudo-arc vertex is 13.

$H_{27}$ = (1 20 24 5 18 14 19 23 10-15-21 3 16 11-12-2 22 9-6-4 **13** 8 17 25 7).

None of the arcs out of 13 satisfy theorem 1.7. We therefore, choose the arc (13 20) because deg(24) has the greatest degree among all possibilities for $H_i(b_j)$. We obtain

$H_{28}$ = (13 20 1 7 25 17 **8** 24 5 18 14 19 23 10-15-21 3 16 11-12-2 22 9-6-4)

None of the arcs chosen satisfy theorem 1.7. We choose (8 18).

$H_{29}$ = (8 18 5 **24** 14 19 23 10-15-21 3 16 11-12-2 22 9-6-4 13 20 1 7 25 17).

None of the arcs chosen satisfy theorem 1.7. We choose (24 3).

$H_{2\,10}$ = (24 3 21-15-10 23 19 **14** 16 11-12-2 22 9-6-4 13 20 1 7 25 17 8 18 5).

If we choose (14 18), (16 5) is an arc of $G'$. Thus, we obtain the hamilton circuit

$H_{2\,11}$ = (14 18 8 17 25 7 1 20 13 4-6-9 22 2-12-11 16 5 24 3 21-15-10 23 19)

Assume $G$ is a graph containing a hamilton circuit. Algorithm $G_{no\ r-vertices}$ obtains a hamilton circuit directly from $G$ without the use of $r$-vertices. If we fail to obtain an $H_i$-admissible permutation using a pseudo-arc vertex of degree 2, $a$, and $(a\ b)$ is the unique arc not on $H_i$ emanating from $a$, we construct the rotation determined by the arc $(b\ a)$. Thus, the probability of obtaining an *edge* on $H_{i+1}$ containing $a$ is always 1. This is reasonable, since the edge $[a\ b]$ lies on every hamilton circuit in $G$. In all other cases, if we fail to obtain an $H_i$-admissible permutation, we backtrack to $H_{i-1}$. Other than its treatment of vertices of degree 2, it is identical to Algorithm D.

**Example 1.5** $G$ is defined by the following set of arcs:

1 : 7, 20, 4, 9
2 : 22, 20, 12
3 : 24, 16, 10, 17, 14, 21
4 : 13, 19, 6, 1
5 : 18, 8, 13, 16, 24
6 : 4, 9
7 : 13, 1, 25, 17
8 : 5, 17, 16, 18
9 : 22, 1, 6
10 : 25, 23, 15, 3
11 : 16, 14, 12
12 : 11, 2
13 : 4, 7, 18. 5
14 : 19, 18, 11, 3
15 : 10, 21
16 : 3, 11, 8, 5
17 : 8, 25, 7, 22, 3
18 : 14, 5, 8, 13
19 : 23, 14, 4
20 : 1, 23, 2, 24
21 : 15, 3, 24
22 : 2, 9, 17
23 : 20, 19, 10
24 : 21, 3, 20, 5
25 : 17, 10, 7

$H_0 = (1\ 2\ 3\ 4\ 5\ 6\ 7\ 8\ 9\ 10\ 11\ 12\ 13\ 14\ 15\ 16\ 17\ 18\ 19\ 20\ 21\ 22\ 23\ 24\ 25)$.

For $i = 1, 2, \ldots, 25$,

$ORD(i) = i,\qquad ORD^{-1}(i) = i$

$PSEUDO\ =\ \{\text{all vertices except 11}\}$.

Since the deg(3) = 6, choose 3. The ordinal values of all vertices are the same as those of the vertices.

(1) (3 9 21): (3 10), (9 22), (21 4).

(2) (3 9 25): (3 10), (9 1), (25 4).

(3) (3 15 20): (3 16), (15 21), (20 4).

(4) (3 13)(5 17): (3 14), (13 4), (5 18), (17 6).

(5) (3 13)(5 15): (3 14), (13 4), (5 16), (15 6).

$SCORE(4) = SCORE(5) = 3$. Deg(17) = 5, deg(15) = 2. If deg(15) were greater than 2, we would choose (4). Since deg(15) = 2, choose (5).

(3 13)(5 15): (3 14), (13 4), (5 16), (15 6)..

$$H_{01} = ABBREV = (1 \ldots 3\ 14\ 15\ 6 \ldots 13\ 4 \ldots 5\ 16 \ldots)$$

$PSEUDO$ = {1, 2, 4 6, 7, 8, 9, 10, 12, 14, 15, 16, 17, 18, 19, 20, 21, 22, 23, 24, 25}.

Choose 15.

(1) (15 20 25): (15 21), (20 1), (25 16).

(2) (15 20 22): (15 21), (20 23), (22 16).

(3) (15 20 1): (15 21), (20 2), (1 6).

(4) (15 9)(7 12): (15 10), (9 6), (7 13), (12 8).

(5) (15 20)(7 16): (15 21), (20 6), (7 17), (16 8).

$SCORE(4) = SCORE(5) = 3$. deg (12) = 2, deg(16) = 4. Choose (4).

(15 9)(7 12): (15 10), (9 6), (7 13), (12 8).

$$H_{02} = ABBREV = (1 \ldots 3\ 14\ 15\ 10 \ldots 12\ 8\ 9\ 6\ 7\ 13\ 4\ 5\ 16 \ldots)$$

$PSEUDO$ = {6, 12; 2, 19, 21, 22, 23, 25; 1, 4, 8, 10, 14, 16, 18, 20, 24; 17}

(1) (2 17 24): (2 18), (17 25), (24 3).

(2) (4 17 24): (4 18), (17 25), (24 5).

$SCORE(i) = 2$ $(i = 1,2)$. Choose (2).

(4 17 24): (4 18), (17 25), (24 5).

$H_{03}$ = ABBREV = (1 … 3 14 15 10 … 12 8 9 6 7 13 4 18 … 24 5 16 17 25).

PSEUDO = {12, 6; 2, 9, 19, 21, 22, 23, 25; 1, 14, 10, 8, 4, 18, 20, 16}

(1) (4 18 3): (4 19), (18 4), (3 18).

(2) (4 18 24): (4 19), (18 5), (24 18).

(3) (4 18 12): (4 19), (18 8), (12 18)

(4) (4 18 7): (4 19), (18 13), (7 18).

(5) (4 25 6): (4 1), (25 7), (6 18).

(6) (4 25 15): (4 1), (25 10), (15 18).

(7) (16 10)(18 3): (16 11), (10 17), (18 14), (3 19).

(8) (16 12)(18 3): (16 8), (12 17), (18 14), (3 19).

$SCORE(3) = SCORE(5) = 2$. All other $SCORE$ values are 1. $\deg(6) = \deg(12) = 2$. Choose (3).

(4 18 12): (4 19), (18 8), (12 18).

$H_{04}$ = ABBREV = {1 … 3 14 15 10 … 12 18 8 9 6 7 13 4 19 … 24 5 16 17 25)

PSEUDO = {6, 12; 2, 9, 19, 21, 22, 23, 25; 1, 14, 10, 8, 20, 16}

(1) (12 1)(24 2): (12 2), (1 8), (24 3), (2 25).

$SCORE(1) = 2$.

$H_{05}$ = ABBREV = (1 8 9 6 7 13 4 5 16 … 24 3 14 15 10 … 12 2 25)

PSEUDO = {6; 2, 19, 21, 22, 23, 25; 1, 4, 8, 10, 14, 16, 18, 20; 17}

Choose 17.

(1) (17 1 6): (17 8), (1 7), (6 18).

(2) (17 1 8): (17 8), (1 9), (8 18).

(3) (17 1 13): (17 8), (1 4), (13 18).

(4) (17 6 13): (17 7), (6 4), (13 18).

(4) (17 21 14): (17 22), (21 15), (14 18).

(5) (17 24 4): (17 3), (24 5), (4 18).

(6) (17 21)(4 18): (17 22), (21 18), (4 19), (18 5).

(7) (17 1)(4 25): (17 8), (1 18), (4 1), (25 5).

(8) (17 6)(4 18): (17 7), (6 18), (4 19), (18 5).

(10) (17 6)(4 9): (17 7), (6 18), (4 6), (9 5)

(11) (17 6)(4 25): (17 7), (6 18), (4 1), (25 5).

$SCORE(3) = SCORE(5) = SCORE(7) = SCORE(9) = 3.$. We choose (9) because 6 is a pseudo-arc vertex of degree 2.

(17 6)(4 18): (17 7), (6 18), (4 19), (18 5).

$H_{06}$ = ABBREV = (1 8 9 6 18 5 16 17 7 13 4 19 20 21 22 23 24 3 14 15 10 11 12 2 5).

$PSEUDO$ = {6; 2, 19, 21, 22, 23, 25; 1, 8, 10, 14, 16, 20}

(1) (6 13)(16 10): (6 4), (13 18), (16 11), (10 17).

(2) (6 13)(16 1): (6 4), (13 18), (16 11), (10 8).

$SCORE(1) = SCORE(2) = 2$. Choose (1).

$H_1 = H_{07}$ = (1 8 9 6 4 19 20 21 22 23 24 3 14 15 10 17 7 13 18 5 16 11 11 12 2 25)

We now give a new ordering to $H_1$.

$H_1 = (1\ 2\ 3\ 4\ 5\ 6\ 7\ 8\ 9\ 10\ 11\ 12\ 13\ 14\ 15\ 16\ 17\ 18\ 19\ 20\ 21\ 22\ 23\ 24\ 25)$

$PSEUDO = \{2,\ 9,\ 19,\ 21,\ 23,\ 25;\ 1,\ 8,\ 10,\ 14,\ 20\}$

| ORD | $ORD^{-1}$ |
|---|---|
| $ORD(1) = 1$ | $ORD^{-1}(1) = 1$ |
| $ORD(2) = 8$ | $ORD^{-1}(2) = 24$ |
| $ORD(3) = 9$ | $ORD^{-1}(3) = 12$ |
| $ORD(4) = 6$ | $ORD^{-1}(4) = 5$ |
| $ORD(5) = 4$ | $ORD^{-1}(5) = 20$ |
| $ORD(6) = 19$ | $ORD^{-1}(6) = 4$ |
| $ORD(7) = 20$ | $ORD^{-1}(7) = 17$ |
| $ORD(8) = 21$ | $ORD^{-1}(8) = 2$ |
| $ORD(9) = 22$ | $ORD^{-1}(9) = 3$ |
| $ORD(10) = 23$ | $ORD^{-1}(10) = 15$ |
| $ORD(11) = 24$ | $ORD^{-1}(11) = 22$ |
| $ORD(12) = 3$ | $ORD^{-1}(12) = 23$ |
| $ORD(13) = 14$ | $ORD^{-1}(13) = 18$ |
| $ORD(14) = 15$ | $ORD^{-1}(14) = 13$ |
| $ORD(15) = 10$ | $ORD^{-1}(15) = 14$ |
| $ORD(16) = 17$ | $ORD^{-1}(16) = 21$ |
| $ORD(17) = 7$ | $ORD^{-1}(17) = 16$ |
| $ORD(18) = 13$ | $ORD^{-1}(18) = 13$ |
| $ORD(19) = 18$ | $ORD^{-1}(19) = 6$ |
| $ORD(20) = 5$ | $ORD^{-1}(20) = 7$ |
| $ORD(21) = 16$ | $ORD^{-1}(21) = 8$ |
| $ORD(22) = 11$ | $ORD^{-1}(22) = 9$ |
| $ORD(23) = 12$ | $ORD^{-1}(23) = 10$ |
| $ORD(24) = 2$ | $ORD^{-1}(24) = 11$ |
| $ORD(25) = 25$ | $ORD^{-1}(25) = 25$ |

Choose 1.

  (1) (1 16)(2 19): (1 17), (16 2), (2 20), (19 3).

  (2) (1 16)(2 18): (1 17), (16 2), (2 19), (18 3).

  (3) (1 16)(2 20): (1 17), (16 2), (2 21), (20 3).

  (4) (1 6)(2 19): (1 7), (6 2), (2 20), (19 3).

  (5) (1 6)(2 20): (1 7), (6 2), (2 21), (20 3).

  (6) (1 6)(2 15): (1 7), (6 2), (2 16), (15 3): (1 20), (19 8), (8 17), (10 9).

  (7) (1 6)(2 18): (1 7), (6 2), (2 19), (18 3).

  (8) (1 4)(2 19): (1 5), (4 2), (2 20), (19 3).

  (9) (1 4)(2 20): (1 5), (4 2), (2 21), (20 3).

  (10) (1 4)(2 15): (1 5), (4 2), (2 16), (15 3).

  (11) (1 4)(2 18): (1 5), (4 2), (2 19), (18 3).

  (12) (1 16 24); (1 17), (16 25), (24 2).

  (13) (1 6 9): (1 7), (6 10), (9 2): (1 20), (19 23), (22 8).

  (14)  (1 4 2): (1 5), (4 22), (21 2)

  (15) (1 2 19): (1 3), (2 20), (19 2).

  (16) (1 2 15): (1 3), (2 16), (15 2): (1 9), (8 17), (10 8).

  (17) (1 2 18): (1 3), (2 19), (18 2).

  (18) (1 2 5): (1 3), (2 6), (5 2).

$SCORE(16) = SCORE(13) = SCORE(17) = 2$. Since deg(10) = 4, choose

(16).

$$H_{11} = ABBREV = (1\ 3\ \ldots\ 8\ 2\ 16\ \ldots\ 24\ 9\ \ldots\ 15\ 25)$$

$$PSEUDO = \{19,\ 21,\ 22,\ 23,\ 25;\ 10,\ 14,\ 20\}$$

Choose 20.

   (1) (7 9 1): (7 10), (9 3), (1 8).

   (2) (7 10 6): (7 11), (10 7), (6 8): (20 24), (23 20), (19 21).

   (3) (7 10 5): (7 11), (10 6), (5 8).

   (4) (7 23)(13 18): (7 24), (23 8), (13 19), (18 14).

   (5) (7 9)(13 18): (7 10), (9 8), (13 19), (18 14).

   (6) (7 10)(13 18): (7 11), (10 8), (13 19), (18 14).

   (7) (7 25)(13 5): (7 1), (25 8), (13 6), (5 14).

SCORE(2) = 2. (7 10 6): (7 11), (10 7), (6 8).

$$H_{12} = ABBREV = (1\ 3\ \ldots\ 6\ 8\ 2\ 16\ \ldots\ 24\ 9\ 10\ 7\ 11\ \ldots\ 15\ 25)$$

$$PSEUDO = \{19, 21, 22, 25;\ 14\}$$

Choose 14.

   (1) (13 5 17): (13 6), (5 18), (17 14).

   (2) (13 18 19): (13 19), (18 20), (19 14): (14 18), (13 5), (18 15).

   (3) (13 21)(6 9): (13 22), (21 14), (6 10), (9 8): (14 18), (13 15), (19 23), (22 21).

   (4) (13 18)(6 12): (13 19), (18 14), (6 13), (12 8): (14 18), (13 15), (19 23), (3 21).

   (5) (13 21)(6 12): (13 22), (21 14), (6 13), (12 8): (14 18), (16 15), (19 23), (3 21).

SCORE(3) = SCORE(4) = SCORE(5) = SCORE(6) = 1. Choose (4).

$$H_{13} = ABBREV = (1\ 3\ 4\ \ldots\ 6\ 10\ 7\ 11\ \ldots\ 13\ 22\ \ldots\ 24\ 9\ 8\ 2\ 16\ \ldots\ 21\ 14\ 15\ 25)$$

$$PSEUDO\ =\ \{21,\ 22;\ 16,\ 25\}.$$

Choose 16.

(1) (21 11)(9 1): (21 12), (11 14), (9 3), (1 8)

(2) (21 11)(9 2): (21 12), (11 14), (9 16), (2 8).

(3) (21 11 9): (21 12), (11 8), (9 14): (16 3), (24 21), (22 15).

(4) (21 8)(9 2): (21 2),(8 14), (9 16), (2 8): (16 8), (21 15), (22 17), (8 21).

$SCORE(4) = 2.$ Choose (4).

(21 8)(9 2): (21 2), (8 14), (9 16), (2 8).

$H_2\ =\ H_{14}\ =\ $(1 3 4 … 6 10 7 11 … 13 22 … 24 9 16 … 21 2 8 14 15 25)

= (1 9 6 4 19 23 20 24 3 14 11 12 2 22 17 7 13 18 5 16 8 21 15 10 25).

$$PSEUDO\ =\ \{8,\ 25\}.$$

| | |
|---|---|
| $ORD(1) = 1$ | $ORD^{-1}(1) = 1$ |
| $ORD(2) = 9$ | $ORD^{-1}(2) = 13$ |
| $ORD(3) = 6$ | $ORD^{-1}(3) = 9$ |
| $ORD(4) = 4$ | $ORD^{-1}(4) = 4$ |
| $ORD(5) = 19$ | $ORD^{-1}(5) = 19$ |
| $ORD(6) = 23$ | $ORD^{-1}(6) = 3$ |
| $ORD(7) = 20$ | $ORD^{-1}(7) = 16$ |
| $ORD(8) = 24$ | $ORD^{-1}(8) = 21$ |
| $ORD(9) = 3$ | $ORD^{-1}(9) = 2$ |
| $ORD(10) = 14$ | $ORD^{-1}(10) = 24$ |
| $ORD(11) = 11$ | $ORD^{-1}(11) = 11$ |
| $ORD(12) = 12$ | $ORD^{-1}(12) = 12$ |
| $ORD(13) = 2$ | $ORD^{-1}(13) = 17$ |
| $ORD(14) = 22$ | $ORD^{-1}(14) = 10$ |
| $ORD(15) = 17$ | $ORD^{-1}(15) = 23$ |
| $ORD(16) = 7$ | $ORD^{-1}(16) = 20$ |
| $ORD(17) = 13$ | $ORD^{-1}(17) = 15$ |
| $ORD(18) = 18$ | $ORD^{-1}(18) = 18$ |
| $ORD(19) = 5$ | $ORD^{-1}(19) = 5$ |
| $ORD(20) = 16$ | $ORD^{-1}(20) = 7$ |
| $ORD(21) = 8$ | $ORD^{-1}(21) = 22$ |
| $ORD(22) = 21$ | $ORD^{-1}(22) = 14$ |
| $ORD(23) = 15$ | $ORD^{-1}(23) = 6$ |
| $ORD(24) = 10$ | $ORD^{-1}(24) = 8$ |
| $ORD(25) = 25$ | $ORD^{-1}(25) = 25$ |

Choose 25 for 3-cycles and 25 and 8 for POTDTC's.

(1) (<u>25 15 20</u>): (<u>25 16</u>), (<u>15 21</u>), (<u>20 1</u>): (25 7), (17 8), (16 1).

(2) (<u>25 15 24</u>): (<u>25 16</u>), (<u>15 25</u>), (<u>24 1</u>).

(3) (<u>25 15</u>)(<u>21 14</u>): (<u>25 16</u>), (<u>15 1</u>), (<u>21 15</u>), (<u>14 22</u>).

$SCORE(i) = 0$ $(i = 1,2,3)$.

No pseudo-arc vertex is of degree 2.
Choose (1)

(25 15 20): (25 16), (15 21), (20 1).

$$H_{21} = ABBREV = (1 \ldots 15\ 21 \ldots 25\ 16 \ldots 20)$$

$$PSEUDO = \{16,\ 8\}$$

Choose 8.

(1) (21 18 9): (21 19), (18 10), (9 22): (8 5), (18 14), (3 21)

$SCORE(1) = 1$.

$$H_{22} = ABBREV = (1 \ldots 9\ 22 \ldots 25\ 16 \ldots 18\ 10 \ldots 15\ 21\ 19\ 20)$$

$$PSEUDO = \{16\}$$

(1) (20 8 9): (20 9), (8 22), (9 1): (16 3), (24 21), (3 1).

(2) (20 8 21): (20 9), (8 19), (21 1): (16 3), (24 5), (8 1).

$SCORE(i) = 0\ (i = 1,2)$.

Neither (1) nor (2) has a pseudo-arc vertex of degree 2. deg(3) = 6. Thus, choose (1).

$$H_{23} = ABBREV = (1 \ldots 8\ 22 \ldots 25\ 16 \ldots 18\ 10\ 11 \ldots 15\ 21\ 19\ 20\ 9).$$

$$PSEUDO = \{3\}$$

(1) (9 18 15): (9 10), (18 21), (15 1): (3 14), (18 8), (17 1).

$SCORE(1) = 0$.

$$H_{24} = ABBREV = (1 \ldots 8\ 22 \ldots 25\ 16 \ldots 18\ 21\ 19\ 20\ 9 \ldots 15)$$

$$PSEUDO = \{17\}$$

(1) (15 24 20): (15 25), (24 9), (20 1).

(2) (15 25 14): (15 16), (25 15), (14 1).

(3) (15 20 10): (15 9), (20 11), (10 1): (17 3), (16 11), (14 1).

$SCORE(i) = 0$ $(i = 1,2,3)$.

$\deg(14) = \deg(16) = 4$. Choose (3).

$H_{25} = ABBREV = (1 \ldots 8\ 22 \ldots 25\ 16 \ldots 18\ 21\ 19\ 20\ 11 \ldots 15\ 9\ 10)$.

$PSEUDO = \{14\}$.

(10 4 16): (10 5), (4 17), (16 1): (14 19), (4 13), (7 1).

There are no pseudo-arc vertices in this permutation. Thus,

$H_3 = H_{25} = (1 \ldots 4\ 17\ 18\ 21\ 19\ 20\ 11 \ldots 15\ 9\ 10\ 5 \ldots 8\ 22 \ldots 25\ 16)$

$= (1\ 9\ 6\ 4\ 13\ 18\ 8\ 5\ 16\ 11\ 12\ 2\ 22\ 17\ 3\ 14\ 19\ 23\ 20\ 24\ 21\ 15\ 10\ 25\ 7)$

is a hamilton circuit.